\documentclass[12pt]{article}
\usepackage{amsmath,amssymb,amscd,enumerate,eucal,bm}
\numberwithin{equation}{section}
\newcommand{\R}{\mathbb{R}}
\newcommand{\Z}{\mathbb{Z}}
\newcommand{\N}{\mathbb{N}}
\newcommand{\C}{\mathbb{C}}

\newcommand{\T}{\mathbb{T}}

\newcommand{\D}{\mathsf{D}}
\newcommand{\Spin}{\mathbb{S}}
\newcommand{\1}{\bf{1}}

\newcommand{\Lag}{\mathcal{L}}

\newcommand{\Spec}{\mathsf{Spec}}

\newcommand{\BH}{\mathbb{H}}

\begin{document}
\setlength{\baselineskip}{15pt}
\title{Morse-Floer Theory for Super-quadratic Dirac-Geodesics}
\author{Takeshi Isobe$^{(1)}$ \& Ali Maalaoui$^{(2)}$}

\addtocounter{footnote}{1}
\footnotetext{Graduate School of Economics, Hitotsubashi University, 2-1 Naka, Kunitachi, Tokyo 186-8601, Japan:
e-mail:t.isobe@r.hit-u.ac.jp
{\tt{t.isobe@r.hit-u.ac.jp}}}

\addtocounter{footnote}{1}
\footnotetext{Department of mathematics and natural sciences, American University of Ras Al Khaimah, PO Box 10021, Ras Al Khaimah, UAE. E-mail address:
{\tt{ali.maalaoui@aurak.ac.ae}}}

\date{}
\maketitle
\begin{abstract}
In this paper we present the full details of the construction of a Morse-Floer type homology related to the super-quadratic perturbation of the Dirac-geodesic model. This homology is computed explicitly using a Leray-Serre type spectral sequence and this computation leads us to several existence results of Dirac-geodesics.
{\footnotesize \baselineskip 4mm }
\end{abstract}

\medskip

{\footnotesize\textit{}}

\section{Introduction}
The nonlinear sigma model in quantum field theory, consists of two Riemannian manifolds $M$ and $N$ and a Lagrangian $L$ defined by $$L(f)=\frac{1}{2}\int_{M}|df(x)|^{2}dx,$$ where $f:M\to N$ is a map chosen in an adequate space. The critical points of this Lagrangian are solutions to the nonlinear Sigma model and are the well studied harmonic maps. Now, if we set $M=S^{1}$, then $L$ becomes the energy functional of loops which is well studied and whose critical points are geodesics on $N$ (see \cite{K} for more details and results regarding the geodesic problem). Furthermore, through the addition of a supersymmetric structure we obtain an extended model called the supersymmetric Sigma model see \cite{Del},~\cite{Freed} and \cite{J}. 

Specifically, this model is obtained by adding a fermionic action, that is, by assuming that $M$ is a spin manifold, we consider the Lagrangian
$$L(f,\psi)=\frac{1}{2}\int_{M}|df|^{2}+\frac{1}{2}\int_{M}\langle \D_{f} \psi,\psi\rangle,$$
where $\D_{f}$ is the Dirac operator associated with the canonical connection $\nabla^{\Spin(M)}\otimes f^{\ast}\nabla^{TN}$ on $\Spin(M)\otimes f^{\ast}TN$, where $\nabla^{\Spin(M)}$ is the canonical lift of the Levi-Civita connection on $TM$ to the spinor bundle $\Spin(M)\to M$ and $\nabla^{TN}$ is the Levi-Civita connection on $TN$.
There are extensive works on the existence of non-trivial critical points of this functional, called Dirac-Harmonic Maps. One can see for instance \cite{CJLW}, \cite{Amann} and the references therein.

Dirac-geodesics on a compact manifold $N$ are the spinorial analogues of supersymmetric mechanical particles moving on $N$. They are $1$-dimensional versions of the Dirac-harmonic maps first introduced by Chen et al. in~\cite{CJLW}. 

Indeed, if we let  $N$ be a compact Riemannian manifold, then a Dirac-geodesic is obtained as a critical point of the Dirac-geodesic action functional. It is defined on a configuration space $\mathcal{F}^{1,1/2}(S^{1},N)$ defined by
$$\mathcal{F}^{1,1/2}(S^{1},N)=\{(\phi,\psi):\phi\in H^{1}(S^{1},N),\,\psi\in H^{1/2}(S^{1},\Spin(S^{1})\otimes\phi^{\ast}TN)\},$$
where $H^{1}(S^{1},N)=\{\phi\in H^{1}(S^{1},\R^{k}):\text{$\phi(s)\in N$ a.e. $s\in S^{1}$}\}$ is the set of $H^{1}$-loops on $N$ (we may assume without loss of generality $N\subset\R^{k}$ for some $k>1$) and $H^{1/2}(S^{1},\Spin(S^{1})\otimes\phi^{\ast}TN)=\{\psi\in H^{1/2}(S^{1},\Spin(S^{1})\otimes\underline{\R}^{k}:\text{$\psi(s)\in\Spin\otimes T_{\phi(s)}N$ a.e. $s\in S^{1}$}\}$ is the set of $H^{1/2}$-spinors along the loop $\phi$. Here we denote by $\Spin\cong\C$ the $1$-dimensional spinor module, $\Spin(S^{1})\to S^{1}$ is a spinor bundle on $S^{1}$ and $\underline{\R}^{k}\to S^{1}$ is the trivial $\R^{k}$ bundle over $S^{1}$.
Throughout this paper, we assume as in \cite{I} that the metric structure on $\mathcal{F}^{1,1/2}(S^{1},N)$ is induced from the canonical embedding $\mathcal{F}^{1,1/2}(S^{1},N)\subset H^{1}(S^{1},\R^{k})\times H^{1/2}(S^{1},\Spin(S^{1})\otimes\underline{\R}^{k})$. See \cite{I} and $\S3$ for details.
Recall that there are two spin structures on $S^{1}$. Throughout this paper, we fix one of them. For any loop $\phi\in H^{1}(S^{1},N)$, there is a canonical connection $\nabla^{\phi}=\nabla^{\Spin(S^{1})}\otimes\phi^{\ast}\nabla^{TN}$ on $\Spin(S^{1})\otimes\phi^{\ast}TN$. 
Associated to this is a Dirac operator $\D_{\phi}=\frac{\partial}{\partial s}\cdot\nabla^{\phi}_{\frac{\partial}{\partial s}}$, where $\frac{\partial}{\partial s}\cdot$ is the Clifford multiplication by the angular vector field $\frac{\partial}{\partial s}\in TS^{1}$. Under the identification $\Spin\cong\C$, it is simply given by multiplication by $\pm i$. (The choice of the sign is not relevant in this paper). For more details, see~\cite[$\S2$]{I}. The Dirac-geodesic action functional $\Lag$ is defined by
$$\Lag(\phi,\psi)=\frac{1}{2}\int_{S^{1}}\Big|\frac{d\phi}{ds}\Big|^{2}\,ds+\frac{1}{2}\int_{S^{1}}\langle\psi,\D_{\phi}\psi\rangle\,ds,$$
where $\langle\,\cdot\,,\,\cdot\,\rangle$ is the canonical metric on $\Spin(S^{1})\otimes\phi^{\ast}TN$.

Observe that $\Lag$ restricted to $H^{1}(S^{1},N)\times\{0\}$ is again the geodesic energy functional $E(\phi)=\frac{1}{2}\int_{S^{1}}\big|\frac{d\phi}{ds}\big|^{2}\,ds$ whose critical points are closed geodesics on $N$, while for a fixed loop $\phi$, the second term of $\Lag$ is the pure Dirac fermion action and its critical points are harmonic spinors. In this paper, we focus our study by considering a perturbed Dirac-geodesic action functional $\Lag_{H}$ of the following form
$$\Lag_{H}(\phi,\psi)=\frac{1}{2}\int_{S^{1}}\Big|\frac{d\phi}{ds}\Big|^{2}\,ds+\frac{1}{2}\int_{S^{1}}\langle\psi,\D_{\phi}\psi\rangle\,ds-\int_{S^{1}}H(s,\phi(s),\psi(s))\,ds,$$
where the perturbation $H:S^{1}\times\Spin(S^{1})\otimes TN\to\R$ is a smooth function (we write $H=H(s,\phi,\psi)$, where $s\in S^{1}$ and $(\phi,\psi)\in\Spin(S^{1})\otimes TN$, i.e., $\phi\in N$ is a base point and $\psi\in\Spin(S^{1})\otimes T_{\phi}N$ is a point on the fiber over $\phi\in N$).
We call critical points of $\Lag_{H}$ perturbed Dirac-geodesics.

The extended supersymmetric model in quantum mechanics consists of finding supersymmetric particles $(\phi,\psi)$ on $N$ that are critical points of the supersymmetric Dirac-geodesic action $\Lag_{H}$, where we take $H(s,\phi,\psi)=\frac{1}{12}R^{N}_{ikjl}(\phi)\langle\psi^{i},\psi^{j}\rangle\langle\psi^{k},\psi^{l}\rangle$, with $R^{N}_{ikjl}$ a curvature tensor of $N$. See \cite{Freed},~\cite{Del},~\cite{J} 

As far as we know, the only existence results of (perturbed) Dirac-geodesics in the literature was studied in ~\cite{I} for the case of $N=\T$ a flat torus and the case of a compact manifold $N$ with a bumpy metric and ``big" superquadratic perturbation. The pure unperturbed case was studied in \cite{C}, in the case of surfaces via a heat flow approach and explicit formula for Dirac-geodesics was given in the particular cases of the two sphere and the hyperbolic plane. Other than these two cases, the existence of a non-trivial perturbed Dirac-geodesic is widely open. In this paper, we aim to settle this question using a Morse-Floer theoretical approach for superquadratic Dirac-geodesics. Since the spectrum of the Dirac operator $\D_{\phi}$ at any loop $\phi$ is unbounded from below and above, the functional $\Lag_{H}$ is strongly indefinite in the sense that the Morse index and coindex at any critical point are infinite. Thus, we shall construct Floer type homology theory for the functional $\Lag_{H}$ leading us to many existence results. In our setting, we will consider perturbations $H$ satisfying

\begin{equation}
|d^{2}_{\psi\psi}H(s,\phi,\psi)|\le C_{1}(1+|\psi|^{p-1}),
\end{equation}
\begin{equation}
C_{2}|\psi|^{p+1}+2H(s,\phi,\psi)-C_{3}\le\langle\nabla_{\psi}H(s,\phi,\psi),\psi\rangle,
\end{equation}
\begin{equation}
|d^{2}_{\psi\phi}H(s,\phi,\psi)|\le C_{4}(1+|\psi|^{q_{1}}),\quad
|d^{2}_{\phi\phi}H(s,\phi,\psi)|\le C_{4}(1+|\psi|^{q_{2}})
\end{equation}
where in (1.3), $0<q_{i}<p$ for $i=1,2$.

We note that, by integrating (1.1) and (1.3), we have
\begin{equation}
|H(s,\phi,\psi)|\le C(1+|\psi|^{p+1}),\quad|d_{\psi}H(s,\phi,\psi)|\le C(1+|\psi|^{p})
\end{equation}
and 
\begin{equation}
|d_{\phi}H(s,\phi,\psi)|\le C(1+|\psi|^{q_{1}+1})
\end{equation}
for some constant $C>0$.

We call conditions (1.1) and (1.2), superquadratic conditions since they imply that $H(s,\phi,\psi)$ behaves like $|\psi|^{p+1}$ as $|\psi|\to\infty$. We denote by $\BH^{\ell}_{p+1}$ the class of perturbations $H=H(s,\phi,\psi)\in C^{\ell}(S^{1}\times\Spin(S^{1})\otimes TN)$ which satisfies (1.1)--(1.3) above.

Our main result is summarized as follows:

\newtheorem{theorem}{Theorem}[section]
\newtheorem{lemma}{Lemma}[section]
\newtheorem{corollary}{Corollary}[section]
\newtheorem{definition}{Definition}[section]
\newtheorem{proposition}{Proposition}[section]
\newtheorem{remark}{Remark}[section]
\newtheorem{claim}{Claim}[section]

\begin{theorem}
Let us assume that $p\ge 3$, then we have:
\begin{enumerate}[(1)]
\item For $H\in\BH_{p+1}^{3}$, the Morse-Floer homology $HF_{\ast}(\Lag_{H},\mathcal{F}^{1,1/2}(S^{1},N);\Z_{2})$ is well-defined.
\item For any $H,H'\in\BH_{p+1}^{3}$, there is a natural isomorphism \[\Phi_{HH'}:HF_{\ast}(\Lag_{H},\mathcal{F}^{1,1/2}(S^{1},N);\Z_{2})\to HF_{\ast}(\Lag_{H'},\mathcal{F}^{1,1/2}(S^{1},N);\Z_{2}).\] 

Thus the isomorphism class of $HF_{\ast}(\Lag_{H},\mathcal{F}^{1,1/2}(S^{1},N);\Z_{2})$ is independent of $H\in\BH_{p+1}^{3}$. We define the (p+1)-Dirac-geodesic homology with $\Z_{2}$-coefficient, denoted by $DG^{p+1}H_{\ast}(\mathcal{F}^{1,1/2}(S^{1},N);\Z_{2})$, as the isomorphism class of $HF_{\ast}(\Lag_{H},\mathcal{F}^{1,1/2}(S^{1},N);\Z_{2})$ for (any one of) $H\in\BH_{p+1}^{3}$.
\item 
We have a vanishing result: $DG^{p+1}H_{\ast}(\mathcal{F}^{1,1/2}(S^{1},N);\Z_{2})=0$
\item In the special case of flat torus $N=\T_{\Lambda}:=\R^{n}/\Lambda$ ($\Lambda\subset\R^{n}$ is a lattice of rank $n$),
all the above hold under the assumption $p>2$.
\end{enumerate}
\end{theorem}

In the course of our proof, we will prove much more structural results about the homology of $HF_{\ast}(\Lag_{H},\mathcal{F}^{1,1/2}(S^{1},N);\Z_{2})$ defined for our functional $\Lag_{H}$, see $\S8$ for details.
As observed in~\cite{I}, the configuration space $\mathcal{F}^{1,1/2}(S^{1},N)$ has a natural Hilbert bundle structure over the loop space $H^{1}(S^{1},N)$ with projection map $\pi:\mathcal{F}^{1,1/2}(S^{1},N)\to H^{1}(S^{1},N)$ defined by $\pi(\phi,\psi)=\phi$. 
The fiber over $\phi\in H^{1}(S^{1},N)$ is then $\pi^{-1}(\phi)=H^{1/2}(S^{1},\Spin(S^{1})\otimes\phi^{\ast}TN)$. Note that the loop space $H^{1}(S^{1},N)$ is canonically embedded as the zero section in $\mathcal{F}^{1,1/2}(S^{1},N)$, $\displaystyle H^{1}(S^{1},N)\ni\phi\mathop{\to}^{\sim}(\phi,0)\in\mathcal{F}^{1,1/2}(S^{1},N)$. By this identification, we can naturally consider $H^{1}(S^{1},N)$ as a submanifold of $\mathcal{F}^{1,1/2}(S^{1},N)$.
Observe that $\Lag_{H}$ restricted to $H^{1}(S^{1},N)$ is a perturbed geodesic action $$E_{1,H}(\phi)=\frac{1}{2}\int_{S^{1}}|\dot{\phi}|^{2}\,ds-\int_{S^{1}}H(s,\phi(s),0)\,ds,$$
 while its restriction to the fiber over $\phi\in H^{1}(S^{1},N)$ is a perturbed Dirac fermion action $$E_{2,H;\phi}(\psi)=\frac{1}{2}\int_{S^{1}}\langle\D_{\phi}\psi,\psi\rangle\,ds-\int_{S^{1}}H(s,\phi(s),\psi(s))\,ds.$$
 These two functionals were investigated from a Morse theoretical point of view and they are relatively well understood. Indeed, Morse theory for $E_{1,H}$ is a classical study, see~\cite{Palais} and \cite{K} for details and more recently~\cite{AM3}. For $E_{2,H;\phi}$, in~\cite{I3} and~\cite{I2} Morse-Floer homology is constructed and computed for the superquadratic Dirac fermion action. The same problem was also investigated using a Rabinowitz-Floer homology approach in \cite{M}.

 Because of the fiber bundle structure of $\mathcal{F}^{1,1/2}(S^{1},N)\to H^{1}(S^{1},N)$, and the splitting of the functional $\Lag_{H}$ into two parts, one geodesic and one fermionic, one can distinguish a filtration structure induced by the geodesic part via Morse index  making it natural to conjecture that there is a Leray-Serre type spectral sequence converging to the Morse-Floer homology $HF_{\ast}(\Lag_{H},\mathcal{F}^{1,1/2}(S^{1},N);\Z_{2})$. 
 
 The main difference from a regular spectral sequence, is the lack of topological objects obtained from the functionals in the underlying space of variations and this is mainly because of the absence of the classical Morse index and the cell attaching procedure. In this work, in order to compute the homology $HF_{\ast}(\Lag_{H},\mathcal{F}^{1,1/2}(S^{1},N);\Z_{2})$ we construct a Leray-Serre type spectral sequence whose 2nd page is given by  
$$MH_{\ast}(E_{1,H},H^{1}(S^{1},N);\mathcal{HF}_{\ast}(E_{2,H},H^{1/2}(S^{1},\Spin(S^{1})\otimes TN);\Z_{2})),$$ 
where $\mathcal{HF}_{\ast}(E_{2,H},H^{1/2}(S^{1},\Spin(S^{1})\otimes TN);\Z_{2}))$ is a local coefficient system $$\{HF_{\ast}(E_{2,H;\phi},H^{1/2}(S^{1},\Spin(S^{1})\otimes\phi^{\ast}TN);\Z_{2}))\}_{\phi\in H^{1}(S^{1},N)}$$ on $H^{1}(S^{1},N)$, where $HF_{\ast}(E_{2,H;\phi},H^{1/2}(S^{1},\Spin(S^{1})\otimes\phi^{\ast}TN);\Z_{2})$ is the Morse-Floer homology of $E_{2,H;\phi}$ on $H^{1/2}(S^{1},\Spin(S^{1})\otimes\phi^{\ast}TN)$ with $\Z_{2}$-coefficient and
$$MH_{\ast}(E_{1,H},H^{1}(S^{1},N);\mathcal{HF}_{\ast}(E_{2,H},H^{1/2}(S^{1},\Spin(S^{1})\otimes TN);\Z_{2}))$$
is the Morse homology of $E_{1,H}$ on $H^{1}(S^{1},N)$ with local coefficient $$\mathcal{HF}_{\ast}(E_{2,H},H^{1/2}(S^{1},\Spin(S^{1})\otimes TN);\Z_{2}).$$ 
But as we mentioned earlier, in~\cite{I2} the author proved that $$HF_{\ast}(E_{2,H;\phi},H^{1/2}(S^{1},\Spin(S^{1})\otimes\phi^{\ast}TN);\Z_{2})$$ is well-defined and does not depend on the particular choices of $\phi$ and $H\in\mathbb{H}_{p+1}$ (see $\S8$). Thus, the local coefficient system is in fact trivial. Moreover, it vanishes identically. Thus the spectral sequence collapses at the 2nd page and the vanishing of the homology of $\Lag_{H}$ on $\mathcal{F}^{1,1/2}(S^{1},N)$ follows.

Since in the classical use of Morse theory, the non-vanishing of the homology is used to extract existence and multiplicity results, the vanishing of our homology might lead the reader to think that this approach is not adequate for proving existence of critical points. However, in our case, the vanishing is seen as a topological obstruction, and leads to many new existence results for the perturbed Dirac-geodesics.
\begin{theorem}
Assume that $H\in \mathbb{H}^{3}_{p+1}$ satisfies  $\nabla_{\psi}H(\phi,0)=\nabla_{\psi,\psi}H(\phi,0)=\nabla_{\phi,\psi}H(\phi,0)=0$, then either $(N,g)$ has three (perturbed) prime geodesics, two minimal and one of index one, or we have the existence of at least one non-trivial perturbed Dirac Geodesic.
\end{theorem}

\begin{theorem} 
Assume that $H\in \mathbb{H}^{3}_{p+1}$ is even in $\psi$ and satisfies  $\nabla_{\psi}H(\phi,0)=\nabla_{\psi,\psi}H(\phi,0)=\nabla_{\phi,\psi}H(\phi,0)=0$, then we have the existence of at least one non-trivial perturbed Dirac Geodesic. Moreover, if there exists a sequence of integers $a_{k}\to \infty$ such that the Betti numbers $b_{a_{k}}(\Lambda N,\Z_{2})\not =0$, then we have infinitely many non-trivial Dirac-geodesics.
\end{theorem}

As seen in the following theorem, the case of the torus presents easier features.
\begin{theorem}
For the special case $N=\T$, where $\T$ is a flat torus, under the same assumption $\nabla_{\psi}H(\phi,0)=\nabla_{\psi,\psi}H(\phi,0)=\nabla_{\phi,\psi}H(\phi,0)=0$ but without assuming the evenness of $H$ in $\psi$, there exists a non-trivial perturbed Dirac geodesics in each connected component.
\end{theorem}

Our paper is structured into four main parts:

\begin{itemize}

\item The first four Sections are devoted to the construction of the Floer homology while assuming transversality. This is done by introducing an adequate grading to replace the classical Morse index and then by studying the negative gradient flow lines.

\item In Sections 5, 6 and 7, we investigate the transversality conditions and we show that the homology can be defined for a generic set of perturbations. Then, we show that the homology is stable under perturbations by exhibiting a natural isomorphism between the different perturbed Floer homologies.

\item In Section 8, we focus on computing the homology by introducing a Leray-Serre type spectral sequence as described above.

\item In Section 9, we use our homology to prove some existence results for perturbed Dirac-geodesics.
\end{itemize}

\bigskip

\noindent
{\bf Acknowledgement.} The first author (T.I) is partially supported by JSPS KAKENHI Grant Number 15K04947.

\section{Relative Morse index}

We define a relative Morse index of a critical point $(\phi,\psi)\in\mathcal{F}^{1,1/2}(S^{1},N)$ via the spectral flow of the linearization of the Euler-Lagrange equation. So, let us write $\mathcal{A}_{(\phi,\psi),H}$ the bounded self-adjoint realization of the Hessian $d^{2}\Lag_{H}(\phi,\psi)$ with respect to the $H^{1}\times H^{1/2}$ metric on $\mathcal{F}^{1,1/2}(S^{1},N)$:
$$d^{2}\Lag_{H}(\phi,\psi)\Big[\begin{pmatrix}X\\ \xi\end{pmatrix},\begin{pmatrix}Y\\ \zeta\end{pmatrix}\Big]=\Big\langle\mathcal{A}_{(\phi,\psi),H}\begin{pmatrix}X\\ \xi\end{pmatrix},\begin{pmatrix}Y\\ \zeta\end{pmatrix}\Big\rangle_{T_{(\phi,\psi)}\mathcal{F}^{1,1/2}(S^{1},N)},$$
where $(X,\xi),(Y,\zeta)\in T_{(\phi,\psi)}\mathcal{F}^{1,1/2}(S^{1},N)=H^{1}(S^{1},\phi^{\ast}TN)\times H^{1/2}(S^{1},\Spin(S^{1})\otimes\phi^{\ast}TN)$ and the metric on $\mathcal{F}^{1,1/2}(S^{1},N)$ is defined in (3.4).
Thus $\mathcal{A}_{(\phi,\psi),H}:T_{(\phi,\psi)}\mathcal{F}^{1,1/2}(S^{1},N)\to T_{(\phi,\psi)}\mathcal{F}^{1,1/2}(S^{1},N)$ is a bounded self-adjoint operator.

Let us assume that $(\phi,\psi)\in\mathcal{F}^{1,1/2}(S^{1},N)$ is a non-degenerate critical point. We fix a smooth pair $(\phi_{0},\psi_{0})\in\mathcal{F}^{1,1/2}(S^{1},N)$ in the component containing $(\phi,\psi)$
and take a smooth path $\{(\phi_{t},\psi_{t})\}_{t\in[0,1]}$ connecting $(\phi_{0},\psi_{0})$ to $(\phi,\psi)$ in $\mathcal{F}^{1,1/2}(S^{1},N)$, $(\phi_{t=0},\psi_{t=0})=(\phi_{0},\psi_{0})$, $(\phi_{t=1},\psi_{t=1})=(\phi,\psi)$. By the regularity theory, $(\phi,\psi)$ is smooth and we can always assume that $(\phi_{t},\psi_{t})$ is also a smooth pair for all $t$. We then define:

\begin{definition}
We define relative Morse index $\mu_{H}(\phi,\psi)$ by
$$\mu_{H}(\phi,\psi)=-\mathsf{sf}\{\mathcal{A}_{(\phi_{t},\psi_{t}),H}\}_{0\le t\le 1},$$
where $\mathsf{sf}\{\mathcal{A}_{(\phi_{t},\psi_{t}),H}\}_{0\le t\le 1}$ is the spectral flow of a path of operators $\{\mathcal{A}_{(\phi_{t},\psi_{t}),H}\}_{0\le t\le 1}$.
\end{definition}

We need to prove that the quantity above is well-defined, i.e., $\mu_{H}(\phi,\psi)\in\Z$ and does not depend on the choice of the path $\{(\phi_{t},\psi_{t})\}$. 
In fact, we have the following:

\begin{lemma}
The spectral flow $\mathsf{sf}\{\mathcal{A}_{(\phi_{t},\psi_{t}),H}\}_{0\le t\le 1}$ is well-defined and does not depend on the choice of the path $\{(\phi_{t},\psi_{t})\}_{0\le t\le 1}$.
\end{lemma}
\textit{Proof.} Well-definedness follows since, as we will see shortly below, $\{\mathcal{A}_{(\phi_{t},\psi_{t}),H}\}_{0\le t\le 1}$ is a Fredholm family. See~\cite{BW} for general results about the spectral flow for Fredholm family.
As for the independence of the path, the assertion follows if we show that the spectral flow $\mathsf{sf}\{\mathcal{A}_{(\phi_{t},\psi_{t}),H_{t}}\}_{t\in S^{1}}$ is $0$ for any closed path $\{(\phi_{t},\psi_{t})\}_{t\in S^{1}}$, where $H_{t}$ is a loop of perturbations which satisfy (1.1)--(1.3) uniformly with respect to $t$. (This general case is not necessary for the proof of our assertion. However, for later purpose we treat a slightly general case). For simplicity, we write $\mathcal{A}_{t}=\mathcal{A}_{(\phi_{t},\psi_{t}),H_{t}}$.\\
 Let us denote by $P_{t}(s):T_{\phi_{0}(s)}N\to T_{\phi_{t}(s)}N$ the parallel translation along the path $[0,t]\ni\tau\mapsto\phi_{\tau}(s)\in N$. We define the operators $\mathcal{B}_{t}=P_{t}^{-1}\circ\mathcal{A}_{t}\circ P_{t}$ parametrized by $t\in S^{1}$. Note that $$\mathcal{B}_{t}:H^{1}(S^{1},\phi_{0}^{\ast}TN)\times H^{1/2}(S^{1},\Spin(S^{1})\otimes\phi_{0}^{\ast}TN)\to H^{1}(S^{1},\phi_{0}^{\ast}TN)\times H^{1/2}(S^{1},\Spin(S^{1})\otimes\phi_{0}^{\ast}TN),$$
for all $t\in S^{1}$ and $\mathsf{sf}\{\mathcal{B}_{t}\}_{t\in S^{1}}=\mathsf{sf}\{\mathcal{A}_{t}\}_{t\in S^{1}}$. We thus consider $\mathsf{sf}\{\mathcal{B}_{t}\}_{t\in S^{1}}$ instead of considering $\mathsf{sf}\{\mathcal{A}_{t}\}_{t\in S^{1}}$.
Note that $\mathcal{A}_{t}$ takes the following block decomposition with respect to the canonical decomposition $T_{(\phi,\psi)}\mathcal{F}^{1,1/2}(S^{1},N)=H^{1}(S^{1},\phi^{\ast}TN)\times H^{1/2}(S^{1},\Spin(S^{1})\otimes\phi^{\ast}TN)$ (see $\S3.1$):
$$\mathcal{A}_{t}=\begin{pmatrix}\mathcal{A}_{11,t}&\mathcal{A}_{12,t}\\
\mathcal{A}_{21,t}&\mathcal{A}_{22,t}\end{pmatrix}.$$
The components $\mathcal{A}_{ij,t}$ are given as follows:
\begin{align}
\mathcal{A}_{11,t}u=(-\Delta+1)^{-1}\Big(&-\nabla_{s}\nabla_{s}X-R(X,\partial_{s}\phi_{t})\partial_{s}\phi_{t}+\frac{1}{2}\nabla_{X}R(\phi_{t})\langle\psi_{t},\partial_{s}\phi_{t}\cdot\psi_{t}\rangle\notag\\
&+\frac{1}{2}R(\phi_{t})\langle\psi_{t},\nabla_{s}X\cdot\psi_{t}\rangle-H_{t,\phi\phi}(\,\cdot\,,\phi_{t},\psi_{t})[X]\Big),
\end{align}
\begin{align}
\mathcal{A}_{12,t}\xi=\Big(-\Delta+1)^{-1}(\frac{1}{2}R(\phi_{t})\langle\xi,\partial_{s}\phi_{t}\cdot\psi_{t}\rangle+\frac{1}{2}R(\phi_{t})\langle\psi_{t},\partial_{s}\phi\cdot\xi\rangle-H_{t,\phi\psi}(\,\cdot\,,\phi_{t},\psi_{t})[\xi]\Big),
\end{align}
\begin{align}
&\mathcal{A}_{21,t}X=(1+|\D|)^{-1}\Big(e_{1}\cdot\partial_{s}\psi^{i}X^{m}\Gamma_{mj}^{k}(\phi_{t})\otimes\frac{\partial}{\partial y^{k}}(\phi_{t})+\nabla_{X}\Gamma_{jk}^{i}(\phi_{t})\partial_{s}\phi_{t}^{j}e_{1}\cdot\psi_{t}^{k}\otimes\frac{\partial}{\partial y^{i}}(\phi_{t})\notag\\
&+\Gamma_{jk}^{i}(\phi_{t})\partial_{s}X^{j}e_{1}\cdot\psi^{k}_{t}\otimes\frac{\partial}{\partial y^{i}}(\phi_{t})+\Gamma_{jk}^{i}(\phi_{t})\partial_{s}\phi_{t}^{j}e_{1}\cdot\psi_{t}^{k}\Gamma_{li}^{m}(\phi_{t})X^{l}\otimes\frac{\partial}{\partial y^{m}}(\phi_{t})-H_{t,\psi\phi}(\,\cdot\,,\phi,\psi)[X]\Big),
\end{align}
\begin{align}
\mathcal{A}_{22,t}\xi=(1+|\D|)^{-1}(\D_{\phi_{t}}\xi-H_{t,\psi\psi}(\,\cdot\,\phi,\psi)[\xi]),
\end{align}
where $\Delta=\frac{\partial^{2}}{\partial s^{2}}$, $\D=\D^{\Spin(S^{1})}\otimes\bm{1}_{\underline{\R}^{k}}$ and $\nabla_{s}=\nabla_{\frac{\partial}{\partial s}}$.

By the above form of $\mathcal{A}_{ij}$, it is easy to see that $\mathcal{A}_{t}$ can be written as
\begin{equation}
\mathcal{A}_{t}=\begin{pmatrix}\mathcal{A}^{0}_{11,t}&O\\
O&\mathcal{A}_{22,t}^{0}\end{pmatrix}+\mathcal{K}_{t},
\end{equation}
where
$$\mathcal{A}_{11,t}^{0}=(-\Delta+1)^{-1}(-\nabla_{s}\nabla_{s}),$$
$$\mathcal{A}^{0}_{22,t}=(1+|\D|)^{-1}\D_{\phi_{t}}$$
and $\mathcal{K}_{t}$ is compact.

Now, we consider $\{\mathcal{B}_{t}\}_{t\in S^{1}}$. Observe that $\mathcal{B}_{t}$ is a $0$-th order elliptic operator acting on sections of the bundle $\phi_{0}^{\ast}TN\oplus\Spin(S^{1})\otimes\phi_{0}^{\ast}TN\to S^{1}$ and its principal symbol is independent of $t$. Considering $S^{1}$ as the quotient $S^{1}=[0,1]/\{0\}\sim\{1\}$, we also have $\mathcal{B}_{1}=P_{1}^{-1}\circ\mathcal{B}_{0}\circ P_{1}$, where $P_{1}$ is considered as an automorphism of the bundle $\phi_{0}^{\ast}TN\oplus\Spin(S^{1})\otimes\phi_{0}^{\ast}TN\to S^{1}$ (acting $\phi_{0}^{\ast}TN$ parts only). In such a case, the spectral flow $\mathsf{sf}\{\mathcal{B}_{t}\}_{t\in S^{1}}$ is given by the Booss-Wojciechowski's desuspension formula~\cite[Theorem 17.17]{BW} as follows:
\begin{equation}
\mathsf{sf}\{\mathcal{B}_{t}\}_{t\in S^{1}}=\text{ind}(P_{\ge 0}(\mathcal{B}_{0})-P_{1}P_{<0}(\mathcal{B}_{0})),
\end{equation}
where $P_{\ge 0}(\mathcal{B}_{0})$ and $P_{<0}(\mathcal{B}_{0})$ are spectral projections of $\mathcal{B}_{0}$ to the spectral set $[0,\infty)$ and $(-\infty,0)$, respectively; $P_{\ge 0}(\mathcal{B}_{0})=\bm{1}_{[0,\infty)}(\mathcal{B}_{0})$, $P_{<0}(\mathcal{B}_{0})=\bm{1}_{(-\infty,0)}(\mathcal{B}_{0})$.

Define $\mathcal{A}_{t}^{0}=\begin{pmatrix}\mathcal{A}^{0}_{11,t}&O\\
O&\mathcal{A}_{22,t}^{0}\end{pmatrix}$ and $\mathcal{B}_{t}^{0}=P_{t}^{-1}\circ\mathcal{A}_{t}^{0}\circ P_{t}=:\begin{pmatrix}\mathcal{B}_{t,11}^{0}&O\\
O&\mathcal{B}_{t,22}^{2}\end{pmatrix}$. Since the difference $\mathcal{B}_{0}^{0}-\mathcal{B}_{0}$ is a compact operator, $P_{\ge 0}(\mathcal{B}_{0})-P_{\ge 0}(\mathcal{B}_{0}^{0})$ and $P_{<0}(\mathcal{B}_{0})-P_{<0}(\mathcal{B}_{0}^{0})$ are compact. Thus by the index formula (2.6) and $P_{\ge 0}(\mathcal{B}_{0}^{0})=P_{\ge 0}(\mathcal{B}_{0,11}^{0})\oplus P_{\ge 0}(\mathcal{B}_{0,22}^{0})$, $P_{<0}(\mathcal{B}_{0}^{0})=P_{<0}(\mathcal{B}_{0,11}^{0})\oplus P_{<0}(\mathcal{B}_{0,22}^{0})$, we have
\begin{align}
\mathsf{sf}\{\mathcal{B}_{t}\}_{t\in S^{1}}&=\text{ind}(P_{\ge 0}(\mathcal{B}_{0}^{0})-P_{1}P_{<0}(\mathcal{B}_{0}^{0}))\notag\\
&=\text{ind}(P_{\ge 0}(\mathcal{B}_{0,11}^{0})-P_{1}P_{<0}(\mathcal{B}_{0,11}^{0})+\text{ind}(P_{\ge 0}(\mathcal{B}_{0,22}^{0})-P_{1}P_{<0}(\mathcal{B}_{0,22}^{0}))\notag\\
&=0+\text{ind}(P_{\ge 0}(\mathcal{B}_{0,22}^{0})-P_{1}P_{<0}(\mathcal{B}_{0,22}^{0}))\notag\\
&=\text{ind}(P_{\ge 0}(\D_{\phi_{0}})-P_{1}P_{<0}(\D_{\phi_{0}}))\notag\\
&=\mathsf{sf}\{P_{t}^{-1}\circ\D_{\phi_{t}}\circ P_{t}\}_{t\in S^{1}}\notag\\
&=\mathsf{sf}\{\D_{\phi_{t}}\}_{t\in S^{1}},
\end{align}
where we have used in the third line the fact that $\mathcal{B}_{0,11}^{0}$ is a non-negative operator so that $P_{\ge 0}(\mathcal{B}_{0,11}^{0})=\text{id}$ and $P_{<0}(\mathcal{B}_{0,11}^{0})=0$.

On the other hand, by the Atiyah-Patodi-Singer index theorem \cite{APS},~\cite{BW}, the spectral flow $\mathsf{sf}\{\D_{\phi_{t}}\}_{t\in S^{1}}$ is given by the Fredholm index of the operator $\frac{\partial}{\partial t}+\D_{\phi_{t}}$ defined on $S^{1}\times S^{1}$, $\mathsf{sf}\{\D_{\phi_{t}}\}_{t\in S^{1}}=\text{ind}\big(\frac{\partial}{\partial t}+\D_{\phi_{t}}\big)$, where
\begin{equation}
\frac{\partial}{\partial t}+\D_{\phi_{t}}:H^{1}(S^{1}\times S^{1},\Spin(S^{1}\times S^{1})^{+}\otimes\phi_{\cdot}^{\ast}TN)\to L^{2}(S^{1}\times S^{1},\Spin(S^{1}\times S^{1})^{-}\otimes\phi_{\cdot}^{\ast}TN)
\end{equation}
is the Dirac operator, where $\Spin(S^{1}\times S^{1})^{\pm}$ is the $\pm$-spinor bundle associated with a spin structure on $S^{1}\times S^{1}$ whose restriction to the first summand gives the trivial spin structure on $S^{1}$, while the restriction to the second summand gives the one of $\Spin(S^{1})\to S^{1}$.
By the Atiyah-Singer index theorem, we have
\begin{align*}
\text{ind}\Big(\frac{\partial}{\partial t}+\D_{\phi_{t}}\Big)&=\int_{S^{1}\times S^{1}}\mathsf{ch}(\phi^{\ast}TN\otimes\C)\hat{\mathsf{A}}(S^{1}\times S^{1})\\
&=\int_{S^{1}\times S^{1}}c_{1}(\phi^{\ast}TN\otimes\C)=\int_{S^{1}\times S^{1}}\phi^{\ast}c_{1}(TN\otimes\C)=0,
\end{align*}
where $\mathsf{ch}$ is the Chern character and $\hat{\mathsf{A}}$ is the $\hat{A}$-genus. 
The last equation follows from $c_{1}(TN\otimes\C)=0$ as the first Chern class of the complexification of the real bundle $TN$.
Thus, by (2.7) we have $\mathsf{sf}\{\mathcal{A}_{t}\}_{t\in S^{1}}=\mathsf\{\mathcal{B}_{t}\}_{t\in S^{1}}=0$. This completes the proof.
\hfill$\Box$

\section{The negative gradient flow of $\Lag_{H}$ on $\mathcal{F}^{1/2}(S^{1},N)$}
\subsection{The gradient flow equation}

In the following, we assume that $H$ satisfies (1.1)--(1.3) for some $p\ge 3$. For the case of the flat tori $N=\T_{\Lambda}$ as in the statement of Theorem 1.1 (4), we only need to assume $p>2$.

We assume that $N$ is isometrically embedded in $\R^{k}$. Thus, we consider $H^{1}(S^{1},N)\subset H^{1}(S^{1},\R^{k})$ and $H^{1/2}(S^{1},\Spin(S^{1})\otimes\phi^{\ast}TN)\subset H^{1/2}(S^{1},\Spin(S^{1})\otimes\underline{\R}^{k})$, $\underline{\R}^{k}=N\times\R^{k}\to N$, the trivial bundle over $N$. The metric on $\mathcal{F}^{1,1/2}(S^{1},N)$ is induced from these imbeddings. To describe this metric, recall that there is a canonical identification
\begin{equation}
T_{(\phi,\psi)}\mathcal{F}^{1,1/2}(S^{1},N)=H^{1}(S^{1},\phi^{\ast}TN)\times H^{1/2}(S^{1},\Spin(S^{1})\otimes\phi^{\ast}TN).
\end{equation}
It is described as follows (see~\cite[$\S 3.3$]{I} for more details): Let $(-1,1)\ni t\mapsto (\phi_{t},\psi_{t})\in\mathcal{F}^{1,1/2}(S^{1},N)$ be a smooth path passing through $(\phi,\psi)\in\mathcal{F}^{1,1/2}(S^{1},N)$ at $t=0$, $(\phi_{t},\psi_{t})\big|_{t=0}=(\phi,\psi)$. We set $X:=\frac{d}{dt}\big|_{t=0}\phi_{t}\in H^{1}(S^{1},\phi^{\ast}TN)$. Using local coordinate $y=(y^{k})$ on $N$, we write $\psi_{t}=\psi_{t}^{i}\otimes\frac{\partial}{\partial y^{i}}(\phi_{t})$, where $\psi_{t}^{i}\in H^{1/2}(S^{1},\Spin(S^{1}))$. Then the variation vector field $\frac{d}{dt}\big|_{t=0}\psi_{t}$ is given by
\begin{align}
\frac{d}{dt}\Big|_{t=0}\psi_{t}&=\frac{d}{dt}\Big|_{t=0}\psi_{t}^{i}\otimes\frac{\partial}{\partial y^{i}}+\psi^{i}\otimes\nabla_{\frac{\partial}{\partial t}}\frac{\partial}{\partial y^{i}}(\phi_{t})\notag\\
&=\Big(\frac{d}{dt}\Big|_{t=0}\psi_{t}^{k}+X^{j}\Gamma_{ji}^{k}(\phi)\psi^{i}\Big)\otimes\frac{\partial}{\partial y^{k}}(\phi),
\end{align}Christoffel
where $\Gamma_{ij}^{k}$ is the Christopher symbol of the metric on $N$.
We set $\xi^{k}=\frac{d}{dt}\Big|_{t=0}\psi_{t}^{k}$ and $\xi=\xi^{k}\otimes\frac{\partial}{\partial y^{k}}(\phi)$. Then the identification (3.1) is given by
\begin{equation}
\frac{d}{dt}\Big|_{t=0}(\phi_{t},\psi_{t})\mapsto(X,\xi).
\end{equation}
Using the identification (3.1), the metric on $\mathcal{F}^{1,1/2}(S^{1},N)$ is defined by
\begin{equation}
\Big\langle\begin{pmatrix}X\\ \xi\end{pmatrix},\begin{pmatrix}Y\\ \zeta\end{pmatrix}\Big\rangle_{T_{(\phi,\psi)}\mathcal{F}^{1,1/2}}=\langle X,Y\rangle_{H^{1}}+\langle\xi,\zeta\rangle_{H^{1/2}},
\end{equation}
where $X,Y\in H^{1}(S^{1},\phi^{\ast}TN)\subset H^{1}(S^{1},\R^{k})$, $\langle X,Y\rangle_{H^{1}}=((-\Delta+1)X,Y)_{L^{2}}$ ($\Delta=\frac{\partial^{2}}{\partial s^{2}}$) and $\xi,\zeta\in H^{1/2}(S^{1},\Spin(S^{1})\otimes\phi^{\ast}TN)\subset H^{1/2}(S^{1},\Spin(S^{1})\otimes\underline{\R}^{k})$, $\langle\xi,\zeta\rangle_{H^{1/2}}=((1+|\D|)\xi,\zeta)_{L^{2}}$ ($\D=\D^{\Spin(S^{1})}\otimes{\bm 1}_{\underline{\R}^{k}}$). With this metric, $\mathcal{F}^{1,1/2}(S^{1},N)$ becomes a complete Hilbert manifold.
The gradient of $\Lag_{H}$ with respect to this metric is denoted by $\nabla_{1,1/2}\Lag_{H}$. According to the decomposition (3.1), $\nabla_{1,1/2}\Lag_{H}$ reads as
$$\nabla_{1,1/2}\Lag_{H}(\phi,\psi)=\nabla_{\phi}\Lag_{H}(\phi,\psi)\oplus\nabla_{\psi}\Lag_{H}(\phi,\psi).$$
The components are given in the following Proposition:

\begin{proposition}[\cite{I}] For $(\phi,\psi)\in\mathcal{F}^{1,1/2}(S^{1},N)$, we have the following:
\begin{equation}
\nabla_{\phi}\Lag_{H}(\phi,\psi)=(-\Delta+1)^{-1}\Big(-\nabla_{s}\partial_{s}\phi
+\frac{1}{2}R(\phi)\langle\psi,\partial_{s}\phi\cdot\psi\rangle-\nabla_{\phi}H(s,\phi,\psi)\Big),
\end{equation}
\begin{equation}
\nabla_{\psi}\Lag_{H}(\phi,\psi)=(1+|\D|)^{-1}(\D_{\phi}-\nabla_{\psi}H(s,\phi,\psi)),
\end{equation}
where 
$$R(\phi)\langle\psi,\partial_{s}\phi\cdot\psi\rangle=\Big\langle\psi,\partial_{s}\cdot\psi^{i}\otimes\frac{\partial}{\partial y^{j}}(\phi)\Big\rangle\partial_{s}\phi^{l}R^{j}_{iml}(\phi)g^{ms}(\phi)\frac{\partial}{\partial y^{s}}(\phi)$$
and $R_{ijk}^{l}$ is the curvature tensor of the metric on $N$.
\end{proposition}

See~\cite{I} for the details of the derivation of the above formula.

Thus the negative gradient flow equation $\frac{d}{dt}\begin{pmatrix}\phi\\ \psi\end{pmatrix}=-\nabla_{1,1/2}\Lag_{H}(\phi,\psi)$ takes the following form:
\begin{equation}
\frac{\partial\phi}{\partial t}=(-\Delta+1)^{-1}\Big(\nabla_{s}\partial_{s}\phi
-\frac{1}{2}R(\phi)\langle\psi,\partial_{s}\phi\cdot\psi\rangle+\nabla_{\phi}H(s,\phi,\psi)\Big),
\end{equation}
\begin{equation}
\frac{\partial\psi}{\partial t}=-(1+|\D|)^{-1}(\D_{\phi}-\nabla_{\psi}H(s,\phi,\psi)).
\end{equation}

\subsection{Regularity and estimates for the negative gradient flow}

In this section, we establish regularity properties of solutions to the negative gradient flow equation (3.7), (3.8). We first prove the following:

\begin{proposition}
Let $(\phi(t),\psi(t))\in C^{1}(\R,\mathcal{F}^{1,1/2}(S^{1},N))$ be a solution to the negative gradient flow equation
\begin{equation}
\frac{d}{dt}\begin{pmatrix}\phi\\ \psi\end{pmatrix}=-\nabla_{1,1/2}\Lag_{H}(\phi,\psi)
\end{equation}
which satisfies the condition $\sup_{t\in\R}|\Lag_{H}(\phi(t),\psi(t))|=:C_{0}<+\infty$. Then there exists $C(C_{0})>0$ such that
$$\sup_{t\in\R}\|\partial_{s}\phi(t)\|_{L^{2}(S^{1})}+\sup_{t\in\R}\|\psi(t)\|_{H^{1/2}(S^{1})}\le C(C_{0}).$$
\end{proposition}

\bigskip

To prove the above proposition, we prepare some preliminary materials. We denote by $\mathsf{H}_{0}^{-}$, $\mathsf{H}_{0}^{0}$ and $\mathsf{H}_{0}^{+}$ the negative, zero and positive subspaces of $H^{1/2}(S^{1},\Spin(S^{1})\otimes\underline{\R}^{k})$ with respect to the Dirac operator $\D=\D^{\Spin(S^{1})}\otimes\bm{1}_{\underline{\R}^{k}}$. Thus we have
\begin{equation}
H^{1/2}(S^{1},\Spin(S^{1})\otimes\underline{\R}^{k})=\mathsf{H}_{0}^{-}\oplus\mathsf{H}^{0}_{0}\oplus\mathsf{H}^{+}_{0}.
\end{equation}
We denote by $P_{0}^{-}$, $P_{0}^{0}$ and $P_{0}^{+}$ the spectral projections onto $\mathsf{H}_{0}^{-}$, $\mathsf{H}_{0}^{0}$ and $\mathsf{H}_{0}^{+}$, respectively. For a spinor $\psi\in H^{1/2}(S^{1},\Spin(S^{1})\otimes\underline{\R}^{k})$, we denote $\psi^{-}_{0}=P_{0}^{-}\psi$, $\psi_{0}^{0}=P_{0}^{0}\psi$ and $\psi_{0}^{+}=P_{0}^{+}\psi$. Thus, we have
$$\psi=\psi_{0}^{-1}+\psi_{0}^{0}+\psi_{0}^{+}.$$
By the embedding $\iota:N\hookrightarrow\R^{k}$, for $\psi=\psi^{k}\otimes\frac{\partial}{\partial y^{k}}(\phi)\in H^{1/2}(S^{1},\Spin(S^{1})\otimes\phi^{\ast}TN)$ we have $\tilde{\psi}:=i_{\ast}\psi=\psi^{k}\otimes\iota_{\ast}\big(\frac{\partial}{\partial y^{k}}(\phi)\big)\in H^{1/2}(S^{1},\Spin(S^{1})\otimes\underline{\R}^{k})$. Note that, we have
\begin{align*}
\D\tilde{\psi}&=\frac{\partial}{\partial s}\cdot D_{\frac{\partial}{\partial s}}\tilde{\psi}\\
&=\frac{\partial}{\partial s}\cdot\Big(\nabla_{\frac{\partial}{\partial s}}\psi^{k}\otimes\iota_{\ast}\Big(\frac{\partial}{\partial y^{k}}(\phi)\Big)+\psi^{k}\otimes(\nabla_{\frac{\partial}{\partial s}}i_{\ast})\Big(\frac{\partial}{\partial y^{k}}(\phi)\Big)+\psi^{k}\otimes\iota_{\ast}\Big(\nabla_{\frac{\partial}{\partial s}}\frac{\partial}{\partial y^{k}}(\phi)\Big)\Big)\\
&=\iota_{\ast}\Big(\frac{\partial}{\partial s}\cdot\nabla_{\frac{\partial}{\partial s}}\psi^{k}\otimes\frac{\partial}{\partial y^{k}}(\phi)+\frac{\partial}{\partial s}\cdot\psi^{k}\otimes\nabla_{\frac{\partial}{\partial s}}\Big(\frac{\partial}{\partial y^{k}}(\phi)\Big)\Big)+\frac{\partial}{\partial s}\cdot\psi^{k}\otimes(\nabla_{\frac{\partial}{\partial s}}\iota_{\ast})\Big(\frac{\partial}{\partial y^{k}}(\phi)\Big)\\
&=\iota_{\ast}(\D_{\phi}\psi)+A\Big(\frac{\partial\phi}{\partial s},\frac{\partial}{\partial s}\cdot\psi\Big),
\end{align*}
where $D$ is the canonical connection (i.e., flat one) on $\Spin(S^{1})\otimes\underline{\R}^{k}$, $A=\nabla\iota_{\ast}$ is the second fundamental form of the embedding $\iota:N\hookrightarrow\R^{k}$ and $A\big(\frac{\partial\phi}{\partial s},\frac{\partial}{\partial s}\cdot\psi\big):=\frac{\partial}{\partial s}\cdot\psi^{k}\otimes(\nabla_{\frac{\partial}{\partial s}}\iota_{\ast})\big(\frac{\partial}{\partial y^{k}}(\phi)\big)$.
Therefore, $\iota_{\ast}(\D_{\phi}\psi)$ is expressed, in terms of the Dirac operator $\D$ as
\begin{equation}
\iota_{\ast}(\D_{\phi}\psi)=\D\tilde{\psi}-A(\partial_{s}\phi,\partial_{s}\cdot\psi).
\end{equation}
In the following, we identify $\tilde{\psi}$ with $\psi$ by $\iota_{\ast}$ and simply write $\tilde{\psi}=\psi$.

We first prove the following:

\begin{lemma} Under the assumption of Proposition 3.2, there exists a positive constant $C(C_{0})$ depending only on $C_{0}$ such that the following holds for all $\tau\in\R$:
\begin{equation}
\int_{\tau}^{\tau+1}\|\phi(t)\|_{H^{1}(S^{1})}^{2}\,dt\le C(C_{0}),
\end{equation}
\begin{equation}
\int_{\tau}^{\tau+1}\|\psi(t)\|_{H^{1/2}(S^{1})}\,dt\le C(C_{0}),
\end{equation}
\begin{equation}
\int_{\tau}^{\tau+1}\|\psi(t)\|^{p+1}_{L^{p+1}(S^{1})}\,dt\le C(C_{0}).
\end{equation}
\end{lemma}
{\it Proof.} Let $-\infty<a<b<+\infty$ be arbitrary. Since $(\phi(t),\psi(t))$ is a solution to the negative gradient flow equation (3.9), we have
\begin{align*}
&\int_{a}^{b}\|\partial_{t}\phi(t)\|_{H^{1}(S^{1})}^{2}+\|\partial_{t}\psi(t)\|_{H^{1/2}(S^{1})}^{2}\,dt\notag\\
&=-\int_{a}^{b}\Big\langle\nabla_{1,1/2}\Lag_{H}(\phi(t),\psi(t)),\begin{pmatrix}\partial_{t}\phi(t)\\ \partial_{t}\psi(t)\end{pmatrix}\Big\rangle_{T_{(\phi(t),\psi(t))}\mathcal{F}^{1,1/2}}\,dt\\
&=-\int_{a}^{b}\frac{d}{dt}\Lag_{H}(\phi(t),\psi(t))\,dt\\
&=\Lag_{H}(\phi(a),\psi(a))-\Lag_{H}(\phi(b),\psi(b))\le 2C_{0}
\end{align*}
by our assumption $\sup_{t\in\R}|\Lag_{H}(\phi(t),\psi(t))|=:C_{0}$. Since $-\infty<a<b<+\infty$ were arbitrary, we have
\begin{equation}
\int_{-\infty}^{\infty}\|\partial_{t}\phi(t)\|_{H^{1/2}(S^{1})}^{2}\,dt+\int_{-\infty}^{\infty}\|\partial_{t}\psi(t)\|_{H^{1/2}(S^{1})}^{2}\,dt\le 2C_{0}.
\end{equation}

Under the identification $\begin{pmatrix}0\\ \psi(t)\end{pmatrix}\in H^{1}(S^{1},\phi(t)^{\ast}TN)\times H^{1/2}(S^{1},\Spin(S^{1})\otimes\phi(t)^{\ast}TN)=T_{(\phi(t),\psi(t))}\mathcal{F}^{1,1/2}$, taking the inner product with (3.9), we have
\begin{align*}
(\partial_{t}\psi(t),\psi(t))_{H^{1/2}(S^{1})}&=-\Big\langle\nabla_{1,1/2}\Lag_{H}(\phi(t),\psi(t)),\begin{pmatrix}0\\ \psi(t)\end{pmatrix}\Big\rangle_{T_{(\phi(t),\psi(t))}\mathcal{F}^{1,1/2}}\\
&=-\int_{S^{1}}\langle\psi(t),\D_{\phi(t)}\psi(t)\rangle\,ds+\int_{S^{1}}\langle\nabla_{\psi}H(s,\phi(t),\psi(t)),\psi(t)\rangle\,ds.
\end{align*}
By (1.2), we thus have
\begin{align}
(\partial_{t}\psi(t),\psi(t))_{H^{1/2}(S^{1})}+2\Lag_{H}(\phi(t),\psi(t))&=\int_{S^{1}}|\partial_{s}\phi(t)|^{2}\,ds-2\int_{S^{1}}H(s,\phi(t),\psi(t))\,ds\notag\\
&\quad+\int_{S^{1}}\langle\nabla_{\psi}H(s,\phi(t),\psi(t)),\psi(t)\rangle\,ds\notag\\
&\ge\int_{S^{1}}|\partial_{s}\phi(t)|^{2}\,ds+C_{2}\int_{S^{1}}|\psi(t)|^{p+1}\,ds-C.
\end{align}
Combining (3.16) with the assumption $\sup_{t\in\R}|\Lag_{H}(\phi(t),\psi(t))|=C_{0}<+\infty$, we have
\begin{equation}
\int_{S^{1}}|\partial_{s}\phi(t)|^{2}\,ds+\int_{S^{1}}|\psi(t)|^{p+1}\,ds\le 2C_{0}+C\|\psi(t)\|_{H^{1/2}(S^{1})}\|\partial_{t}\psi(t)\|_{H^{1/2}(S^{1})}.
\end{equation}
On the other hand, there exists a positive constant $\lambda_{+}>0$ (depending only on the smallest positive eigenvalue of $\D$) such that
\begin{align}
\|\psi_{0}^{+}(t)\|_{H^{1/2}(S^{1})}^{2}&\le \lambda_{+}\int_{S^{1}}\langle\psi_{0}^{+},\D\psi\rangle\,ds\notag\\
&\le\lambda_{+}\int_{S^{1}}\langle\psi_{0}^{+},\D_{\phi}\psi+A(\partial_{s}\cdot\phi,\partial_{s}\cdot\psi)\rangle\,ds\notag\\
&\le\lambda_{+}\int_{S^{1}}\langle\psi_{0}^{+},\D_{\phi}\psi\rangle\,ds+C\int_{S^{1}}|\partial_{s}\phi||\psi_{0}^{+}||\psi|\,ds\notag\\
&\le\lambda_{+}\int_{S^{1}}\langle P_{\phi(t)}\psi_{0}^{+}(t),\D_{\phi(t)}\psi(t)\rangle\,ds+C\int_{S^{1}}|\partial_{s}\phi||\psi_{0}^{+}||\psi|\,ds,
\end{align}
where we have used (3.11) in the second inequality. Note that, since $(\D_{\phi(t)}\psi(t))(s)\in\Spin(S^{1})\otimes T_{\phi(t)(s)}N$, the last inequality of (3.18) follows, where $P_{y}:\R^{k}\to T_{y}N$ is the orthogonal projection. We also note that, by Lemma 10.2 in the Appendix, $P_{\phi(t)}\psi_{0}^{+}(t)\in H^{1/2}(S^{1},\Spin(S^{1})\otimes\phi(t)^{\ast}TN)$ so that the integral $\int_{S^{1}}\langle P_{\phi(t)}\psi_{0}^{+}(t),\D_{\phi(t)}\psi(t)\rangle\,ds$ is well-defined.

To estimate the integral $\int_{S^{1}}\langle P_{\phi(t)}\psi(t),\D_{\phi(t)}\psi(t)\rangle\,ds$, we observe that $\begin{pmatrix}0\\P_{\phi(t)}\psi_{0}^{+}(t)\end{pmatrix}\in T_{(\phi(t),\psi(t))}\mathcal{F}^{1,1/2}$ so that we have
\begin{align}
&\Big|\Big\langle d\Lag_{H}(\phi(t),\psi(t)),\begin{pmatrix}0\\P_{\phi(t)}\psi_{0}^{+}(t)\end{pmatrix}\Big\rangle\Big|\notag\\
=&~\biggl|\int_{S^{1}}\langle P_{\phi(t)}\psi_{0}^{+}(t),\D_{\phi(t)}\psi(t)\rangle\,ds-\int_{S^{1}}\langle\nabla_{\psi}H(s,\phi(t),\psi(t)),P_{\phi(t)}\psi_{0}^{+}(t)\rangle\,ds\biggr|\notag\\
=&~-\Big(\begin{pmatrix}\partial_{t}\phi(t)\\ \partial_{t}\psi(t)\end{pmatrix},\begin{pmatrix}0\\ P_{\phi(t)}\psi_{0}^{+}(t)\end{pmatrix}\Big)_{T_{(\phi(t),\psi(t))}\mathcal{F}^{1,1/2}}\notag\\
=&~-\langle\partial_{t}\psi(t),P_{\phi(t)}\psi_{0}^{+}(t)\rangle_{H^{1/2}(S^{1})}\notag\\
\le&~\|\partial_{t}\psi(t)\|_{H^{1/2}(S^{1})}\|P_{\phi(t)}\psi_{0}^{+}(t)\|_{H^{1/2}(S^{1})}\notag\\
\le&~C(\|\partial_{t}\psi(t)\|_{H^{1/2}(S^{1})}\|\psi(t)\|_{H^{1/2}(S^{1})}+\|\partial_{t}\psi(t)\|_{H^{1/2}(S^{1})}\|\partial_{s}\phi(t)\|_{L^{2}(S^{1})}^{1/2}\|\psi(t)\|_{L^{4}(S^{1})}),
\end{align}
where in the last inequality, we have used (10.8) in the Appendix.

Combining (3.18) and (3.19), we obtain
\begin{align}
\|\psi_{0}^{+}(t)\|_{H^{1/2}(S^{1})}^{2}&\le C(\|\partial_{t}\psi(t)\|_{H^{1/2}(S^{1})}\|\psi(t)\|_{H^{1/2}(S^{1})}+\|\partial_{t}\psi(t)\|_{H^{1/2}(S^{1})}\|\partial_{s}\phi(t)\|_{L^{2}(S^{1})}^{1/2}\|\psi(t)\|_{L^{4}(S^{1})})\notag\\
&+\int_{S^{1}}|\nabla_{\psi}H(s,\phi(t),\psi(t))||P_{\phi(t)}\psi_{0}^{+}(t)|\,ds+C\int_{S^{1}}|\partial_{s}\phi(t)||\psi_{0}^{+}(t)||\psi(t)|\,ds\notag\\
&\le C(\|\partial_{t}\psi(t)\|_{H^{1/2}(S^{1})}\|\psi(t)\|_{H^{1/2}(S^{1})}+\|\partial_{t}\psi(t)\|_{H^{1/2}(S^{1})}\|\partial_{s}\phi(t)\|_{L^{2}(S^{1})}^{1/2}\|\psi(t)\|_{L^{4}(S^{1})})\notag\\
&+C\int_{S^{1}}(1+|\psi(t)|^{p})|\psi_{0}^{+}(t)|\,ds+C\|\partial_{s}\phi(t)\|_{L^{2}(S^{1})}\|\psi_{0}^{+}(t)\|_{L^{4}(S^{1})}\|\psi(t)\|_{L^{4}(S^{1})},
\end{align}
where in the last inequality, we have used $|\nabla_{\psi}H(s,\phi,\psi)|\le C(1+|\psi|^{p})$ which follows from (1.4).

On the other hand, by (3.17), we have
\begin{equation}
\|\partial_{s}\phi(t)\|_{L^{2}(S^{1})}\le C(1+\|\psi(t)\|_{H^{1/2}(S^{1})}\|\partial_{t}\psi(t)\|_{H^{1/2}(S^{1})})^{1/2}
\end{equation}
and
\begin{equation}
\|\psi(t)\|_{L^{4}(S^{1})}\le C\|\psi(t)\|_{L^{p+1}(S^{1})}\le C(1+\|\psi(t)\|_{H^{1/2}(S^{1})}\|\partial_{t}\psi(t)\|_{H^{1/2}(S^{1})})^{1/p+1},
\end{equation}
where by our assumption $p\ge 3$, we have used $\|\psi(t)\|_{L^{4}(S^{1})}\le C\|\psi(t)\|_{L^{p+1}(S^{1})}$ which follows from the H\"older's inequality.
Plugging (3.21) and (3.22) into (3.20) and noting the inequality $\|\psi_{0}^{+}\|_{L^{q}(S^{1})}\le C\|\psi\|_{L^{q}(S^{1})}$ which holds for $1<q<\infty$ (since $P_{0}^{+}$ is a pseudo differential operator of order $0$), we obtain
\begin{align}
\|\psi_{0}^{+}(t)\|^{2}_{H^{1/2}(S^{1})}&\le C\|\partial_{t}\psi(t)\|_{H^{1/2}(S^{1})}\|\psi(t)\|_{H^{1/2}(S^{1})}\notag\\
&~+C(1+\|\partial_{t}\psi(t)\|_{H^{1/2}(S^{1})}\|\psi(t)\|_{H^{1/2}(S^{1})})^{\frac{1}{4}+\frac{1}{p+1}}\|\partial_{t}\psi(t)\|_{H^{1/2}(S^{1})}\notag\\
&~+C\|\psi(t)\|_{H^{1/2}(S^{1})}+C\|\psi(t)\|_{L^{p+1}(S^{1})}^{p+1}\notag\\
&~+C(1+\|\partial_{t}\psi(t)\|_{H^{1/2}(S^{1})}\|\psi(t)\|_{H^{1/2}(S^{1})})^{\frac{1}{2}+\frac{2}{p+1}}.
\end{align}
Here, $\frac{1}{2}+\frac{2}{p+1}\le 1$ for $p\ge3$ implies that 
$$(1+\|\partial_{t}\psi(t)\|_{H^{1/2}(S^{1})}\|\psi(t)\|_{H^{1/2}(S^{1})})^{\frac{1}{2}+\frac{2}{p+1}}\le 1+\|\partial_{s}\psi(t)\|_{H^{1/2}(S^{1})}\|\psi(t)\|_{H^{1/2}(S^{1})}.$$
Also, by (3.17), we have $\|\psi(t)\|_{L^{p+1}(S^{1})}^{p+1}\le C(1+\|\partial_{t}\psi(t)\|_{H^{1/2}(S^{1})}\|\psi(t)\|_{H^{1/2}(S^{1})})$. Thus, we have
\begin{align}
\|\psi_{0}^{+}(t)\|^{2}_{H^{1/2}(S^{1})}&\le C\|\partial_{t}\psi(t)\|_{H^{1/2}(S^{1})}\|\psi(t)\|_{H^{1/2}(S^{1})}+C\|\psi(t)\|_{H^{1/2}(S^{1})}+C\notag\\
&~+C(1+\|\partial_{t}\psi(t)\|_{H^{1/2}(S^{1})}\|\psi(t)\|_{H^{1/2}(S^{1})})^{\frac{1}{4}+\frac{1}{p+1}}\|\partial_{t}\psi(t)\|_{H^{1/2}(S^{1})}\notag\\
&\le C\|\partial_{t}\psi(t)\|_{H^{1/2}(S^{1})}\|\psi(t)\|_{H^{1/2}(S^{1})}+C\|\psi(t)\|_{H^{1/2}(S^{1})}+C\notag\\
&~+C(1+\|\partial_{t}\psi(t)\|_{H^{1/2}(S^{1})}\|\psi(t)\|_{H^{1/2}(S^{1})})^{\frac{1}{2}}\|\partial_{t}\psi(t)\|_{H^{1/2}(S^{1})},
\end{align}
where in the last inequality, we have used $\frac{1}{4}+\frac{1}{p+1}\le\frac{1}{2}$ for $p\ge 3$.

Integrating (3.24) over $[\tau,\tau+1]$ with respect to $t$, we have, by the H\"older's inequality
\begin{align}
&\int_{\tau}^{\tau+1}\|\psi_{0}^{+}(t)\|_{H^{1/2}(S^{1})}^{2}\,dt\le C\biggl(\int_{\tau}^{\tau+1}\|\partial_{t}\psi(t)\|_{H^{1/2}(S^{1})}^{2}\,dt\biggr)^{1/2}\biggl(\int_{\tau}^{\tau+1}\|\psi(t)\|_{H^{1/2}(S^{1})}^{2}\,dt\biggr)^{1/2}\notag\\
&~+C\biggl(\int_{\tau}^{\tau+1}\|\psi(t)\|^{2}_{H^{1/2}(S^{1})}\,dt\biggr)^{1/2}+C\notag\\
&+~C\biggl(\int_{\tau}^{\tau+1}(1+\|\partial_{t}\psi(t)\|_{H^{1/2}(S^{1})}\|\psi(t)\|_{H^{1/2}(S^{1})})\,dt\biggr)^{1/2}\biggl(\int_{\tau}^{\tau+1}\|\partial_{t}\psi(t)\|_{H^{1/2}(S^{1})}^{2}\,dt\biggr)^{1/2}\notag\\
&\le C+C\biggl(\int_{\tau}^{\tau+1}\|\psi(t)\|_{H^{1/2}(S^{1})}^{2}\,dt\biggr)^{1/2}+C\biggl(\int_{\tau}^{\tau+1}\|\psi(t)\|_{H^{1/2}(S^{1})}^{2}\,dt\biggr)^{1/4}\notag\\
&\le C+C\biggl(\int_{\tau}^{\tau+1}\|\psi(t)\|_{H^{1/2}(S^{1})}^{2}\,dt\biggr)^{1/2},
\end{align}
where in the second inequality, we have used $\int_{-\infty}^{\infty}\|\partial_{t}\psi(t)\|_{H^{1/2}(S^{1})}^{2}\,dt\le 2C_{0}$ (see (3.15)) and the H\"older's inequality once again.

By the same argument, we have
\begin{equation}
\int_{\tau}^{\tau+1}\|\psi_{0}^{-}(t)\|_{H^{1/2}(S^{1})}^{2}\,dt\le C\biggl(1+\biggl(\int_{\tau}^{\tau+1}|\psi(t)\|_{H^{1/2}(S^{1})}^{2}\,dt\biggr)^{\frac{1}{2}}\biggr).
\end{equation}
As for the component $\psi_{0}^{0}$, since the $H^{1/2}$-norm is equivalent to the $L^{p+1}$-norm on the finite dimensional space $\mathsf{H}^{0}_{0}$, we have
\begin{align}
\|\psi_{0}^{0}(t)\|_{H^{1/2}(S^{1})}^{2}&\le C\|\psi_{0}^{0}(t)\|_{L^{p+1}(S^{1})}^{2}\notag\\
&\le C(1+\|\partial_{t}\psi(t)\|_{H^{1/2}(S^{1})}\|\psi(t)\|_{H^{1/2}(S^{1})})^{\frac{2}{p+1}}\notag\\
&\le C(1+\|\partial_{t}\psi(t)\|_{H^{1/2}(S^{1})}\|\psi(t)\|_{H^{1/2}(S^{1})}),
\end{align}
where we have used (3.17) and $\frac{2}{p+1}<1$.

Integrating (3.27) over $[\tau,\tau+1]$ and using H\"older's inequality and (3.15), we have
\begin{equation}
\int_{\tau}^{\tau+1}\|\psi_{0}^{0}(t)\|_{H^{1/2}(S^{1})}^{2}\,dt\le C\biggl(1+\biggl(\int_{\tau}^{\tau+1}\|\psi(t)\|^{2}_{H^{1/2}(S^{1})}\,dt\biggr)^{\frac{1}{2}}\biggr).
\end{equation}

By (3.25), (3.26) and (3.28), we have
\begin{equation}
\int_{\tau}^{\tau+1}\|\psi(t)\|_{H^{1/2}(S^{1})}^{2}\,dt\le C\biggr(1+\biggr(\int_{\tau}^{\tau+1}\|\psi(t)\|_{H^{1/2}(S^{1})}^{2}\,dt\biggr)^{\frac{1}{2}}\biggr).
\end{equation}
(3.29) implies that
\begin{equation}
\int_{\tau}^{\tau+1}\|\psi(t)\|_{H^{1/2}(S^{1})}^{2}\,dt\le C
\end{equation}
for some $C>0$ depending only on $C_{0}$ and $N$.

On the other hand, integrating (3.17) over $[\tau,\tau+1]$ and using (3.15), (3.30) and the H\"older's inequality, we obtain
\begin{align}
&\int_{\tau}^{\tau+1}\|\phi(t)\|_{H^{1}(S^{1})}^{2}\,dt+\int_{\tau}^{\tau+1}\|\psi(t)\|_{L^{p+1}(S^{1})}^{p+1}\,dt\notag\\
&\le C+C\biggl(\int_{\tau}^{\tau+1}\|\partial_{t}\psi(t)\|_{H^{1/2}(S^{1})}^{2}\,dt\biggr)^{\frac{1}{2}}\biggl(\int_{\tau}^{\tau+1}\|\psi(t)\|_{H^{1/2}(S^{1})}^{2}\,dt\biggr)^{\frac{1}{2}}\notag\\
&\le C,
\end{align}
where $C>0$ depends only on $C_{0}$ and $N$. This completes the proof.
\hfill$\Box$

\medskip

We now prove Proposition 3.2.

\medskip

{\it Proof of Proposition 3.2.}
Recall that the negative gradient flow takes the form (3.7), (3.8). 
As before, we identify $\phi$ and $\psi$ with $\iota\circ\phi$ and $\tilde{\psi}=\iota_{\ast}\psi$ and denote them by $\phi$ and $\psi$, respectively. Thus, we think of $\phi$ as taking its values in $\R^{k}$ and $\psi$ takes its values in $\Spin(S^{1})\otimes\underline{\R}^{k}$. We first treat the spinor part of the equation (3.8). Recall, by (3.11), the second equation (3.8) is written as follows
\begin{align}
\partial_{t}\psi(t)&=-(1+|\D|)^{-1}(\D\psi-A(\partial_{s}\phi(t),\partial_{s}\cdot\psi(t))-\nabla_{\psi}H(s,\phi(t),\psi(t)))\notag\\
&=-(1+|\D|)^{-1}((\D+\lambda)\psi(t)-(\lambda\psi(t)+A(\partial_{s}\phi(t),\partial_{s}\cdot\psi(t))-\nabla_{\psi}H(s,\phi(t),\psi(t)))\notag\\
&=-\mathsf{L}_{\lambda}\psi(t)+(1+|\D|)^{-1}(\lambda\psi(t)+A(\partial_{s}\phi(t),\partial_{s}\cdot\psi(t))+\nabla_{\psi}H(s,\phi(t),\psi(t))),
\end{align}
where for $\lambda\in\R$, we set $\mathsf{L}_{\lambda}=(1+|\D|)^{-1}\D_{-\lambda}$ and $\D_{-\lambda}=\D+\lambda$.

We take $\lambda\in\R$ such that $-\lambda\not\in\Spec(\D)$. The fundamental solution of the differential operator $\frac{d}{dt}+\mathsf{L}_{\lambda}:C^{1}(\R,H^{1/2}(S^{1},\Spin(S^{1})\otimes\underline{\R}^{k}))\to C^{0}(\R,H^{1/2}(S^{1},\Spin(S^{1})\otimes\underline{\R}^{k})$ is given by
\begin{equation}
G_{\lambda}(t)=-e^{-t\mathsf{L}_{\lambda}^{-}}\bm{1}_{(-\infty,0]}(t)P_{\lambda}^{-}+e^{-t\mathsf{L}_{\lambda}^{+}}\bm{1}_{[0,\infty)}(t)P_{\lambda}^{+},
\end{equation}
where as in (3.10), $P_{\lambda}^{\pm}:H^{1/2}(S^{1},\Spin(S^{1})\otimes\underline{\R}^{k})\to\mathsf{H}_{\lambda}^{\pm}$ are spectral projections with respect to the operator $\D_{\lambda}$ onto its positive/negative subspaces of $H^{1/2}(S^{1},\Spin(S^{1})\otimes\underline{\R}^{k})=\mathsf{H}_{\lambda}^{-}\oplus\mathsf{H}^{+}_{\lambda}$ and $\mathsf{L}_{\lambda}^{\pm}=\mathsf{L}_{\lambda}|_{\mathsf{H}^{\pm}_{\lambda}}$.
Since $-\lambda\not\in\Spec(\D)$, there exists $\kappa>0$ such that the following holds:
\begin{equation}
\|G_{\lambda}(t)\|_{\mathsf{op}(H^{1/2}(S^{1},\Spin(S^{1})\otimes\underline{\R}^{k})}\le e^{-\kappa|t|},
\end{equation}
where for $T:H^{1/2}(S^{1},\Spin(S^{1})\otimes\underline{\R}^{k})\to H^{1/2}(S^{1},\Spin(S^{1})\otimes\underline{\R}^{k})$, $\|T\|_{\mathsf{op}(H^{1/2}(S^{1},\Spin(S^{1})\otimes\underline{\R}^{k}))}$ denotes the operator norm of $T$.

Using the fundamental solution $G_{\lambda}$, the solution $\psi(t)$ of equation (3.32) is expressed as
\begin{equation}
\psi(t)=\int_{\R}G_{\lambda}(\tau)(1+|\D|)^{-1}(\lambda\psi(t-\tau)+A(\partial_{s}\phi(t-\tau),\partial_{s}\cdot\psi(t-\tau))+\nabla_{\psi}H(s,\phi(t-\tau),\psi(t-\tau)))\,d\tau.
\end{equation}
From this, we have
\begin{align}
\|\psi(t)\|_{H^{1/2}(S^{1})}&\le\int_{\R}e^{-\kappa|\tau|}\|(1+|\D|)^{-1}(\lambda\psi(t-\tau)+A(\partial_{s}\phi(t-\tau),\partial_{s}\cdot\psi(t-\tau))\notag\\
&\qquad\qquad+\nabla_{\psi}H(s,\phi(t-\tau),\psi(t-\tau)))\|_{H^{1/2}(S^{1})}\,d\tau\notag\\
&=\sum_{n\in\Z}\int_{n}^{n+1}e^{-\kappa|\tau|}\|(1+|\D|)^{-1}(\lambda\psi(t-\tau)+A(\partial_{s}\phi(t-\tau),\partial_{s}\cdot\psi(t-\tau))\notag\\
&\qquad\qquad+\nabla_{\psi}H(s,\phi(t-\tau),\psi(t-\tau)))\|_{H^{1/2}(S^{1})}\,d\tau\notag\\
&\le\Big(\sum_{n\in\Z}e^{-\kappa|n|}\Big)\sup_{\eta\in\R}\int_{\eta}^{\eta+1}\|(1+|\D|)^{-1}(\lambda\psi(t-\tau)+A(\partial_{s}\phi(t-\tau),\partial_{s}\cdot\psi(t-\tau))\notag\\
&\qquad\qquad+\nabla_{\psi}H(s,\phi(t-\tau),\psi(t-\tau)))\|_{H^{1/2}(S^{1})}\,d\tau\notag\\
&\le C \sup_{\eta\in\R}\int_{\eta}^{\eta+1}\|(1+|\D|)^{-1}(\lambda\psi(t-\tau)+A(\partial_{s}\phi(t-\tau),\partial_{s}\cdot\psi(t-\tau))\notag\\
&\qquad\qquad+\nabla_{\psi}H(s,\phi(t-\tau),\psi(t-\tau)))\|_{H^{1/2}(S^{1})}\,d\tau.
\end{align}
Here, we note that the following hold:
\begin{equation}
|A(\partial_{s}\phi(t),\partial_{s}\cdot\psi(t))|\le C|\partial_{s}\phi(t)||\psi(t)|,
\end{equation}
\begin{equation}
|\nabla_{\psi}H(s,\phi(t),\psi(t))|\le C(1+|\psi(t)|^{p}).
\end{equation}
From (3.37) and the Sobolev embedding $H^{1/2}(S^{1})\subset L^{q}(S^{1})$ for any $1<q<+\infty$, we have $A(\partial_{s}\phi(t),\partial_{s}\cdot\psi(t))\in L^{r}(S^{1})$ for any $1<r<2$ and $t\in\R$ and
\begin{equation}
(1+|\D|)^{-1}A(\partial_{s}\phi(t),\partial_{s}\cdot\psi(t))\in W^{1,r}(S^{1})
\end{equation}
for any $t\in\R$ by the elliptic regularity. By the Sobolev embedding $W^{1,r}(S^{1})\subset H^{1/2}(S^{1})$ for $r\ge 1$ and $H^{1/2}(S^{1})\subset L^{\frac{2r}{2-r}}(S^{1})$, we therefore have
\begin{align}
\|(1+|\D|)^{-1}A(\partial_{s}\phi(t),\partial_{s}\cdot\psi(t))\|_{H^{1/2}(S^{1})}&\le C\|(1+|\D|)^{-1}A(\partial_{s}\phi(t),\partial_{s}\cdot\psi(t))\|_{W^{1,r}(S^{1})}\notag\\
&\le C\|A(\partial_{s}\phi(t),\partial_{s}\cdot\psi(t))\|_{L^{r}(S^{1})}\notag\\
&\le C\|\partial_{s}\phi(t)\|_{L^{2}(S^{1})}\|\psi(t)\|_{H^{1/2}(S^{1})}.
\end{align}
Similarly, by (3.38), we have $\nabla_{\psi}H(s,\phi(t),\psi(t))\in L^{\frac{p+1}{p}}(S^{1})$ and
\begin{align}
\|(1+|\D|)^{-1}\nabla_{\psi}H(s,\phi(t),\psi(t))\|_{H^{1/2}(S^{1})}&\le C\|(1+|\D|)^{-1}\nabla_{\psi}H(s,\phi(t),\psi(t))\|_{W^{1,\frac{p+1}{p}}(S^{1})}\notag\\
&\le C\|\nabla_{\psi}H(s,\phi(t),\psi(t))\|_{L^{\frac{p+1}{p}}(S^{1})}\notag\\
&\le C(1+\|\psi(t)\|_{L^{p+1}(S^{1})}^{p}).
\end{align}
Combining (3.36), (3.40) and (3.41), we obtain, by the H\"older's inequality
\begin{align}
&\|\psi(t)\|_{H^{1/2}(S^{1})}\notag\\
&\le C\sup_{\eta\in\R}\int_{\eta}^{\eta+1}(1+\|\psi(\tau)\|_{H^{1/2}(S^{1})}+\|\partial_{s}\phi(t)\|_{L^{2}(S^{1})}\|\psi(t)\|_{H^{1/2}(S^{1})}+\|\psi(t)\|_{L^{p+1}(S^{1})}^{p})\,d\tau\notag\\
&\le C\sup_{\eta\in\R}\biggl(1+\biggl(\int_{\tau}^{\tau+1}\|\psi(\tau)\|_{H^{1/2}(S^{1})}^{2}\,d\tau\biggr)^{\frac{1}{2}}\notag\\
&\hskip 2cm+\biggl(\int_{\eta}^{\eta+1}\|\partial_{s}\phi(\tau)\|_{L^{2}(S^{1})}^{2}\,d\tau\biggr)^{\frac{1}{2}}\biggl(\int_{\eta}^{\eta+1}\|\psi(\tau)\|_{H^{1/2}(S^{1})}^{2}\,d\tau\biggr)^{\frac{1}{2}}\notag\\
&\hskip 2cm+\biggl(\int_{\eta}^{\eta+1}\|\psi(\tau)\|_{L^{p+1}(S^{1})}^{p+1}\,d\tau\biggr)^{\frac{p}{p+1}}\biggr).
\end{align}
By (3.12)--(3.14), (3.42) implies that
\begin{equation}
\sup_{t\in\R}\|\psi(t)\|_{H^{1/2}(S^{1})}\le C(C_{0})
\end{equation}
for some $C(C_{0})>0$ depending only on $C_{0}$ and $N$.
This gives the desired estimate for $\psi(t)$.
To estimate $\sup_{t\in\R}\|\partial_{s}\phi(t)\|_{L^{2}(S^{1})}$, by (3.17) and (3.43), we first have
\begin{equation}
\|\partial_{s}\phi(t)\|_{L^{2}(S^{1})}^{2}\le C(1+\|\partial_{t}\psi(t)\|_{H^{1/2}(S^{1})}).
\end{equation}
On the other hand, by the equation (3.32), we have
\begin{align}
\|\partial_{t}\psi(t)\|_{H^{1/2}(S^{1})}&\le C(\|\psi(t)\|_{H^{1/2}(S^{1})}+\|(1+|\D|)^{-1}A(\partial_{s}\phi(t),\partial_{s}\cdot\psi(t))\|_{H^{1/2}(S^{1})}\notag\\
&\qquad+\|(1+|\D|)^{-1}\nabla_{\psi}H(s,\phi(t),\psi(t))\|_{H^{1/2}(S^{1})})\notag\\
&\le C(1+\|\psi(t)\|_{H^{1/2}(S^{1})}+\|\partial_{s}\phi(t)\|_{L^{2}(S^{1})}\|\psi(t)\|_{H^{1/2}(S^{1})}+\|\psi(t)\|_{L^{p+1}(S^{1})}^{p}),
\end{align}
where we have used (3.40) and (3.41).

By the embedding $H^{1/2}(S^{1})\subset L^{p+1}(S^{1})$ and (3.43), we have $\sup_{t\in\R}\|\psi(t)\|_{L^{p+1}(S^{1})}\le C(C_{0})$. Combining this with (3.45), we obtain
\begin{equation}
\|\partial_{t}\psi(t)\|_{H^{1/2}(S^{1})}\le C(1+\|\partial_{s}\phi(t)\|_{L^{2}(S^{1})}).
\end{equation}
By (3.44) and (3.46), we obtain
\begin{equation}
\|\partial_{s}\phi(t)\|_{L^{2}(S^{1})}^{2}\le C(1+\|\partial_{s}\phi(t)\|_{L^{2}(S^{1})}).
\end{equation}
(3.47) implies that
\begin{equation}
\|\partial_{s}\phi(t)\|_{L^{2}(S^{1})}^{2}\le C(C_{0})
\end{equation}
for all $t\in\R$, where $C(C_{0})>0$ is a constant depending only on $C_{0}$ and $N$.
This gives the desired estimate for $\phi(t)$ and completes the proof of Proposition 3.2.
\hfill$\Box$

\bigskip

Under the same assumption of Proposition 3.2, we can further refine the result of Proposition 3.2 by using elliptic regularity theory. In fact,  we have the following refinement:

\begin{proposition}
Let $(\phi(t),\psi(t))\in C^{1}(\R,\mathcal{F}^{1,1/2}(S^{1},N))$ be a solution to the negative gradient flow equation (3.9) which satisfies the same condition as in Proposition 3.2, $C_{0}=\sup_{t\in\R}|\Lag_{H}(\phi(t),\psi(t))|<+\infty$. Then there exists $C(C_{0})>0$ such that the following holds
$$\sup_{t\in\R}\|\phi(t)\|_{C^{2,2/3}(S^{1})}+\sup_{t\in\R}\|\psi(t)\|_{C^{1,2/3}(S^{1})}\le C(C_{0}).$$
\end{proposition}
\textit{Proof.} 
In the course of the proof, $C(C_{0})$ will denote constants which will depend only on $C_{0}$, but may change from line to line. We first consider the spinorial part of the flow equation. By the Sobolev embedding $H^{1/2}(S^{1})\subset L^{r}(S^{1})$ for all $r>1$ and by the result of Proposition 3.2, we have
\begin{equation}
\sup_{t\in\R}\|\lambda\psi(t)+A(\partial_{s}\phi(t),\partial_{s}\cdot\psi(t))+\nabla_{\psi}H(s,\phi(t),\psi(t))\|_{L^{r}(S^{1})}\le C(C_{0})
\end{equation}
for all $1<r<2$.
\bigskip
Thus, by the elliptic regularity theory, we have
\begin{equation}
\sup_{t\in\R}\|(1+|\D|)^{-1}(\lambda\psi(t)+A(\partial_{s}\phi(t),\partial_{s}\cdot\psi(t))+\nabla_{\psi}H(s,\phi(t),\psi(t)))\|_{W^{1,r}(S^{1})}\le C(C_{0}).
\end{equation}
Since $1<r<2$, by the Sobolev embedding theorem, we have $W^{1,r}(S^{1})\subset C^{0,\alpha}(S^{1})$ for some $0<\alpha<1$
and
\begin{equation}
\sup_{t\in\R}\|(1+|\D|)^{-1}(\lambda\psi(t)+A(\partial_{s}\phi(t),\partial_{s}\cdot\psi(t))+\nabla_{\psi}H(s,\phi(t),\psi(t)))\|_{C^{0,\alpha}(S^{1})}\le C(C_{0})
\end{equation}
by (3.50).

We recall the formula (3.35). By the mapping property $G_{\lambda}(t):C^{0,\alpha}(S^{1})\to C^{0,\alpha}(S^{1})$ and the estimate $\|G_{\lambda}(t)\|_{\mathsf{op}(C^{0,\alpha}(S^{1})}\le Ce^{-\kappa|t|}$ for all $t\in\R$, where $\kappa>0$, we have
\begin{align}
&\|\psi(t)\|_{C^{0,\alpha}(S^{1})}\notag\\
\le&~ C\int_{\R}e^{-\kappa|\tau|}\|(1+|\D|)^{-1}((1+|\D|)^{-1}(\lambda\psi(t)+A(\partial_{s}\phi(t),\partial_{s}\cdot\psi(t))+\nabla_{\psi}H(s,\phi(t),\psi(t)))\|_{C^{0,\alpha}(S^{1})}\notag\\
\le&~ C(C_{0})\int_{\R}e^{-\kappa|\tau|}\,d\tau\le C(C_{0})
\end{align}
for all $t\in\R$.

We next look at the $\phi$ part of the flow equation (3.7). We first observe that, in terms of the local coordinated of $N$, $\nabla_{\partial_{s}}\partial_{s}\phi=\partial_{s}^{2}\phi+\partial_{s}\phi^{i}\partial_{s}\phi^{j}\Gamma^{k}_{ij}(\phi)\frac{\partial}{\partial y^{k}}$. For simplicity, we write $\partial_{s}\phi^{i}\partial_{s}\phi^{j}\Gamma(\phi,\partial_{s}\phi,\partial_{s}\phi)=\Gamma^{k}_{ij}(\phi)\frac{\partial}{\partial y^{k}}$. Thus the equation (3.7) takes the form
\begin{align}
\partial_{s}\phi(t)&=(-\Delta+1)^{-1}\Big(\Delta\phi(t)+\Gamma(\phi(t),\partial_{s}\phi(t),\partial_{s}\phi(t))-
-\frac{1}{2}R(\phi)\langle\psi,\partial_{s}\cdot\psi\rangle+\nabla_{\phi}H(s,\phi,\psi)\Big)\notag\\
&=-\phi(t)+(-\Delta+1)^{-1}\Big(\phi(t)+\Gamma(\phi(t),\partial_{s}\phi(t),\partial_{s}\phi(t))-
-\frac{1}{2}R(\phi)\langle\psi,\partial_{s}\cdot\psi\rangle\notag\\
&\qquad+\nabla_{\phi}H(s,\phi,\psi)\Big).
\end{align}
Since the fundamental solution of the operator $\frac{d}{dt}+1$ is $G(t)=e^{-t}{\1}_{[0,\infty)}(t)$, we have the representation
\begin{align}
\phi(t)&=\int_{0}^{\infty}e^{-\tau}(-\Delta+1)^{-1}\Big(\phi(t-\tau)+\Gamma(\phi(t-\tau),\partial_{s}\phi(t-\tau),\partial_{s}\phi(t-\tau))-
\notag\\
&-\frac{1}{2}R(\phi(t-\tau))\langle\psi(t-\tau),\partial_{s}\phi(t-\tau)\cdot\psi(t-\tau)\rangle+\nabla_{\phi}H(s,\phi(t-\tau),\psi(t-\tau))\Big)\,d\tau.
\end{align}
We note that $|\Gamma(\phi,\partial_{s}\phi,\partial_{s}\phi)|\le C|\partial_{s}\phi|^{2}$ for some constant $C>0$ depending only on $N$. By the embedding $L^{1}(S^{1})\subset H^{-s}(S^{1})$ for any $s>\frac{1}{2}$ and Proposition 3.2, we thus have
\begin{equation}
\|(-\Delta+1)^{-1}\Gamma(\phi(t),\partial_{s}\phi(t),\partial_{s}\phi(t))\|_{W^{2-s,2}(S^{1})}\le C\|\Gamma(\phi(t),\partial_{s}\phi(t),\partial_{s}\phi(t))\|_{L^{1}(S^{1})}\le C(C_{0})
\end{equation}
for all $t\in\R$.
We recall that the term $R(\phi)\langle\psi,\partial_{s}\phi\cdot\psi\rangle$ is given by $R(\phi)\langle\psi,\partial_{s}\phi\cdot\psi\rangle=\Big\langle\psi,\partial_{s}\cdot\psi^{i}\otimes\frac{\partial}{\partial y^{j}}(\phi)\Big\rangle\partial_{s}\phi^{l}R^{j}_{iml}(\phi)g^{ms}(\phi)\frac{\partial}{\partial y^{s}}(\phi)$ and it belongs to $L^{r}(S^{1})$ for any $1<r<2$. Thus, by Proposition 3.2, we have the following bound for any $1<r<2$ and for any $t\in\R$:
\begin{equation}
\|(-\Delta+1)^{-1}R(\phi(t))\langle\psi(t),\partial_{s}\phi(t)\cdot\psi(t)\rangle\|_{W^{2,r}(S^{1})}\le C(C_{0}).
\end{equation}
By (1.5) and the Sobolev embedding $H^{1/2}(S^{1})\subset L^{r}(S^{1})$ (for any $1<r<\infty$), we also have $\nabla_{\phi}H(\cdot,\phi,\psi)\in L^{r}(S^{1})$ for any $1<r<\infty$ and
\begin{equation}
\|(-\Delta+1)^{-1}\nabla_{\phi}H(\cdot,\phi(t),\psi(t))\|_{W^{2,r}(S^{1})}\le C(C_{0})
\end{equation}
for all $t\in\R$.
From (3.54)--(5.57), we obtain
\begin{equation}
\sup_{t\in\R}\|\phi(t)\|_{W^{2-s,2}(S^{1})}\le C(C_{0}),
\end{equation}
where $\frac{1}{2}<s<1$ is arbitrary.

If we take $s=\frac{2}{3}$ in (3.59), we obtain $\sup_{t\in\R}\|\partial_{s}\phi(t)\|_{W^{\frac{1}{3},2}(S^{1})}\le C(C_{0})$. Combining this with the embedding $W^{\frac{1}{3},2}(S^{1})\subset L^{6}(S^{1})$, we have
\begin{equation}
\sup_{t\in\R}\|\partial_{s}\phi(t)\|_{L^{6}(S^{1})}\le C(C_{0}).
\end{equation}
(3.52) and (3.60) improve the bounds (3.49) and (3.50) as follows:
\begin{equation}
\sup_{t\in\R}\|\lambda\psi(t)+A(\partial_{s}\phi(t),\partial_{s}\cdot\psi(t))+\nabla_{\psi}H(s,\phi(t),\psi(t))\|_{L^{6}(S^{1})}\le C(C_{0}),
\end{equation}
\begin{equation}
\sup_{t\in\R}\|(1+|\D|)^{-1}(\lambda\psi(t)+A(\partial_{s}\phi(t),\partial_{s}\cdot\psi(t))+\nabla_{\psi}H(s,\phi(t),\psi(t)))\|_{W^{1,6}(S^{1})}\le C(C_{0}).
\end{equation}
From the representation (3.35), (3.62) and the mapping property $G_{\lambda}(t):W^{1,6}(S^{1})\to W^{1,6}(S^{1})$ with the bound $\|G_{\lambda}(t)\|_{\mathsf{op}(W^{1,6}(S^{1})}\le Ce^{-|\kappa|t}$, we have, as in (3.52)
\begin{equation}
\sup_{t\in\R}\|\psi(t)\|_{W^{1,6}(S^{1})}\le C(C_{0}).
\end{equation}
By (3.58), (3.59) and (3.62), we have the following improvement of the estimates (3.55)--(3.57):
\begin{equation}
\|(-\Delta+1)^{-1}\Gamma(\phi(t),\partial_{s}\phi(t),\partial_{s}\phi(t))\|_{W^{2,3}(S^{1})}\le C(C_{0}),
\end{equation}
\begin{equation}
\|(-\Delta+1)^{-1}R(\phi(t))\langle\psi(t),\partial_{s}\phi(t)\cdot\psi(t)\rangle\|_{W^{2,6}(S^{1})}\le C(C_{0})
\end{equation}
and
\begin{equation}
\|(-\Delta+1)^{-1}\nabla_{\phi}H(\cdot,\phi(t),\psi(t))\|_{W^{2,r}(S^{1})}\le C(C_{0})
\end{equation}
for any $1<r<\infty$. Thus, by (3.54), (3.63)--(3.65), we have the following improvement of (3.58):
\begin{equation}
\sup_{t\in\R}\|\phi(t)\|_{W^{2,3}(S^{1})}\le C(C_{0}).
\end{equation}
By the Sobolev embedding $W^{2,3}(S^{1})\subset C^{1,\frac{2}{3}}(S^{1})$, we have 
\begin{equation}
\sup_{t\in\R}\|\partial_{s}\phi(t)\|_{C^{0,\frac{2}{3}}(S^{1})}\le C(C_{0}).
\end{equation}
Once again, we return to the equation of $\psi$, (3.32). By (3.62) and (3.67) and the Sobolev embedding $W^{1,6}(S^{1})\subset C^{0,\frac{2}{3}}(S^{1})$, we have
$$\sup_{t\in\R}\|\lambda\psi(t)+A(\partial_{s}\phi(t),\partial_{s}\cdot\psi(t))+\nabla_{\psi}H(s,\phi(t),\psi(t))\|_{C^{0,\frac{2}{3}}(S^{1})}\le C(C_{0})$$
and 
\begin{equation}
\sup_{t\in\R}\|(1+|\D|)^{-1}(\lambda\psi(t)+A(\partial_{s}\phi(t),\partial_{s}\cdot\psi(t))+\nabla_{\psi}H(s,\phi(t),\psi(t)))\|_{C^{1,\frac{2}{3}}(S^{1})}\le C(C_{0}).
\end{equation}
By (3.35), (3.68) and the mapping property $G_{\lambda}(t):C^{1,\frac{2}{3}}(S^{1})\to C^{1,\frac{2}{3}}(S^{1})$ with the bound $\|G_{\lambda}(t)\|_{\mathsf{op}(C^{1,\frac{2}{3}}(S^{1})}\le Ce^{-|\kappa|t}$ for any $t\in\R$, where $\kappa>0$, as in (3.52) we have
\begin{equation}
\sup_{t\in\R}\|\psi(t)\|_{C^{1,\frac{2}{3}}(S^{1})}\le C(C_{0}).
\end{equation}
This gives the desired estimate for $\psi$.

Now, we again return to the equation for $\phi$, (3.7).
By (3.67) and (3.69), we have
\begin{equation}
\|(-\Delta+1)^{-1}\Gamma(\phi(t),\partial_{s}\phi(t),\partial_{s}\phi(t))\|_{C^{2,\frac{2}{3}}(S^{1})}\le C(C_{0}),
\end{equation}
\begin{equation}
\|(-\Delta+1)^{-1}R(\phi(t))\langle\psi(t),\partial_{s}\phi(t)\cdot\psi(t)\rangle\|_{C^{2,\frac{2}{3}}(S^{1})}\le C(C_{0}),
\end{equation}
\begin{equation}
\|(-\Delta+1)^{-1}\nabla_{\phi}H(\cdot,\phi(t),\psi(t))\|_{C^{2,\frac{2}{3}}(S^{1})}\le C(C_{0}).
\end{equation}
From the representation (3.54), (3.70)--(3.72), we finally obtained the desired estimate
\begin{equation}
\sup_{t\in\R}\|\phi(t)\|_{C^{2,\frac{2}{3}}(S^{1})}\le C(C_{0}).
\end{equation}
This completes the proof.
\hfill$\Box$


\subsection{Fredholm property}

In this section, we first give a functional analytic set up for the moduli problem of gradient flow lines. We then study its Fredholm property.

Let $\mathsf{x}_{-}=(\phi_{-},\psi_{-}),\mathsf{x}_{+}=(\phi_{+},\psi_{+})\in\text{crit}(\Lag_{H})$. We define the $W^{1,2}$-Sobolev space of paths connecting $\mathsf{x}_{-}$ and $\mathsf{x}_{+}$ denoted by $W^{1,2}_{\mathsf{x}_{-},\mathsf{x}_{+}}(\R,\mathcal{F}^{1,1/2}(S^{1},N))$ as follows: Let us denote by $\iota(N)$ the injectivity radius of $N$.
$(\phi,\psi)\in W^{1,2}_{\mathsf{x}_{-},\mathsf{x}_{+}}(\R,\mathcal{F}^{1,1/2}(S^{1},N))$ if and only if the following (i) and (ii) hold: 

(i) $\phi\in W^{1,2}_{\text{loc}}(\R,H^{1}(S^{1},N))$ and there exists $T>0$ such that
\begin{equation}
\phi(t)=\exp_{\phi_{-}}(X_{-}(t))
\end{equation}
for $t\le -T$ for some $X_{-}\in W^{1,2}((-\infty,-T],H^{1}(S^{1},\phi_{-}^{\ast}TN))$ with $|X_{-}(t)(s)|_{T_{\phi_{-}(t)(s)}N}<\iota(N)$ for all $(t,s)\in(-\infty,-T]\times S^{1}$,
\begin{equation}
\phi(t)=\exp_{\phi_{+}}(X_{+}(t))
\end{equation}
for $t\ge T$ for some $X_{+}\in W^{1,2}([T,+\infty),H^{1}(S^{1},\phi_{+}^{\ast}TN))$ with $|X_{+}(t)(s)|_{T_{\phi_{+}(t)(s)}N}<\iota(N)$ for all $(t,s)\in[T,+\infty)\times S^{1}$.

(ii) $\psi\in W^{1,2}_{\text{loc}}(\R,H^{1/2}(S^{1},\Spin(S^{1})\otimes\phi^{\ast}TN)$ and
\begin{equation}
\psi(t)=\mathsf{S}_{-,t}(\psi_{-}+\xi_{-}(t))
\end{equation}
for $t\le -T$ for some $\xi_{-}\in W^{1,2}((-\infty,-T],H^{1/2}(S^{1},\Spin(S^{1})\otimes\phi_{-}^{\ast}TN))$ (where $T>0$ is as in (i)),
\begin{equation}
\psi(t)=\mathsf{S}_{+,t}(\psi_{+}+\xi_{+}(t))
\end{equation}
for $t\ge T$ for some $\xi_{+}\in W^{1,2}([T,+\infty),H^{1/2}(S^{1},\Spin(S^{1})\otimes\phi_{+}^{\ast}TN))$.

In (3.76) and (3.77), $\mathsf{S}_{\pm,t}:\phi_{\pm}^{\ast}TN\to\phi(t)^{\ast}TN$ denote the parallel translation along the path $[0,1]\ni\tau\mapsto\exp_{\phi_{\pm}}(\tau X_{\pm}(t))\in N$ for $\pm t\ge T$.

We then define $\mathcal{F}_{\mathsf{x}_{-},\mathsf{x}_{+}}:W^{1,2}_{\mathsf{x}_{-},\mathsf{x}_{+}}(\R,\mathcal{F}^{1,1/2}(S^{1},N))\to L^{2}(\R,T\mathcal{F}^{1,1/2}(S^{1},N))$ by
\begin{equation}
\mathcal{F}_{\mathsf{x}_{-},\mathsf{x}_{+}}(\ell)=\frac{d\ell}{dt}+\nabla_{1,1/2}\Lag_{H}(\ell)
\end{equation}
for $\ell\in W^{1,2}_{\mathsf{x}_{-},\mathsf{x}_{+}}(\R,\mathcal{F}^{1,1/2}(S^{1},N))$, where $L^{2}(\R,T\mathcal{F}^{1,1/2}(S^{1},N))$ is a fiber bundle defined over $W^{1,2}_{\mathsf{x}_{-},\mathsf{x}_{+}}(\R,\mathcal{F}^{1,1/2}(S^{1},N))$ whose fiber at $\ell\in W^{1,2}_{\mathsf{x}_{-},\mathsf{x}_{+}}(\R,\mathcal{F}^{1,1/2}(S^{1},N))$ is $L^{2}(\R,\ell^{\ast}T\mathcal{F}^{1,1/2}(S^{1},N))$, that is,
$$L^{2}(\R,T\mathcal{F}^{1,1/2}(S^{1},N))_{\ell}=\Big\{V(t)\in T_{\ell(t)}\mathcal{F}^{1,1/2}(S^{1},N),\,\int_{\R}\|V(t)\|_{T_{\ell(t)}\mathcal{F}^{1,1/2}(S^{1},N)}^{2}\,dt<+\infty\Big\}.$$
By our definition of the space $W^{1,2}_{\mathsf{x}_{-},\mathsf{x}_{+}}(\R,\mathcal{F}^{1,1/2}(S^{1},N))$, under the assumption (1.1)--(1.3) in the introduction, it is not difficult to prove that (3.78) is well-defined, i.e., $\frac{d\ell}{dt}+\nabla_{1,1/2}\Lag_{H}(\ell)\in L^{2}(\R,T\mathcal{F}^{1,1/2}(S^{1},N))$ for all $\ell\in W^{1,2}_{\mathsf{x}_{-},\mathsf{x}_{+}}(\R,\mathcal{F}^{1,1/2}(S^{1},N))$. Moreover, $\mathcal{F}_{\mathsf{x}_{-},\mathsf{x}_{+}}$ defines a $C^{1}$-map from $W^{1,2}_{\mathsf{x}_{-},\mathsf{x}_{+}}(\R,\mathcal{F}^{1,1/2}(S^{1},N))$ to $L^{2}(\R,T\mathcal{F}^{1,1/2}(S^{1},N))$. 

We shall prove the following:

\begin{proposition}
Assume that $\mathsf{x}_{-},\mathsf{x}_{+}\in\text{crit}(\Lag_{H})$ are non-degenerate. Then $\mathcal{F}_{\mathsf{x}_{-},\mathsf{x}_{+}}$ defined by (3.78) is Fredholm. Its index at any $\ell\in W^{1,2}_{\mathsf{x}_{-},\mathsf{x}_{+}}(\R,\mathcal{F}^{1,1/2}(S^{1},N))$ is given by ${\rm ind}(d\mathcal{F}_{\mathsf{x}_{-},\mathsf{x}_{+}}(\ell))=\mu_{H}(\mathsf{x}_{-})-\mu_{H}(\mathsf{x}_{+})$.
\end{proposition}
{\it Proof.} We decompose the proof into two steps:

\noindent
{\bf Step 1.} {\it Reduction to the case of $\ell$ which is constant near $t=\pm\infty$.} 

\noindent
We first reduce to the case where $\ell\in W^{1,2}_{\mathsf{x}_{-},\mathsf{x}_{+}}(\R,\mathcal{F}^{1,1/2}(S^{1},N))$ is equal to $\mathsf{x}_{\pm}$ for $\pm t$ large enough. That is, we shall show that it suffices to prove the assertion for the case where $\ell$ satisfies $\ell(t)=\mathsf{x}_{-}$ for $t\le-R$ and $\ell(t)=\mathsf{x}_{+}$ for $t\ge R$ for some $R>0$.

To prove this, let $\ell(t)=(\phi(t),\psi(t))\in W^{1,2}_{\mathsf{x}_{-},\mathsf{x}_{+}}(\R,\mathcal{F}^{1,1/2}(S^{1},N))$ be arbitrary. By definition, we have the representation (3.74)--(3.77). For $R>T$, we define $\phi_{R}$ and $\psi_{R}$ as follows:
$$\phi_{R}(t)=\left\{\begin{array}{lllll}
\phi_{-}&\ (t\le-R-1)\\
\exp_{\phi_{-}}((t+R+1)X_{-}(-R))&\ (-R-1\le t\le-R)\\
\phi(t)&\ (-R\le t\le R)\\
\exp_{\phi_{+}}((R+1-t)X_{+}(R))&\ (R\le t\le R+1)\\
\phi_{+}&\ (R+1\le t),
\end{array}\right.$$
$$\psi_{R}(t)=\left\{\begin{array}{lllll}
\psi_{-}&\ (t\le-R-1)\\
\mathsf{S}_{-,t;R}(\psi_{-}+(t+R+1)\xi_{-}(-R))&\ (-R-1\le t\le -R)\\
\psi(t)&\ (-R\le t\le R)\\
\mathsf{S}_{+,t;R}(\psi_{+}+(R+1-t)\xi_{+}(R))&\ (R\le t\le R+1)\\
\psi_{+}&\ (R+1\le t),
\end{array}\right.$$
where $\mathsf{S}_{\pm,t;R}$ are the parallel translations along paths $[0,1]\ni\tau\mapsto\exp_{\phi_{\pm}}(\tau(\mp t+R+1)X_{\pm}(\pm R))$.

It can be easily checked that $\ell_{R}=:(\phi_{R},\psi_{R})\in W^{1,2}_{\mathsf{x}_{-},\mathsf{x}_{+}}(\R,\mathcal{F}^{1,1/2}(S^{1},N))$ and $\ell_{R}$ converges to $\ell=(\phi,\psi)$ as $R\to\infty$ in the sense that the following holds
\begin{equation}
\|\phi_{R}-\phi\|_{W^{1,2}(\R,H^{1}(S^{1},\R^{k}))}\to 0,
\end{equation}
\begin{equation}
\|\psi_{R}-\psi\|_{W^{1,2}(\R,H^{1}(S^{1},\Spin(S^{1})\otimes\underline{\R^{k}}))}\to 0
\end{equation}
as $R\to\infty$.

We next show that $d\mathcal{F}_{\mathsf{x}_{-},\mathsf{x}_{+}}(\ell_{R})$ converges (in the sense of the operator norm) to $d\mathcal{F}_{\mathsf{x}_{-},\mathsf{x}_{+}}(\ell)$ as $R\to\infty$. For this purpose, we observe that by the condition $X_{-}\in W^{1,2}((-\infty,-T],H^{1}(S^{1},\phi_{-}^{\ast}TN))$, we have $\|X_{-}(t)\|_{H^{1}(S^{1})}\to 0$ as $t\to-\infty$. Similarly, we also have $\|X_{+}(t)\|_{H^{1}(S^{1})}\to 0$ as $t\to\infty$. Combining these with the definition of $\phi_{R}$, it follows that $\|\phi_{R}-\phi\|_{L^{\infty}(\R\times S^{1})}\to 0$ as $R\to\infty$. We take $R>0$ so large that $\|\phi_{R}-\phi\|_{L^{\infty}(\R\times S^{1})}<\iota(N)$. We denote by $\mathsf{S}_{R}(t):T_{\phi(t)}N\to T_{\phi_{R}(t)}N$ the parallel translation along the shortest geodesic starting from $\phi(t)$ and ending at $\phi_{R}(t)$. By~\cite[Lemma 7.4]{I}, $S_{R}(t)$ defines naturally a bounded linear operator which we also denote by $\mathsf{S}_{R}(t)$
\begin{equation}
S_{R}(t):H^{1}(S^{1},\phi(t)^{\ast}TN)\to H^{1}(S^{1},\phi_{R}(t)^{\ast}TN),
\end{equation}
\begin{equation}
S_{R}(t)=\bm{1}_{\Spin(S^{1})}\otimes S_{R}(t):H^{1/2}(S^{1},\Spin(S^{1})\otimes\phi(t)^{\ast}TN)\to H^{1/2}(S^{1},\Spin(S^{1})\otimes\phi_{R}(t)^{\ast}TN).
\end{equation}
We then have:
\begin{lemma} For $\Theta=(\xi,\zeta)\in W^{1,2}(\R,\ell^{\ast}T\mathcal{F}^{1,1/2}(S^{1},N))$, we define $(\mathsf{S}_{R}\Theta)(t)=(\mathsf{S}_{R}(t)\xi(t),\mathsf{S}_{R}(t)\zeta(t))$. Then $\mathsf{S}_{R}:W^{1,2}(\R,\ell^{\ast}T\mathcal{F}^{1,1/2}(S^{1},N))\to W^{1,2}(\R,\ell_{R}^{\ast}T\mathcal{F}^{1,1/2}(S^{1},N))$ is an isomorphism.
Moreover, we have the estimate
$$\|\mathsf{S}_{R}\Theta-\Theta\|_{W^{1,2}(\R,H^{1}(S^{1},\R^{k})\times H^{1/2}(S^{1},\Spin(S^{1})\otimes\underline{\R}^{k}))}\le\delta(R)\|\Theta\|_{W^{1,2}(\R,\ell^{\ast}T\mathcal{F}^{1,1/2}(S^{1},N))}$$
for some $\delta(R)>0$ with $\delta(R)\to 0$ as $R\to\infty$.
\end{lemma}
{\it Proof of Lemma 3.2.} By~\cite[Lemma 7.4]{I}, the operator norms $\|\mathsf{S}_{R}(t)\|_{\text{op}(H^{1}(S^{1},\phi(t)^{\ast}TN),H^{1}(S^{1},\phi_{R}(t)^{\ast}TN))}$ and $\|\mathsf{S}_{R}(t)\|_{\mathsf{op}(H^{1/2}(S^{1},\Spin(S^{1})\otimes\phi(t)^{\ast}TN),H^{1/2}(S^{1},\Spin(S^{1})\otimes\phi_{R}(t)^{\ast}TN))}$ depends only on $\|\phi(t)\|_{H^{1}(S^{1})}$ and $\|\phi_{R}(t)\|_{H^{1}(S^{1})}$. Since $\phi\in W^{1,2}(\R,H^{1}(S^{1},N))$, it follows that $\sup_{t\in\R}\|\phi(t)\|_{H^{1}(S^{1})}<+\infty$. Combining this with the definition of $\phi_{R}$, it is easy to see that we have a uniform bound $$\|\phi_{R}(t)\|_{H^{1}(S^{1})}\le C$$ for some $C>0$ independent of $t\in\R$ and $R>T$.
We thus have uniform bounds of the operator norms
\begin{equation}
\|\mathsf{S}_{R}(t)\|_{\text{op}(H^{1}(S^{1},\phi(t)^{\ast}TN), H^{1}(S^{1},\phi_{R}(t)^{\ast}TN))}\le C,
\end{equation}
\begin{equation}
\|\mathsf{S}_{R}(t)\|_{\mathsf{op}(H^{1/2}(S^{1},\Spin(S^{1})\otimes\phi(t)^{\ast}TN),H^{1/2}(S^{1},\Spin(S^{1})\otimes\phi_{R}(t)^{\ast}TN))}\le C,
\end{equation}
where $C>0$ does not depend on $t\in\R$ and $R>T$.

From (3.83) and (3.84), we see that
\begin{equation}
\mathsf{S}_{R}:L^{2}(\R,\ell^{\ast}T\mathcal{F}^{1,1/2}(S^{1},N))\to L^{2}(\R,\ell_{R}^{\ast}T\mathcal{F}^{1,1/2}(S^{1},N))
\end{equation}
defined by $(\mathsf{S}_{R}\Theta)(t)=(\mathsf{S}_{R}(t)\xi(t),\mathsf{S}_{R}(t)\zeta(t))$ is a bounded linear map.

We next estimate the $H^{1}(S^{1})$-norm of the derivative $\nabla_{t}(\mathsf{S}_{R}(t)\xi(t))=\nabla_{t}\mathsf{S}_{R}(t)\xi(t)+\mathsf{S}_{R}(t)(\nabla_{t}\xi(t))$.
For this purpose, recall that, in the notation of~\cite{I}, $\mathsf{S}_{R}(t)$ is defined by $\mathsf{S}_{R}(t)=\mathsf{P}_{\phi_{R}(t),\phi(t)}$, where for $x,y\in N$ with $d(x,y)<\iota(N)$, $\mathsf{P}_{x,y}:T_{x}N\to T_{y}N$ denotes the parallel translation along the shortest geodesic from $x$ to $y$. Thus
\begin{equation}
\nabla_{t}\mathsf{S}_{R}(t)=\nabla_{x}\mathsf{P}_{\phi_{R}(t),\phi(t)}[\partial_{t}\phi_{R}(t)]+\nabla_{y}\mathsf{P}_{\phi_{R}(t),\phi(t)}[\partial_{t}\phi(t)].
\end{equation}
By (3.79) and (3.86), we have
\begin{equation}
\nabla_{t}\mathsf{S}_{R}(t)\to\nabla_{x}\mathsf{P}_{\phi(t),\phi(t)}[\partial_{t}\phi(t)]+\nabla_{y}\mathsf{P}_{\phi(t),\phi(t)}[\partial_{t}\phi(t)]=0
\end{equation}
in $H^{1}(S^{1})$ as $R\to\infty$ since $\nabla_{x}\mathsf{P}_{x,x}+\nabla_{y}\mathsf{P}_{x,x}=0$. Moreover, the convergence is uniform with respect to $t$.

We define $\epsilon(R)=\sup_{t\in\R}\|\nabla_{t}\mathsf{S}_{R}(t)\|_{H^{1}(S^{1})}$. Note that $\epsilon(R)\to 0$ as $R\to\infty$ and
\begin{equation}
\|\nabla_{t}\mathsf{S}_{R}(t)\xi(t)\|_{H^{1}(S^{1})}\le C\|\nabla_{t}\mathsf{S}_{R}\|_{H^{1}(S^{1})}\|\xi(t)\|_{H^{1}(S^{1})}\le C\epsilon(R)\|\xi(t)\|_{H^{1}(S^{1})},
\end{equation}
where $C>0$ depends only on the constant of the Sobolev embedding $H^{1}(S^{1})\subset L^{\infty}(S^{1})$. (Thus, $H^{1}(S^{1})$ becomes an algebra with respect to the point wise multiplication. This fact has been used in the second inequality of (3.88)). Obviously, the similar estimate holds for $\zeta(t)\in H^{1}(S^{1},\Spin(S^{1})\otimes\phi^{\ast}(t)TN)$,
\begin{equation}
\|\nabla_{t}\mathsf{S}_{R}(t)\zeta(t)\|_{H^{1}(S^{1})}\le C\|\nabla_{t}\mathsf{S}_{R}\|_{H^{1}(S^{1})}\|\zeta(t)\|_{H^{1}(S^{1})}\le C\epsilon(R)\|\zeta(t)\|_{H^{1}(S^{1})}
\end{equation}
On the other hand, for $\zeta\in L^{2}(S^{1},\Spin(S^{1})\otimes\phi^{\ast}(t)TN)$, we have
\begin{align}
\|\nabla_{t}\mathsf{S}_{R}(t)\zeta(t)\|_{L^{2}(S^{1})}&\le\|\nabla_{t}\mathsf{S}_{R}(t)\|_{L^{\infty}(S^{1})}\|\zeta(t)\|_{L^{2}(S^{1})}\notag\\
&\le C\|\nabla_{t}\mathsf{S}_{R}(t)\|_{H^{1}(S^{1})}\|\zeta(t)\|_{L^{2}(S^{1})}\notag\\
&\le C\epsilon(R)\|\zeta(t)\|_{L^{2}(S^{1})},
\end{align}
where we have used the Sobolev embedding $H^{1}(S^{1})\subset L^{\infty}(S^{1})$ again. By interpolating (3.89) and (3.90), we have a similar estimate for the $H^{1/2}(S^{1})$-norm of $\nabla_{t}\mathsf{S}_{R}(t)\zeta(t)$:
\begin{equation}
\|\nabla_{t}\mathsf{S}_{R}(t)\zeta(t)\|_{H^{1/2}(S^{1})}\le C\epsilon(R)\|\zeta(t)\|_{H^{1/2}(S^{1})}.
\end{equation}
Combining (3.88), (3.91) with the boundedness of (3.85), we see that $\mathsf{S}_{R}:W^{1,2}(\R,\ell^{\ast}T\mathcal{F}^{1,1/2}(S^{1},N))\to W^{1,2}(\R,\ell_{R}^{\ast}T\mathcal{F}^{1,1/2}(S^{1},N))$ defines a bounded linear map. Since the inverse of $\mathsf{S}_{R}$ is defined by the parallel translation along the  shortest geodesic starting from $\phi_{R}(t)$ and ending at $\phi(t)$ and by the same reasoning as above it is bounded, the first assertion of the lemma is proved.

To prove the second assertion, we recall the estimates of~\cite[Lemma 7.4]{I}:
\begin{equation}
\|\mathsf{S}_{R}\xi(t)-\xi(t)\|_{H^{1}(S^{1})}\le C(\|\phi_{R}(t)\|_{H^{1}(S^{1})},\|\phi(t)\|_{H^{1}(S^{1})})\|\phi_{R}(t)-\phi(t)\|_{H^{1}(S^{1},\R^{k})}\|\xi(t)\|_{H^{1}(S^{1})},
\end{equation}
\begin{equation}
\|\mathsf{S}_{R}\zeta(t)-\zeta(t)\|_{H^{1/2}(S^{1})}\le C(\|\phi_{R}(t)\|_{H^{1}(S^{1})},\|\phi(t)\|_{H^{1}(S^{1})})\|\phi_{R}(t)-\phi(t)\|_{H^{1}(S^{1},\R^{k})}\|\zeta(t)\|_{H^{1/2}(S^{1})}.
\end{equation}
By (3.92) and (3.88), we have
\begin{align}
&\|\nabla_{t}(\mathsf{S}_{R}(t)\xi(t))-\nabla_{t}\xi(t)\|_{H^{1}(S^{1})}\le\|\nabla_{t}\mathsf{S}_{R}(t)\xi(t)\|_{H^{1}(S^{1})}+\|\mathsf{S}_{R}(t)\nabla_{t}\xi(t)-\nabla_{t}\xi(t)\|_{H^{1}(S^{1})}\notag\\
&\le C\epsilon(R)\|\xi(t)\|_{H^{1}(S^{1})}+C(\|\phi_{R}(t)\|_{H^{1}(S^{1})},\|\phi(t)\|_{H^{1}(S^{1})})\|\phi_{R}(t)-\phi(t)\|_{H^{1}(S^{1})}\|\nabla_{t}\xi(t)\|_{H^{1}(S^{1})}.
\end{align}
Similarly, by (3.93) and (3.91) we have
\begin{align}
&\|\nabla_{t}(\mathsf{S}_{R}(t)\zeta(t))-\nabla_{t}\zeta(t)\|_{H^{1/2}(S^{1})}\le\|\nabla_{t}\mathsf{S}_{R}(t)\zeta(t)\|_{H^{1/2}(S^{1})}+\|\mathsf{S}_{R}(t)\nabla_{t}\zeta(t)-\nabla_{t}\zeta(t)\|_{H^{1/2}(S^{1})}\notag\\
&\le C\epsilon(R)\|\zeta(t)\|_{H^{1/2}(S^{1})}+C(\|\phi_{R}(t)\|_{H^{1}(S^{1})},\|\phi(t)\|_{H^{1}(S^{1})})\|\phi_{R}(t)-\phi(t)\|_{H^{1}(S^{1})}\|\nabla_{t}\zeta(t)\|_{H^{1/2}(S^{1})}.
\end{align}
Since $\sup_{t\in\R}\|\phi_{R}(t)-\phi(t)\|_{H^{1}(S^{1})}\to 0$ as $R\to\infty$, the asserted estimate follows from (3.92)--(3.95). This completes the proof.
\hfill$\Box$

\begin{lemma} We have
$$\|\mathsf{S}_{R}^{-1}\circ d\mathcal{F}_{\mathsf{x}_{-},\mathsf{x}_{+}}(\ell_{R})\circ\mathsf{S}_{R}-d\mathcal{F}_{\mathsf{x}_{-},\mathsf{x}_{+}}(\ell)\|_{\mathsf{op}(W^{1,2}(\R,\ell^{\ast}T\mathcal{F}^{1,1/2}(S^{1},N)),L^{2}(\R,\ell^{\ast}T\mathcal{F}^{1,1/2}(S^{1},N)))}\to 0$$
as $R\to\infty$.
\end{lemma}
{\it Proof.}
For $\Theta=(\xi,\zeta)\in W^{1,2}(\R,\ell^{\ast}T\mathcal{F}^{1,1/2}(S^{1},N))$, we have
\begin{align}
\mathsf{S}_{R}^{-1}\circ d\mathcal{F}_{\mathsf{x}_{-},\mathsf{x}_{+}}(\ell_{R})\circ\mathsf{S}_{R}[\Theta](t)&=\nabla_{t}\Theta(t)+\mathsf{S}_{R}^{-1}(t)\nabla_{t}\mathsf{S}_{R}(t)\Theta(t)\notag\\
&\quad+\mathsf{S}_{R}^{-1}(t)\circ d\nabla_{1,1/2}\Lag_{H}(\ell_{R}(t))\circ\mathsf{S}_{R}(t)[\Theta(t)].
\end{align}
By (3.88), (3.91) and the uniform bound $\sup_{t\in\R}\|\mathsf{S}_{R}^{-1}(t)\|_{H^{1}(S^{1})}\le C$ (which follows from~\cite[Lemma 7.4]{I} as in the proof of Lemma 3.2), we have
\begin{equation}
\|\mathsf{S}_{R}^{-1}(t)\nabla_{t}\mathsf{S}_{R}(t)\Theta(t)\|_{T_{\ell(t)}\mathcal{F}^{1,1/2}(S^{1},N)}\le\epsilon(R)\|\Theta(t)\|_{T_{\ell(t)}\mathcal{F}^{1,1/2}(S^{1},N)}
\end{equation}
for some $\epsilon(R)>0$ with $\epsilon(R)\to 0$ as $R\to\infty$.

On the other hand, by the estimate of Lemma 3.2 and $\ell_{R}\to\ell$ as $R\to\infty$ (in the sense that (3.79), (3.80) hold), we easily see that
\begin{align}
&\|\mathsf{S}_{R}^{-1}(t)\circ d\nabla_{1,1/2}\Lag_{H}(\ell_{R}(t))\circ\mathsf{S}_{R}(t)[\Theta(t)]-d\nabla_{1,1/2}\Lag_{H}(\ell(t))[\Theta(t)]\|_{T_{\ell(t)}\mathcal{F}^{1,1/2}(S^{1},N))}\notag\\
&\le\epsilon'(R)\|\Theta(t)\|_{T_{\ell(t)}\mathcal{F}^{1,1/2}(S^{1},N))}
\end{align}
for some $\epsilon'(R)>0$ with $\epsilon'(R)\to 0$ as $R\to\infty$.
From (3.96)--(3.98), the we have the assertion of the lemma.
\hfill$\Box$

\medskip

Since the space of Fredholm operator is open in the space of bounded linear operators with respect to the operator norm, in view of Lemma 3.3, if suffices to prove that $\mathsf{S}_{R}^{-1}\circ d\mathcal{F}_{\mathsf{x}_{-},\mathsf{x}_{+}}(\ell_{R})\circ\mathsf{S}_{R}$ is Fredholm for all large $R>0$. By Lemma 3.2, this is equivalent to the assertion that $d\mathcal{F}_{\mathsf{x}_{-},\mathsf{x}_{+}}(\ell_{R})$ is Fredholm for all large $R>0$. Notice also that the Fredholm index of $d\mathcal{F}_{\mathsf{x}_{-},\mathsf{x}_{+}}(\ell)$ is equal to that of $d\mathcal{F}_{\mathsf{x}_{-},\mathsf{x}_{+}}(\ell_{R})$ for all large $R>0$. Since $\ell_{R}=\mathsf{x_{\pm}}$ when $\pm t\ge R+1$, we have reduced the problem to the case where $\ell$ is constant near $t=\pm\infty$ as asserted.

\medskip
\noindent
{\bf Step 2.} {\it Proof of Proposition 3.4 for $\ell$ constant near $t=\pm\infty$.}

\noindent
Let us assume that $\ell\in W^{1,2}_{\mathsf{x}_{-},\mathsf{x}_{+}}(\R,\mathcal{F}^{1,1/2}(S^{1},N))$ satisfies $\ell(t)=\mathsf{x}_{-}$ for $t\le -T$ and $\ell(t)=\mathsf{x}_{+}$ for $t\ge T$ for some $T>0$. We set $\ell(t)=(\phi(t),\psi(t))$ and define the parallel translation along the path $\phi(\tau,s)$ ($-T\le\tau\le t$) by $\mathsf{P}_{t}(s):T_{\phi_{-}(s)}N\to T_{\phi(t,s)}N$. Note that, by our assumption, $\mathsf{P}_{t}(s)=\mathsf{P}_{-T}(s)$ for $t\le-T$ and $\mathsf{P}_{t}(s)=\mathsf{P}_{T}(s)$ for $t\ge T$. By~\cite[Lemma 7.2]{I}, $\{\mathsf{P}_{t}\}_{t\in\R}$ defines families of bounded linear operators $$\mathsf{P}_{t}:H^{1}(S^{1},\phi_{-}^{\ast}TN)\to H^{1}(S^{1},\phi(t)^{\ast}TN)\subset H^{1}(S^{1},\underline{\R}^{k})$$ and $$\mathsf{P}_{t}={\bm 1}\otimes\mathsf{P}_{t}:H^{1/2}(S^{1},\Spin(S^{1})\otimes\phi_{-}^{\ast}TN)\to H^{1/2}(S^{1},\Spin(S^{1})\otimes\phi^{\ast}(t)TN)\subset H^{1/2}(S^{1},\Spin(S^{1})\otimes\underline{\R}^{k}),$$ which are continuous with respect to the respective operator norms. Since $\mathsf{P}_{t}$ is independent of $t$ outside of a compact set $[-T,T]$, there exists $C_{T}>0$ such that 
\begin{equation}
\|\mathsf{P}_{t}\|_{\mathsf{op}(H^{1}(S^{1},\phi_{-}^{\ast}TN),H^{1}(S^{1},\underline{\R}^{k}))}\le C_{T}
\end{equation}
and
\begin{equation}
\|\mathsf{P}_{t}\|_{H^{1/2}(S^{1},\Spin(S^{1})\otimes\phi_{-}^{\ast}TN),H^{1/2}(S^{1},\Spin(S^{1})\otimes\underline{\R}^{k}))}\le C_{T}
\end{equation}
for all $t\in\R$.
For $\Theta=(\xi,\zeta)\in W^{1,2}(\R,\mathsf{x}_{-}^{\ast}T\mathcal{F}^{1,1/2}(S^{1},N))$ (i.e., $\xi\in W^{1,2}(\R,H^{1}(S^{1},\phi_{-}^{\ast}TN))$, $\zeta\in W^{1,2}(\R,H^{1/2}(S^{1},\Spin(S^{1})\otimes\phi_{-}^{\ast}TN))$), we define $\mathsf{P}(\Theta)(t)=\mathsf{P}_{t}(\Theta(t))$. By the definition, we have $\mathsf{P}(\Theta)(t)\in\ell^{\ast}(t)T\mathcal{F}^{1,1/2}(S^{1},N)$ for all $t\in\R$. Moreover, by (3.99) and (3.100), we see that $\mathsf{P}$ defines a bounded linear map $\mathsf{P}:W^{1,2}(\R,T_{\mathsf{x}_{-}}\mathcal{F}^{1,1/2}(S^{1},N))\to W^{1,2}(\R,\ell^{\ast}T\mathcal{F}^{1,1/2}(S^{1},N))$. To see this, for $\Theta=(\xi,\zeta)\in W^{1,2}(\R,T_{\mathsf{x}_{-}}\mathcal{F}^{1,1/2}(S^{1},N))$ we have $\mathsf{P}(\Theta)(t)=(\mathsf{P}_{t}\xi(t),\mathsf{P}_{t}\zeta(t))$ and
\begin{equation}
\|\mathsf{P}_{t}\xi(t)\|_{H^{1}(S^{1})}\le C_{T}\|\xi(t)\|_{H^{1}(S^{1})},
\end{equation}
\begin{equation}
\|\nabla_{t}(\mathsf{P}_{t}\xi(t))\|_{H^{1}(S^{1})}=\|\mathsf{P}_{t}\nabla_{t}\xi(t)\|_{H^{1}(S^{1})}\le C_{T}\|\nabla_{t}\xi(t)\|_{H^{1}(S^{1})},
\end{equation}
\begin{equation}
\|\mathsf{P}_{t}\zeta(t)\|_{H^{1/2}(S^{1})}\le C_{T}\|\zeta(t)\|_{H^{1/2}(S^{1})},
\end{equation}
\begin{equation}
\|\nabla_{t}(\mathsf{P}_{t}\zeta(t))\|_{H^{1/2}(S^{1})}=\|\mathsf{P}_{t}\nabla_{t}\zeta(t)\|_{H^{1/2}(S^{1})}\le C_{T}\|\nabla_{t}\zeta(t)\|_{H^{1/2}(S^{1})}
\end{equation}
for all $t\in\R$.
By (3.101)--(3.104), integrating over $\R$ we have
$$\int_{-\infty}^{\infty}(\|\mathsf{P}_{t}\xi(t)\|_{H^{1}(S^{1})}^{2}+\|\nabla_{t}(\mathsf{P}_{t}\xi(t))\|_{H^{1}(S^{1})}^{2})\,dt\le C_{T}^{2}\int_{-\infty}^{\infty}(\|\xi(t)\|_{H^{1}(S^{1})}^{2}+\|\nabla_{t}\xi(t)\|_{H^{1}(S^{1})}^{2})\,dt$$
and
$$\int_{-\infty}^{\infty}(\|\mathsf{P}_{t}\zeta(t)\|_{H^{1/2}(S^{1})}^{2}+\|\nabla_{t}(\mathsf{P}_{t}\zeta(t))\|_{H^{1/2}(S^{1})}^{2})\,dt\le C_{T}^{2}\int_{-\infty}^{\infty}(\|\zeta(t)\|_{H^{1/2}(S^{1})}^{2}+\|\nabla_{t}\zeta(t)\|_{H^{1/2}(S^{1})}^{2})\,dt$$
and the claim is proved.

Since $\mathsf{P}$ is an isomorphism (the inverse is given by the parallel translation along the curve $\phi(-\tau)$ ($-t\le \tau\le T$)), to complete the proof it suffices to prove the Fredholm property and the index formula for the operator
\begin{equation}
\mathsf{P}^{-1}\circ d\mathcal{F}_{\mathsf{x_{-},x_{+}}}(\ell)\circ\mathsf{P}:W^{1,2}(\R,T_{\mathsf{x}_{-}}T\mathcal{F}^{1,1/2}(S^{1},N))\to L^{2}(\R,T_{\mathsf{x}_{-}}T\mathcal{F}^{1,1/2}(S^{1},N)).
\end{equation}
Since $\mathsf{P}_{t}$ is the parallel translation, we have
$$\mathsf{P}^{-1}\circ d\mathcal{F}_{\mathsf{x_{-},x_{+}}}(\ell)\circ\mathsf{P}=\nabla_{t}+\mathsf{P}^{-1}\circ d\nabla_{1,1/2}\Lag_{H}(\ell)\circ\mathsf{P}.$$
We set $A(t):=\mathsf{P}_{t}^{-1}\circ d\nabla_{1,1/2}\Lag_{H}(\ell(t))\circ\mathsf{P}_{t}$. 
We observe that
$$A(-\infty)=A(-T)=\mathsf{P}^{-1}_{-T}\circ d\nabla_{1,1/2}\Lag_{H}(\ell(-T))\circ\mathsf{P}_{-T}=d\nabla_{1,1/2}\Lag_{H}(\mathsf{x}_{-})$$
and
$$A(+\infty)=A(T)=\mathsf{P}_{T}^{-1}\circ d\nabla_{1,1/2}\Lag_{H}(\ell(T))\circ\mathsf{P}_{T}=\mathsf{P}^{-1}_{T}\circ d\nabla_{1,1/2}\Lag_{H}(\mathsf{x}_{+})\circ\mathsf{P}_{T}$$
are invertible hyperbolic operators since we have assumed that $\mathsf{x}_{-},\mathsf{x}_{+}\in\text{crit}(\Lag_{H})$ are non-degenerate. In this sense, the family of operators $\{A(t)\}_{t\in\R}$ is asymptotically hyperbolic. Moreover, we have:

\begin{lemma}
$A(t)$ takes the following form
$$A(t)=d\nabla_{1,1/2}\Lag_{H}(\mathsf{x}_{-})+K(t),$$
where $K(t):T_{\mathsf{x}_{-}}\mathcal{F}^{1,1/2}(S^{1},N)\to T_{\mathsf{x}_{-}}\mathcal{F}^{1,1/2}(S^{1},N)$ is compact for all $t\in\R$.
\end{lemma}
{\it Proof.} By an approximation argument as in~\cite[Lemma 7.6]{I}, we see that $(-\nabla^{2}_{s}+1)^{-1/2}\nabla_{s}:H^{1}(S^{1},\phi^{\ast}TN)\to H^{1}(S^{1},\phi^{\ast}TN)\subset H^{1}(S^{1},\underline{\R}^{k})$ is a compact perturbation of $(-\Delta+1)^{-1/2}D:H^{1}(S^{1},\phi^{\ast}TN)\to H^{1}(S^{1},\underline{\R}^{k})$, where the latter operator acts on $H^{1}(S^{1},\phi^{\ast}TN)$ by composition with the canonical embedding $H^{1}(S^{1},\phi^{\ast}TN)\subset H^{1}(S^{1},\underline{\R}^{k})$. Thus, the difference of their squares $(-\nabla^{2}_{s}+1)^{-1}\nabla_{s}^{2}-(-\Delta+1)^{-1}\Delta$ is also compact. Since $\nabla_{s}$ is a compact perturbation of $D$ (the later acts on vector fields along $\phi$ after composing with the canonical embedding $\phi^{\ast}TN\subset\underline{\R}^{k}$ as above), $(-\nabla^{2}_{s}+1)^{-1}\nabla_{s}^{2}$ is a compact perturbation of $(-\Delta+1)^{-1}\nabla_{s}^{2}$. By a similar reasoning, $(1+|\D_{\phi}|)^{-1}\D_{\phi}$ is a compact perturbation of $(1+|\D|)^{-1}\D_{\phi}$ (again, $(1+|\D|)^{-1}$ acts on spinors along $\phi$ after composing with the canonical embedding $\Spin(S^{1})\otimes\phi^{\ast}TN\subset\Spin(S^{1})\otimes\underline{\R}^{k}$). Thus, by the formula (2.5) in $\S2$, we see that $d\nabla_{1,1/2}\Lag_{H}(\ell(t))$ takes the following form
\begin{equation}
d\nabla_{1,1/2}\Lag_{H}(\ell(t))=\begin{pmatrix}(-\nabla_{s}^{2}+1)^{-1}\nabla_{s}^{2}&O\\
O&(1+|\D_{\phi(t)}|)^{-1}\D_{\phi(t)}\end{pmatrix}+\mathcal{K}_{t}
\end{equation}
for some compact operator $\mathcal{K}_{t}:T_{\ell(t)}\mathcal{F}^{1,1/2}(S^{1},N)\to T_{\ell(t)}\mathcal{F}^{1,1/2}(S^{1},N)$.
Similarly, we have
\begin{equation}
d\nabla_{1,1/2}\Lag_{H}(\mathsf{x}_{-})=\begin{pmatrix}(-\nabla_{s}^{2}+1)^{-1}\nabla_{s}^{2}&O\\
O&(1+|\D_{\mathsf{x}_{-}}|)^{-1}\D_{\mathsf{x}_{-}}\end{pmatrix}+\mathcal{K}
\end{equation}
for some compact operator $\mathcal{K}:T_{\mathsf{x}_{-}}\mathcal{F}^{1,1/2}(S^{1},N)\to T_{\mathsf{x}_{-}}\mathcal{F}^{1,1/2}(S^{1},N)$.

Defining $\tilde{\nabla}_{s}=\mathsf{P}^{-1}\circ\nabla_{s}\circ\mathsf{P}$ and $\tilde{\D}_{\phi}=\mathsf{P}^{-1}\circ\D_{\phi}\circ\mathsf{P}$, by (3.106) we have
\begin{equation}
\mathsf{P}_{t}^{-1}\circ d\nabla_{1,1/2}\Lag_{H}(\ell(t))\circ\mathsf{P}_{t}=\begin{pmatrix}(-\tilde{\nabla}_{s}^{2}+1)^{-1}\tilde{\nabla}_{s}^{2}&O\\
O&(1+|\tilde{\D}_{\phi(t)}|)^{-1}\tilde{\D}_{\phi(t)}\end{pmatrix}+\tilde{\mathcal{K}}_{t},
\end{equation}
where $\tilde{\mathcal{K}}_{t}=\mathsf{P}_{t}^{-1}\circ\mathcal{K}_{t}\circ\mathsf{P}_{t}:T_{\mathsf{x}_{-}}\mathcal{F}^{1,1/2}(S^{1},N)\to T_{\mathsf{x}_{-}}\mathcal{F}^{1,1/2}(S^{1},N)$ is compact. By~\cite[Lemma 7.6]{I}, $\{(1+|\tilde{\D}_{\phi(t)}|)^{-1}\tilde{\D}_{\phi(t)}-(1+|\D_{\mathsf{x}_{-}}|)^{-1}\D_{\mathsf{x}_{-}}\}_{t\in\R}$ defines a continuous family of compact operators. By the same reasoning, $\{(-\tilde{\nabla}_{s}^{2}+1)^{-1/2}\tilde{\nabla}_{s}-(-\nabla_{s}^{2}+1)^{-1/2}\nabla_{s}\}_{t\in\R}$ is a continuous family of compact operators and therefore the difference of their squares $\{(-\tilde{\nabla}_{s}^{2}+1)^{-1}\tilde{\nabla}_{s}^{2}-(-\nabla_{s}^{2}+1)^{-1}\nabla_{s}^{2}\}_{t\in\R}$ also defines a continuous family of compact operators. From this observation, comparing (3.107) and (3.108), the assertion of the lemma follows.
\hfill$\Box$

\medskip

We now complete the proof of Step 2. By Lemma 3.4, we are now in a position to apply the result of Abbondandolo-Majer~\cite[Theorem 3.4]{AM1}, \cite[Theorem B]{AM2}. In fact, using the non-degeneracy of $\mathsf{x}_{-}\in\text{crit}(\Lag_{H})$, $d\nabla_{1,1/2}\Lag_{H}(\mathsf{x}_{-})$ is a hyperbolic operator. Combining this with Lemma 3.4, we can apply Theorem 3.4 of~\cite{AM1} (see also Theorem B of~\cite{AM1}) and conclude that $\mathsf{P}^{-1}\circ d\mathcal{F}_{\mathsf{x_{-},x_{+}}}(\ell)\circ\mathsf{P}=\nabla_{t}+A(t):W^{1,2}(\R,T_{\mathsf{x}_{-}}T\mathcal{F}^{1,1/2}(S^{1},N))\to L^{2}(\R,T_{\mathsf{x}_{-}}T\mathcal{F}^{1,1/2}(S^{1},N))$ is Fredholm. In this setting, by~\cite[$\S4.3$]{AV}, its index is given by the spectral flow of the family of operators $\{A(t)\}_{t\in\R}$:
\begin{equation}
\text{index}(d\mathcal{F}_{\mathsf{x}_{-},\mathsf{x}_{+}}(\ell))=\text{index}(\mathsf{P}^{-1}\circ d\mathcal{F}_{\mathsf{x}_{-},\mathsf{x}_{+}}(\ell)\circ\mathsf{P})=\mathsf{sf}\{A(t)\}_{-\infty\le t\le +\infty}.
\end{equation}
By our definition of the relative Morse index in Definition 2.1 in $\S2$ (we can choose the base point of the path as $(\phi_{0},\psi_{0})=(\phi(0),\psi(0))$, we have
\begin{align}
\mathsf{sf}\{A(t)\}_{-\infty\le t\le+\infty}&=\mathsf{sf}\{A(t)\}_{-T\le t\le T}\notag\\
&=\mathsf{sf}\{A(t)\}_{-T\le t\le 0}+\mathsf{sf}\{A(t)\}_{0\le t\le T}\notag\\
&=\mu_{H}(\mathsf{x}_{-})-\mu_{H}(\mathsf{x}_{+}).
\end{align}
By (3.109) and (3.110), we have the desired index formula.
\hfill$\Box$

\medskip

For $\mathsf{x}_{-}=(\phi_{-},\psi_{-}),\mathsf{x}_{+}=(\phi_{+},\psi_{+})\in\text{crit}(\Lag_{H})$, we denote by $\mathcal{M}(\mathsf{x}_{-},\mathsf{x}_{+})$ the space of solutions to the negative gradient flow equations starting from $\mathsf{x}_{-}$ and ending at $\mathsf{x}_{+}$,
$$\mathcal{M}(\mathsf{x}_{-},\mathsf{x}_{+})=\{\ell\in W^{1,2}_{\mathsf{x}_{-},\mathsf{x}_{+}}(\R,\mathcal{F}^{1,1/2}(S^{1},N)):\mathcal{F}_{\mathsf{x}_{-},\mathsf{x}_{+}}(\ell)=0\}.$$
By Proposition 3.4 and the implicit function theorem, if $0\in L^{2}(\R,T\mathcal{F}^{1,1/2}(S^{1},N))$ is a regular value of $\mathcal{F}_{\mathsf{x}_{-},\mathsf{x}_{+}}$, $\mathcal{M}(\mathsf{x}_{-},\mathsf{x}_{+})$ is a manifold of dimension $\mu_{H}(\mathsf{x}_{-})-\mu_{H}(\mathsf{x}_{+})$. When $\mathsf{x}_{-}\ne\mathsf{x}_{+}$, $\R$ acts freely on $\mathcal{M}(\mathsf{x}_{-},\mathsf{x}_{+})$ by time shift: $\R\times\mathcal{M}(\mathsf{x}_{-},\mathsf{x}_{+})\ni(a,\ell)\mapsto\ell(\cdot+a)\in\mathcal{M}(\mathsf{x}_{-},\mathsf{x}_{+})$. The quotient of $\mathcal{M}(\mathsf{x}_{-},\mathsf{x}_{+})$ by this action is denoted by $\hat{\mathcal{M}}(\mathsf{x}_{-},\mathsf{x}_{+}):=\mathcal{M}(\mathsf{x}_{-},\mathsf{x}_{+})/\R$. This is the set of negative gradient flow lines and if $0\in L^{2}(\R,T\mathcal{F}^{1,1/2}(S^{1},N))$ is a regular value of $\mathcal{F}_{\mathsf{x}_{-},\mathsf{x}_{+}}$, $\hat{\mathcal{M}}(\mathsf{x}_{-},\mathsf{x}_{+})$ is a manifold of dimension $\mu_{H}(\mathsf{x}_{-})-\mu_{H}(\mathsf{x}_{+})-1$.

\subsection{Compactness of flow lines}

In this section, we prove compactness properties of the spaces $\mathcal{M}(\mathsf{x}_{-},\mathsf{x}_{+})$ and $\hat{\mathcal{M}}(\mathsf{x}_{-},\mathsf{x}_{+})$ for $\mathsf{x}_{-}=(\phi_{-},\psi_{-}),\mathsf{x}_{+}=(\phi_{+},\psi_{+})\in\text{crit}(\Lag_{H})$.
The first result is a consequence of Proposition 3.3.

\begin{proposition}
$\mathcal{M}(\mathsf{x}_{-},\mathsf{x}_{+})\subset C^{0}_{\text{loc}}(\R,\mathcal{F}^{1,1/2}(S^{1},N))$ is relatively compact.
\end{proposition}
{\it Proof.} Let $\ell=(\psi,\phi)\in\mathcal{M}(\mathsf{x}_{-},\mathsf{x}_{+})$ be arbitrary. $\Lag_{H}(\ell(t))$ is non-increasing along the negative flow line $\ell$ and we have
$$\Lag_{H}(\mathsf{x}_{-})\le\Lag_{H}(\ell(t))\le\Lag_{H}(\mathsf{x}_{+})$$
for all $t\in\R$. By Proposition 3.3, this implies the following for some $C=C(\mathsf{x}_{-},\mathsf{x}_{+})$ depending only on $\Lag_{H}(\mathsf{x}_{-})$ and $\Lag_{H}(\mathsf{x}_{+})$,
\begin{equation}
\sup_{t\in\R}\|\phi(t)\|_{C^{2,\frac{2}{3}}(S^{1})}+\sup_{t\in\R}\|\psi(t)\|_{C^{1,\frac{2}{3}}(S^{1})}\le C(\mathsf{x}_{-},\mathsf{x}_{+}).
\end{equation}
By (3.111) and the compactness of the embedding $C^{2,\frac{2}{3}}(S^{1})\times C^{1,\frac{2}{3}}(S^{1})\subset H^{1}(S^{1})\times H^{1/2}(S^{1})$, $\{\ell(t):\ell\in\mathcal{M}(\mathsf{x}_{-},\mathsf{x}_{+}),~t\in\R\}\subset\mathcal{F}^{1,1/2}(S^{1},N)$ is relatively compact.
By the gradient flow equation (3.7), (3.8) and the elliptic regularity for the elliptic operators $-\Delta+1$ and $1+|\D|$, (3.111) implies that there exists another constant $C=C(\mathsf{x}_{-},\mathsf{x}_{+})$ depending only on $\Lag_{H}(\mathsf{x}_{-})$ and $\Lag_{H}(\mathsf{x}_{+})$ such that 
\begin{equation}
\sup_{t\in\R}\|\partial_{t}\phi(t)\|_{C^{2,\frac{2}{3}}(S^{1})}+\sup_{t\in\R}\|\partial_{t}\psi(t)\|_{C^{1,\frac{2}{3}}(S^{1})}\le C(\mathsf{x}_{-},\mathsf{x}_{+}).
\end{equation}
Using the mean value theoremS and the compactness of the embedding $C^{2,\frac{2}{3}}(S^{1})\times C^{1,\frac{2}{3}}(S^{1})\subset H^{1}(S^{1})\times H^{1/2}(S^{1})$, (3.112) implies that $\mathcal{M}(\mathsf{x}_{-},\mathsf{x}_{+})\subset C^{0}(\R,\mathcal{F}^{1,1/2}(S^{1},N))$ is an equicontinuous family. Therefore, by the Ascoli-Arzel\`a theorem, $\mathcal{M}(\mathsf{x}_{-},\mathsf{x}_{+})\subset C^{0}_{\text{loc}}(\R,\mathcal{F}^{1,1/2}(S^{1},N))$ is relatively compact as asserted. This completes the proof.
\hfill$\Box$

\medskip

Under the conditions (1.2) and (1.4), $\Lag_{H}$ satisfies the Palais-Smale condition:

\begin{proposition} Assume that $H\in C^{1}(S^{1}\times\Spin(S^{1})\otimes TN)$ satisfies (1.2) and (1.4). Then $\Lag_{H}$ satisfies the Palais-Smale condition in the following sense:\\
 Let us assume that $\{\ell_{n}\}=\{(\phi_{n},\psi_{n})\}\in\mathcal{F}^{1,1/2}(S^{1},N)$ satisfies $\sup_{n\ge 1}\Lag_{H}(\ell_{n})<+\infty$ and $\|\nabla_{1,1/2}\Lag_{H}(\ell_{n})\|_{T_{\ell_{n}}\mathcal{F}^{1,1/2}(S^{1},N)}\to 0$ as $n\to\infty$. Then there exists a subsequence $\{\ell_{n_{k}}\}\subset\{\ell_{n}\}$ and $\ell_{\infty}=(\phi_{\infty},\psi_{\infty})\in\text{crit}(\Lag_{H})$ such that $\|\ell_{n}-\ell_{\infty}\|_{H^{1}(S^{1},\R^{k})\times H^{1/2}(S^{1},\Spin(S^{1})\otimes\underline{\R}^{k})}\to 0$ as $n\to\infty$.
\end{proposition}
\textit{Proof.}
Note that, in~\cite[Lemma 5.1]{I}, the Palais-Smale condition for $\Lag_{H}$ has been verified under slightly different condition on $H$. To the present case, essentially the same argument in that paper applies and we omit the details here.
\hfill$\Box$

\medskip

We next turn to the compactness issue of the set of flow lines. Recall from the last subsection, for $\mathsf{x}_{-}=(\phi_{-},\psi_{-}),\mathsf{x}_{+}=(\phi_{+},\psi_{+})\in\text{crit}(\Lag_{H})$, $\hat{\mathcal{M}}(\mathsf{x}_{-},\mathsf{x}_{+})$ denotes the set of negative gradient flow lines connecting $\mathsf{x}_{-}$ and $\mathsf{x}_{+}$. For $\ell\in\mathcal{M}(\mathsf{x}_{-},\mathsf{x}_{+})$, the associated flow line (compactified at both ends) $\ell(\overline{\R}):=\{\mathsf{x}_{-}\}\cup\ell(\R)\cup\{\mathsf{x}_{+}\}\subset\mathcal{F}^{1,1/2}(S^{1},N)$ is compact. We define on $\hat{\mathcal{M}}(\mathsf{x}_{-},\mathsf{x}_{+})$ a metric $d_{H}$ via
$$d_{H}(\hat{\ell}_{1},\hat{\ell}_{2})=d_{\text{Hausdorff}}(\ell_{1}(\overline{\R}),\ell_{2}(\overline{\R})),$$
where $\hat{\ell}_{i}\in\hat{\mathcal{M}}(\mathsf{x}_{-},\mathsf{x}_{+})$ ($i=1,2$) are flow lines corresponding to $\ell_{i}\in\mathcal{M}(\mathsf{x}_{-},\mathsf{x}_{+})$ and $d_{\text{Hausdorff}}$ is the Hausdorff distance defined on compact sets in $\mathcal{F}^{1,1/2}(S^{1},N)$.
In the next proposition, we show that the metric space $(\hat{\mathcal{M}}(\mathsf{x}_{-},\mathsf{x}_{+}),d_{H})$ is relatively compact. We also give a characterization of its closure $\overline{\hat{\mathcal{M}}(\mathsf{x}_{-},\mathsf{x}_{+})}$.

\begin{proposition}
The metric space $(\hat{\mathcal{M}}(\mathsf{x}_{-},\mathsf{x}_{+}),d_{H})$ is relatively compact, i.e., for any sequence $\{\hat{\ell}_{n}\}\subset\hat{\mathcal{M}}(\mathsf{x}_{-},\mathsf{x}_{+})$, there exists a subsequence $\{\hat{\ell}_{n_{k}}\}$ and a compact subset $K_{\infty}\in\mathcal{F}^{1,1/2}(S^{1},N)$ such that $d_{{\rm Hausdorff}}(\ell_{n_{k}}(\overline{\R}),K_{\infty})\to 0$ as $n_{k}\to\infty$. The limit $K_{\infty}$ has the following properties:
\begin{enumerate}[(1)]
\item $K_{\infty}$ is invariant under the negative gradient flow of $\Lag_{H}$,
\item If $\Lag_{H}$ is a Morse function on $\mathcal{F}^{1,1/2}(S^{1},N)$, there exist finitely many critical points $\mathsf{x}_{-}=\mathsf{x}_{0},\mathsf{x}_{1},\ldots,\mathsf{x}_{k-1},\mathsf{x}_{k}=\mathsf{x}_{+}\in\text{crit}(\Lag_{H})$ and $\ell_{i}\in\mathcal{M}(\mathsf{x}_{i},\mathsf{x}_{i+1})$ for $0\le i\le k-1$ such that $K_{\infty}=\bigcup_{i=0}^{k-1}\ell_{i}(\overline{\R})$.
\end{enumerate}
\end{proposition}
{\it Proof.} Since $\Lag_{H}$ is non-increasing along the negative gradient flow, we have $\Lag_{H}(\mathsf{x}_{+})\le\Lag_{H}(\ell_{n}(t))\le\Lag_{H}(\mathsf{x}_{-})$ for all $n\ge 1$ and $t\in\R$. Using Proposition 3.3, this implies the following uniform estimate
\begin{equation}
\sup_{t\in\R}\|\phi_{n}(t)\|_{C^{2,2/3}(S^{1})}+\sup_{t\in\R}\|\psi_{n}(t)\|_{C^{1,2/3}(S^{1})}\le C
\end{equation}
for some $C>0$ independent of $n$, where $\ell_{n}=(\phi_{n},\psi_{n})$. Thus, there exists $R>0$ such that $$\ell_{n}(\overline{\R})\subset\overline{B}_{R}(C^{2,2/3}(S^{1})\times C^{1,2/3}(S^{1}))\cap\mathcal{F}^{1,1/2}(S^{1},N),$$ where $\overline{B}_{R}(C^{2,2/3}(S^{1})\times C^{1,2/3}(S^{1}))$ is the closed ball of radius $R$ with center at $(0,0)$ in $C^{2,2/3}(S^{1},\R^{k})\times C^{1,2/3}(S^{1},\Spin(S^{1})\otimes\underline{\R}^{k})$. Since the embedding $C^{2,2/3}(S^{1})\times C^{1,2/3}(S^{1})\subset H^{1}(S^{1})\times H^{1/2}(S^{1})$ is compact, the inclusion $B_{R}(C^{2,2/3}(S^{1})\times C^{1,2/3}(S^{1}))\cap\mathcal{F}^{1,1/2}(S^{1},N)\subset\mathcal{F}^{1,1/2}(S^{1},N)$ is compact. Since the set of compact subsets in a compact metric space is compact with respect to the Hausdorff distance, there exists a subsequence $\{\ell_{n_{k}}\}$ and a compact subset $K_{\infty}\subset B_{R}(C^{2,2/3}(S^{1})\times C^{1,2/3}(S^{1}))\cap\mathcal{F}^{1,1/2}(S^{1},N)$ such that 
$$d_{\text{Hausdorff}}(\ell_{n_{k}}(\overline{\R}),K_{\infty})\to 0 \text{ as }n_{k}\to\infty.$$ Because all the flow lines $\ell_{n_{k}}(\overline{\R})$ are invariant under the negative gradient flow, by the Hausdorff convergence, $K_{\infty}$ is also invariant under the negative gradient flow and this finishes the proof of (1).\\

 To prove the second assertion, assume $\Lag_{H}$ is Morse. Since $\Lag_{H}(\mathsf{x}_{-})\le\Lag_{H}\le\Lag_{H}(\mathsf{x}_{+})$ on $K_{\infty}$, by the Hausdorff convergence, the Palais-Smale condition (Proposition 3.6) and the Morse property, there are at most finitely many critical points on $K_{\infty}$. Since $\mathsf{x}_{\pm}\in\ell_{n}(\overline{\R})$ for all $n$, we also have $\mathsf{x}_{\pm}\in K_{\infty}$ by the Hausdorff convergence. To complete the proof, we need to show the following:
\begin{enumerate}[(i)]
\item For any $c\in\R$ with $\Lag_{H}(\mathsf{x}_{+})\le c\le\Lag_{H}(\mathsf{x}_{-})$, the set $K_{\infty}\cap\Lag_{H}^{-1}(c)$ consists of exactly one point.
\item For any $c\in\R$ with $c<\Lag_{H}(\mathsf{x}_{+})$ or $c>\Lag_{H}(\mathsf{x}_{-})$, we have $K_{\infty}\cap\Lag_{H}^{-1}(c)=\emptyset$.
\end{enumerate}
As for (ii), observe that for any $n$, $\ell_{n}(\overline{\R})\cap\Lag_{H}^{-1}(c)\ne\emptyset$ if and only if $\Lag_{H}(\mathsf{x}_{+})\le c\le\Lag_{H}(\mathsf{x}_{-})$. Thus, the assertion (ii) follows form the Hausdorff convergence $d_{\text{Hausdorff}}(\ell_{n_{k}}(\overline{\R}),K_{\infty})\to 0$. Therefore, we also have $K_{\infty}\cap\Lag_{H}^{-1}(c)\ne\emptyset$ for $\Lag_{H}(\mathsf{x}_{+})\le c\le\Lag_{H}(\mathsf{x}_{-})$. Thus, it remains to prove $\#(K_{\infty}\cap\Lag_{H}^{-1}(c))=1$ for any $c\in\R$ with $\Lag_{H}(\mathsf{x}_{+})\le c\le\Lag_{H}(\mathsf{x}_{-})$. In order to prove this fact, we assume by contradiction that there are two distinct points $\mathsf{y,z}\in K_{\infty}\cap\Lag_{H}^{-1}(c)$ for some $c\in\R$ with $\Lag_{H}(\mathsf{x}_{+})\le c\le\Lag_{H}(\mathsf{x}_{-})$. Again, from the Hausdorff convergence $\ell_{n_{k}}(\overline{\R})\to K_{\infty}$ and the fact that $\ell_{n_{k}}(\overline{\R})$ is a flow line, there exists $\mathsf{y}_{k}\in\ell_{n_{k}}(\overline{\R})$ and $t_{k}\in\R$ such that $\mathsf{y}_{k}\to\mathsf{y}$ and $\mathsf{z}_{k}:=\ell(t_{k},\mathsf{y}_{k})\to\mathsf{z}$ as $k\to\infty$, where $\ell(t,\mathsf{y})$ is the solution to the negative gradient flow equation (3.9) with $\ell(0,\mathsf{y})=\mathsf{y}$. We may assume $t_{k}\ge 0$, otherwise we simply reverse the roles of $\mathsf{x}_{k}$ and $\mathsf{y}_{k}$. Under this assumption, we have
\begin{align}
&\int_{0}^{t_{k}}\|\nabla_{1,1/2}\Lag_{H}(\ell_{n_{k}}(t,\mathsf{y}_{k}))\|_{T\mathcal{F}^{1,1/2}(S^{1},N)}^{2}\,dt\notag\\
&=\Lag_{H}(\mathsf{y}_{k})-\Lag_{H}(\ell(t_{k},\mathsf{y}_{k}))\to\Lag_{H}(\mathsf{y})-\Lag_{H}(\mathsf{z})=0.
\end{align}
Knowing that $\Lag_{H}$ is Morse and satisfies the Palais-Smale condition (see Proposition 3.6), the set $\text{crit}(\Lag_{H})\cap\Lag_{H}^{-1}([\Lag_{H}(\mathsf{z}),\Lag_{H}(\mathsf{x})])$ is finite, (3.114) implies that $t_{k}\to 0$ or $\ell([0,t_{k}],\mathsf{y}_{k})$ converges to some critical point of $\Lag_{H}$. In any case, we arrive at the contradiction $\mathsf{y}=\mathsf{z}$. This completes the proof.
\hfill$\Box$

\medskip

We assume that $\Lag_{H}$ is a Morse function on $\mathcal{F}^{1,1/2}(S^{1},N)$ and $$\mathcal{F}_{\mathsf{x}_{-},\mathsf{x}_{+}}:W^{1,2}_{\mathsf{x}_{-},\mathsf{x}_{+}}(\R,\mathcal{F}^{1,1/2}(S^{1},N))\to L^{2}(\R,T\mathcal{F}^{1,1/2}(S^{1},N))$$ defined by (3.78) has $0\in L^{2}(\R,T\mathcal{F}^{1,1/2}(S^{1},N))$ as a regular value. Under this assumption, if $\mathsf{x}_{-}\ne\mathsf{x}_{+}$, $\hat{\mathcal{M}}(\mathsf{x}_{-},\mathsf{x}_{+})$ is a manifold of dimension $\mu_{H}(\mathsf{x}_{-})-\mu_{H}(\mathsf{x}_{+})-1$ as remarked at the end of the last subsection. Under these assumptions, as a corollary of Proposition 3.7, we have:

\begin{corollary}
Under the above assumption, we have
\begin{enumerate}[(1)]
\item If $\mu_{H}(\mathsf{x}_{-})-\mu_{H}(\mathsf{x}_{+})=1$, $\hat{\mathcal{M}}(\mathsf{x}_{-},\mathsf{x}_{+})$ is a compact $0$-dimensional manifold, hence it consists of at most finitely many points. 
\item If $\mu_{H}(\mathsf{x}_{-})-\mu_{H}(\mathsf{x}_{+})=2$, $\hat{\mathcal{M}}(\mathsf{x}_{-},\mathsf{x}_{+})$ is a 1-dimensional manifold and has a compactification $\overline{\hat{\mathcal{M}}(\mathsf{x}_{-},\mathsf{x}_{+})}$.
The boundary $\partial\overline{\hat{\mathcal{M}}(\mathsf{x}_{-},\mathsf{x}_{+})}$ of this manifold (if non-empty) consists of broken flow lines of the form $\ell_{1}(\overline{\R})\cup\ell_{2}(\overline{\R})$ for some $\ell_{1}\in\mathcal{M}(\mathsf{x}_{-},\mathsf{y})$ and $\ell_{2}\in\mathcal{M}(\mathsf{y},\mathsf{x}_{+})$, where $\mathsf{y}\in\text{crit}(\Lag_{H})$ with $\mu_{H}(\mathsf{x}_{-})-\mu_{H}(\mathsf{y})=1$.
\end{enumerate}
\end{corollary}

\subsection{Gluing negative  gradient flow lines}

We shall show that the converse of the compactness result (Corollary 3.1 (2)) holds under a certain transversality assumption.

Let $\mathsf{x,y,z}\in\text{crit}(\Lag_{H})$ be such that $\mu_{H}(\mathsf{x})=\mu_{H}(\mathsf{z})+1=\mu_{H}(\mathsf{y})+2$ and consider $\ell_{1}\in\mathcal{M}(\mathsf{x},\mathsf{z})$ and $\ell_{2}\in\mathcal{M}(\mathsf{z},\mathsf{y})$. We assume that $d\mathcal{F}_{\mathsf{x,z}}(\ell_{i}):W^{1,2}(\R,\ell_{i}^{\ast}T\mathcal{F}^{1,1/2}(S^{1},N))\to L^{2}(\R,\ell_{i}^{\ast}T\mathcal{F}^{1,1/2}(S^{1},N))$ are surjective for $i=1,2$. We shall prove in what follows that under the above assumption, there exists a family of solutions $\{\ell_{1,2;R}\}_{R}$ to the negative gradient flow equation $\frac{d\ell_{1,2;R}}{dt}+\nabla_{1,1/2}\Lag_{H}(\ell_{1,2;R})=0$ parametrized by large $R>>1$ such that the flow line $\ell_{1,2;R}(\overline{\R})$ converges to $\ell_{1}(\overline{\R})\cup\ell_{2}(\overline{\R})$ in the Hausdorff metric as $R\to\infty$. To construct such solutions, we first construct pregluings $\ell_{1}\#_{R}\ell_{2}$ of $\ell_{1}$ and $\ell_{2}$ for all large $R>>1$ which are good approximate solutions to the negative gradient flow equation. 

We set $\ell_{i}=(\phi_{i},\psi_{i})$ for $i=1,2$ and $\mathsf{z}=(\phi_{\mathsf{z}},\psi_{\mathsf{z}})$.

\medskip

\noindent
\underline{Construction of pregluing $\phi_{1}\#_{R}\phi_{2}$:} 

Let $\alpha\in C^{\infty}(\R)$ be such that $\alpha(t)=0$ for $t\le -1$, $\alpha(t)=1$ for $t\ge 1$ and $0\le\alpha(t)\le 1$ for $-1\le t\le 1$. For $R>0$, we set $\alpha_{R}(t)=\alpha(R^{-1}t)$. By definition (see (3.74), (3.75)), $\phi_{1}$ and $\phi_{2}$ are expressed as
$$\phi_{1}(t)=\exp_{\phi_{\mathsf{z}}}(X_{1,+}(t))\quad\text{for $t\ge T_{1}$},$$
$$\phi_{2}(t)=\exp_{\phi_{\mathsf{z}}}(X_{2,-}(t))\quad\text{for $t\le-T_{2}$},$$
where $T_{1},T_{2}>0$, $X_{1,+}\in W^{1,2}([T_{1},\infty),H^{1}(S^{1},\phi_{\mathsf{z}}^{\ast}TN))$ and $X_{2,-}\in W^{1,2}((-\infty,-T_{2}],H^{1}(S^{1},\phi_{\mathsf{z}}^{\ast}TN))$.

For $R>\max\{T_{1},T_{2}\}$, we define  pregluing $\phi_{1}\#_{R}\phi_{2}$ by
$$\phi_{1}\#_{R}\phi_{2}=\left\{\begin{array}{lll}
\phi_{1}(t+2R)&\ (t\le-R)\\
\exp_{\phi_{\mathsf{z}}}((1-\alpha_{R}(t))X_{1,+}(t+2R)+\alpha_{R}(t)X_{2,-}(t-2R))&\ (-R\le t\le R)\\
\phi_{2}(t-2R)&\ (t\ge R),
\end{array}\right.$$

\noindent
\underline{Construction of pregluing $\psi_{1}\#_{R}\psi_{2}$:}

By definition (see (3.76), (3.77)), $\psi_{1}$ and $\psi_{2}$ are expressed as
$$\psi_{1}(t)=\mathsf{S}_{1,+,t}(\psi_{\mathsf{z}}+\xi_{1,+}(t))\quad\text{for $t\ge T_{1}$},$$
$$\psi_{2}(t)=\mathsf{S}_{2,-,t}(\psi_{\mathsf{z}}+\xi_{2,-}(t))\quad\text{for $t\le-T_{2}$},$$
where $\xi_{1,+}\in W^{1,2}([T_{1},\infty),H^{1/2}(S^{1},\Spin(S^{1})\otimes\phi_{\mathsf{z}}^{\ast}TN))$, $\xi_{2,-}\in W^{1,2}((-\infty,-T_{2}],H^{1/2}(S^{1},\Spin(S^{1})\otimes\phi_{\mathsf{z}}^{\ast}TN))$, where $\mathsf{S}_{1,+,t}:\phi_{\mathsf{z}}(t)^{\ast}TN\to\phi_{1}(t)^{\ast}TN$ is the parallel translation along the path $[0,1]\ni\tau\mapsto\exp_{\phi_{\mathsf{z}}}(\tau X_{1,+}(t))\in N$ for $t\ge T_{1}$ and $\mathsf{S}_{2,-,t}:\phi_{\mathsf{z}}(t)^{\ast}TN\to\phi_{2}(t)^{\ast}TN$ is the parallel translation along path $[0,1]\ni\tau\mapsto\exp_{\phi_{\mathsf{z}}}(\tau X_{2,-}(t))$ for $t\le-T_{2}$. 

For $R>\max\{T_{1},T_{2}\}$, we define pregluing $\psi_{1}\#_{R}\psi_{2}$ by
$$\psi_{1}\#_{R}\psi_{2}=\left\{\begin{array}{lll}
\psi_{1}(t+2R)&\ (t\le-R)\\
\mathsf{S}_{1,2,R;t}(\psi_{\mathsf{z}}+(1-\alpha_{R}(t))\xi_{1,+}(t+2R)+\alpha_{R}(t)\xi_{2,-}(t-2R))&\ (-R\le t\le R)\\
\psi_{2}(t-2R)&\ (t\ge R),
\end{array}\right.$$
where $\mathsf{S}_{1,2,R;t}:\phi_{\mathsf{z}}(t)^{\ast}TN\to(\phi_{1}\#_{R}\phi_{2})(t)^{\ast}TN$ is the parallel translation along path $[0,1]\ni\tau\mapsto\exp_{\phi_{\mathsf{z}}}(\tau X_{1,2,R;t})$, where we set $X_{1,2,R;t}=(1-\alpha_{R}(t))X_{1,+}(t+2R)+\alpha_{R}(t)X_{2,-}(t-2R)$.

\medskip

We then define pregluing $\ell_{1}\#_{R}\ell_{2}$ by $\ell_{1}\#_{R}\ell_{2}=(\phi_{1}\#_{R}\phi_{2},\psi_{1}\#_{R}\psi_{2})$.
The following properties of the pregluing can be checked from the construction:

\noindent
$\bullet$ The flow line $\ell_{1}\#_{R}\ell_{2}(\overline{\R})$ converges to the broken flow line $\ell_{1}(\overline{\R})\cup\ell_{2}(\overline{\R})$ in the Hausdorff metric as $R\to\infty$.

\noindent
$\bullet$ $\mathcal{F}_{\mathsf{x,y}}(\ell_{1}\#_{R}\ell_{2})\to 0$ as $R\to\infty$. Thus, for large $R>0$, $\ell_{1}\#_{R}\ell_{2}$ is a good approximate solution to the negative gradient flow equation.

A genuine solution to the equation $\mathcal{F}_{\mathsf{x,y}}(\ell)=0$ can be constructed as a small perturbation of $\ell_{1}\#_{R}\ell_{2}$.

\begin{proposition}
Let us assume that $\mathsf{x,y,z}\in\text{crit}(\Lag_{H})$ with $\mu_{H}(\mathsf{x})=\mu_{H}(\mathsf{z})+1=\mu_{H}(\mathsf{y})+2$. Let $\ell_{1}\in\mathcal{M}(\mathsf{x,z})$ and $\ell_{2}\in\mathcal{M}(\mathsf{z,y})$ be regular solutions to the negative gradient flow equation in the sense that $d\mathcal{F}_{\mathsf{x,z}}(\ell_{i}):W^{1,2}(\R,\ell_{i}^{\ast}T\mathcal{F}^{1,1/2}(S^{1},N))\to L^{2}(\R,\ell_{i}^{\ast}T\mathcal{F}^{1,1/2}(S^{1},N))$ are surjective for $i=1,2$. Then there exists $R_{0}>0$ such that the following holds:\\

 For any $R\ge R_{0}$, there exists a unique small ${}^{t}(X_{1,2;R},\xi_{1,2;R})$ with $${}^{t}(X_{1,2;R},\xi_{1,2;R})\in\begin{pmatrix}1&0\\-\Gamma(\ell_{1}\#_{R}\ell_{2})&1\end{pmatrix}\big(\ker(d\mathcal{F}_{\mathsf{x,y}}(\ell_{1}\#_{R}\ell_{2}))^{\perp}\big),$$
 such that $\ell_{1,2;R}:=(\exp_{\phi_{1}\#_{R}\phi_{2}}(X_{1,2;R}),S(X_{1,2;R})(\psi_{1}\#_{R}\psi_{2}+\xi_{1,2;R}))\in W^{1,2}(\R,\mathcal{F}^{1,1/2}(S^{1},N))$ solves the equation $\mathcal{F}_{\mathsf{x,y}}(\ell_{1,2;R})=0,$ where $$\Gamma(\ell_{1}\#_{R}\ell_{2})[X]=-(\psi_{1}\#_{R}\psi_{2})^{k}X^{j}\Gamma_{jk}^{i}(\phi_{1}\#_{R}\phi_{2})\otimes\frac{\partial}{\partial y^{i}}(\phi_{1}\#_{R}\phi_{2})$$
for $X=X^{j}(t)\frac{\partial}{\partial y^{j}}((\phi_{1}\#_{R}\phi_{2})(t))\in W^{1,2}(\R,H^{1}(S^{1},(\phi_{1}\#_{R}\phi_{2})^{\ast}TN))$, $\ker(d\mathcal{F}_{\mathsf{x,y}}(\ell_{1}\#_{R}\ell_{2}))^{\perp}$ is the orthogonal complement of $\ker(d\mathcal{F}_{\mathsf{x,y}}(\ell_{1}\#_{R}\ell_{2}))\subset W^{1,2}(\R,(\ell_{1}\#_{R}\ell_{2})^{\ast}T\mathcal{F}^{1,1/2}(S^{1},N))$ and $S(X_{1,2;R}):T_{\phi_{1}\#_{R}\phi_{2}}N\to T_{\exp_{\phi_{1}\#_{R}\phi_{2}}(X_{1,2;R})}N$ is the parallel translation along the path $[0,1]\ni\tau\mapsto\exp_{\phi_{1}\#_{R}\phi_{2}}(\tau X_{1,2;R})\in N$.
\end{proposition}
\textit{Proof.} 
By the regularity of $\ell_{1}$ and $\ell_{2}$, we first observe that for large $R>0$, 
$$d\mathcal{F}_{\mathsf{x,y}}(\ell_{1}\#_{R}\ell_{2}):W^{1,2}(\R,(\ell_{1}\#_{R}\ell_{2})^{\ast}T\mathcal{F}^{1,1/2}(S^{1},N))\to L^{2}(\R,(\ell_{1}\#_{R}\ell_{2})^{\ast}T\mathcal{F}^{1,1/2}(S^{1},N))$$
 is surjective. To see this, we notice that Proposition 3.4 and the surjectivity of $d\mathcal{F}_{\mathsf{x,z}}(\ell_{1})$ and $d\mathcal{F}_{\mathsf{z,y}}(\ell_{2})$ implies the existence of a bounded right inverses $\mathcal{R}_{\mathsf{x,z};\ell_{1}}:L^{2}(\R,\ell_{1}^{\ast}T\mathcal{F}^{1,1/2}(S^{1},N))\to W^{1,2}(\R,\ell_{1}^{\ast}\mathcal{F}^{1,1/2}(S^{1},N))$ and $\mathcal{R}_{\mathsf{z,y};\ell_{2}}:L^{2}(\R,\ell_{2}^{\ast}T\mathcal{F}^{1,1/2}(S^{1},N))\to W^{1,2}(\R,\ell_{2}^{\ast}\mathcal{F}^{1,1/2}(S^{1},N))$ of $d\mathcal{F}_{\mathsf{x,z}}(\ell_{1})$ and $d\mathcal{F}_{\mathsf{z,y}}(\ell_{2})$, respectively. By gluing these right inverses appropriately, we obtain an approximate right inverse of $d\mathcal{F}_{\mathsf{x,y}}(\ell_{1}\#_{R}\ell_{2})$ as follows. Let $\beta,\gamma\in C^{\infty}(\R)$ be such that $\beta(t)=1$ for $t\le -1$, $\beta(t)=0$ for $t\ge 1$, $0\le\beta(t)\le 1$ for $-1\le t\le 1$, $\gamma(t)=\beta(-t)$ for all $t\in\R$ and $\beta(t)^{2}+\gamma(t)^{2}=1$ for all $t\in\R$. For $R>0$, we set $\beta_{R}(t)=\beta(R^{-1}t)$ and $\gamma_{R}(t)=\gamma(R^{-1}t)$. We define
\begin{align*}
\mathcal{R}_{\mathsf{x,z,y};\ell_{1}\#_{R}\ell_{2}}&=\left\{\begin{array}{lll}
\tau_{2R}\mathcal{R}_{\mathsf{x,z};\ell_{1}}\tau_{-2R}&\ (t\le-R)\\
\beta_{R}S_{1}^{-1}\tau_{2R}\mathcal{R}_{\mathsf{x,z};\ell_{1}}\tau_{-2R}S_{1}\beta_{R}+\gamma_{R}S_{2}^{-1}\tau_{-2R}\mathcal{R}_{\mathsf{z,y};\ell_{2}}\tau_{2R}S_{2}\gamma_{R}&\ (-R\le t\le R)\\
\tau_{-2R}\mathcal{R}_{\mathsf{z,y};\ell_{2}}\tau_{2R}&\ (t\ge R),
\end{array}\right.
\end{align*}
where for $a\in\R$, $\tau_{a}$ is the time shift operator defined by $\tau_{a}\ell(t)=\ell(t+a)$, $S_{1}$ and $S_{2}$ are parallel translations along paths $[0,1]\ni\tau\mapsto\exp_{\phi_{\mathsf{z}}}(\tau X_{1,+}(t+2R)+(1-\tau)X_{1,2,R;t})$ and $[0,1]\ni\tau\mapsto\exp_{\phi_{\mathsf{z}}}(\tau X_{2,-}(t-2R)+(1-\tau)X_{1,2,R;t})$, respectively. Observe that $\mathcal{R}_{\mathsf{x,z,y};\ell_{1}\#_{R}\ell_{2}}$ defines a bounded linear operator between $L^{2}(\R,(\ell_{1}\#_{R}\ell_{2})^{\ast}T\mathcal{F}^{1,1/2}(S^{1},N))$ and $W^{1,2}(\R,(\ell_{1}\#_{R}\ell_{2})^{\ast}T\mathcal{F}^{1,1/2}(S^{1},N))$.
Moreover, form the definition of $\mathcal{R}_{\mathsf{x,z,y};\ell_{1}\#_{R}\ell_{2}}$ it can be easily checked that the operator norm $$\|\mathcal{R}_{\mathsf{x,z,y};\ell_{1}\#_{R}\ell_{2}}\|_{\text{op}(L^{2}(\R,(\ell_{1}\#_{R}\ell_{2})^{\ast}T\mathcal{F}^{1,1/2}),W^{1,2}(\R,(\ell_{1}\#_{R}\ell_{2})^{\ast}T\mathcal{F}^{1,1/2}))}$$ is uniformly bounded for large $R>0$ and
\begin{equation}
\|d\mathcal{F}_{\mathsf{x,y}}(\ell_{1}\#_{R}\ell_{2})\circ\mathcal{R}_{\mathsf{x,z,y};\ell_{1}\#_{R}\ell_{2}}-\bm{1}_{L^{2}(\R,(\ell_{1}\#_{R}\ell_{2})^{\ast}T\mathcal{F}^{1,1/2})}\|_{\text{op}(L^{2}(\R,(\ell_{1}\#_{R}\ell_{2})^{\ast}T\mathcal{F}^{1,1/2}))}\to 0
\end{equation}
as $R\to\infty$. 

By (3.115), $d\mathcal{F}_{\mathsf{x,y}}(\ell_{1}\#_{R}\ell_{2}):W^{1,2}(\R,(\ell_{1}\#_{R}\ell_{2})^{\ast}T\mathcal{F}^{1,1/2}(S^{1},N))\to L^{2}(\R,(\ell_{1}\#_{R}\ell_{2})^{\ast}T\mathcal{F}^{1,1/2}(S^{1},N))$ is onto for all large $R>0$ and there exists a right inverse $$\mathcal{R}_{\ell_{1},\ell_{2};R}:L^{2}(\R,(\ell_{1}\#_{R}\ell_{2})^{\ast}T\mathcal{F}^{1,1/2}(S^{1},N))\to W^{1,2}(\R,(\ell_{1}\#_{R}\ell_{2})^{\ast}T\mathcal{F}^{1,1/2}(S^{1},N))$$
 whose operator norm is uniformly bounded for all large $R>0$.

We next observe that the differential of the map
\begin{align*}
H^{1}(S^{1},(\phi_{1}\#_{R}\phi_{2})^{\ast}TN)&\times H^{1/2}(S^{1},\Spin(S^{1})\otimes(\phi_{1}\#_{R}\phi_{2})^{\ast}TN)\ni(X,\xi)\\
&\mapsto(\exp_{\phi_{1}\#_{R}\phi_{2}}(X),S(X)(\psi_{1}\#_{R}\psi_{2}+\xi))\in\mathcal{F}^{1,1/2}(S^{1},N)
\end{align*}
at $(X,\xi)=(0,0)$ is given by
\begin{align}
&H^{1}(S^{1},(\phi_{1}\#_{R}\phi_{2})^{\ast}TN)\times H^{1/2}(S^{1},\Spin(S^{1})\otimes(\phi_{1}\#_{R}\phi_{2})^{\ast}TN)\ni\begin{pmatrix}X\\ \xi\end{pmatrix}\notag\\
&\mapsto\begin{pmatrix}1&0\\ \Gamma(\ell_{1}\#_{R}\ell_{2})&1\end{pmatrix}\begin{pmatrix}X\\ \xi\end{pmatrix}\in H^{1}(S^{1},(\phi_{1}\#_{R}\phi_{2})^{\ast}TN)\times H^{1/2}(S^{1},\Spin(S^{1})\otimes(\phi_{1}\#_{R}\phi_{2})^{\ast}TN)\notag\\
&\cong T_{\ell_{1}\#_{R}\ell_{2}}\mathcal{F}^{1,1/2}(S^{1},N),
\end{align}
where $\Gamma(\ell_{1}\#_{R}\ell_{2}):H^{1}(S^{1},(\phi_{1}\#_{R}\phi_{2})^{\ast}TN)\to H^{1/2}(S^{1},\Spin(S^{1})\otimes(\phi_{1}\#_{R}\phi_{2})^{\ast}TN)$ is given by
$$\Gamma(\ell_{1}\#_{R}\ell_{2})[X]=-(\psi_{1}\#_{R}\psi_{2})^{k}X^{j}\Gamma_{jk}^{i}(\phi_{1}\#_{R}\phi_{2})\otimes\frac{\partial}{\partial y^{i}}(\phi_{1}\#_{R}\phi_{2})$$
and the identification $T_{\ell_{1}\#_{R}\ell_{2}}\mathcal{F}^{1,1/2}(S^{1},N)\cong H^{1}(S^{1},(\phi_{1}\#_{R}\phi_{2})^{\ast}TN)\times H^{1/2}(S^{1},\Spin(S^{1})\otimes(\phi_{1}\#_{R}\phi_{2})^{\ast}TN)$ is as in (3.2) (see also~\cite[(3.17), (3.19)]{I}). In fact, (3.116) is the differential of the local coordinate map of $\mathcal{F}^{1,1/2}(S^{1},N)$ at $\ell_{1}\#_{R}\ell_{2}$.

We then find a genuine solution to the equation $\mathcal{F}_{\mathsf{x,y}}(\ell)=0$ in the following form, $$\ell=\ell_{1,2;R}:=(\exp_{\phi_{1}\#_{R}\phi_{2}}(X),S(X)(\psi_{1}\#_{R}\psi_{2}+\xi)),$$
 where $(X,\xi)\in W^{1,2}(\R,H^{1}(S^{1},(\phi_{1}\#_{R}\phi_{2})^{\ast}TN)\times H^{1/2}(S^{1},\Spin(S^{1})\otimes(\phi_{1}\#_{R}\phi_{2})^{\ast}TN))$ takes the form
\begin{equation}\begin{pmatrix}1&0\\ \Gamma(\ell_{1}\#_{R}\ell_{2})&1\end{pmatrix}\begin{pmatrix}X\\ \xi\end{pmatrix}=\mathcal{R}_{\ell_{1},\ell_{2};R}(\eta)
\end{equation}
for some $\eta\in L^{2}(\R,(\ell_{1}\#_{R}\ell_{2})^{\ast}T\mathcal{F}^{1,1/2}(S^{1},N))$.

Note that (3.117) implies that 
\begin{align}
&\quad\|X\|_{W^{1,2}(\R,H^{1}(S^{1},(\phi_{1}\#_{R}\phi_{2})^{\ast}TN)}+\|\xi\|_{H^{1/2}(S^{1},\Spin(S^{1})\otimes(\phi_{1}\#_{R}\phi_{2})^{\ast}TN))}\notag\\
&\le C\|\eta\|_{L^{2}(\R,(\ell_{1}\#_{R}\ell_{2})^{\ast}T\mathcal{F}^{1,1/2})}
\end{align}
for some $C>0$ independent of large $R>0$.

By the Taylor expansion, when $X,\xi$ is given by (3.117) we have 
\begin{align*}
\mathcal{F}(X,\xi):&=\mathcal{F}_{\mathsf{x,y}}(\exp_{\phi_{1}\#_{R}\phi_{2}}(X),S(X)(\psi_{1}\#_{R}\psi_{2}+\xi))\notag\\
&=\mathcal{F}_{\mathsf{x,y}}(\ell_{1}\#_{R}\ell_{2})+d\mathcal{F}_{\mathsf{x,y}}(\ell_{1}\#_{R}\ell_{2})\begin{pmatrix}1&0\\ \Gamma(\ell_{1}\#_{R}\ell_{2})&1\end{pmatrix}\begin{pmatrix}X\\ \xi\end{pmatrix}+r(X,\xi)\notag\\
&=\mathcal{F}_{\mathsf{x,y}}(\ell_{1}\#_{R}\ell_{2})+d\mathcal{F}_{\mathsf{x,y}}(\ell_{1}\#_{R}\ell_{2})\circ\mathcal{R}_{\ell_{1},\ell_{2};R}(\eta)+r(X,\xi)\notag\\
&=\mathcal{F}_{\mathsf{x,y}}(\ell_{1}\#_{R}\ell_{2})+\eta+r(X,\xi),
\end{align*}
where the remainder term $r(X,\xi)$ satisfies $r(X,\xi)=o(\|\eta\|_{L^{2}(\R,(\ell_{1}\#_{R}\ell_{2})^{\ast}T\mathcal{F}^{1,1/2})})$ by (3.118).

Since $\mathcal{F}_{\mathsf{x,y}}(\ell_{1}\#_{R}\ell_{2})\to 0$ as $R\to\infty$, by the standard Banach fixed point argument, there is a unique small $\eta\in L^{2}(\R,(\ell_{1}\#_{R}\ell_{2})^{\ast}T\mathcal{F}^{1,1/2}(S^{1},N))$ such that the equation $\mathcal{F}(X,\xi)=0$ holds for $X,\xi$ defined through (3.117). This completes the proof.
\hfill$\Box$

\section{Definition of the Morse-Floer homology under the transversality condition}
Throughout this section, we assume that the following conditions are satisfied for $\Lag_{H}$, where $H\in\BH^{3}_{p+1}$:
\begin{enumerate}[(i)]
\item $\Lag_{H}$ is a Morse function on $\mathcal{F}^{1,1/2}(S^{1},N)$. 
\item The negative gradient flow system $\frac{d\ell}{dt}=-\nabla_{1,1/2}\Lag_{H}(\ell)$ defines a Morse-Smale system.
This means that $0\in L^{2}(\R,\mathcal{C}^{0,\alpha}(M))$ is a regular value of the map $\ell\mapsto\frac{d\ell}{dt}+\nabla_{1,1/2}\Lag_{H}(\ell)$.
\end{enumerate}
Under the above assumption (i), (ii), for $p\in\Z$ we define
\begin{equation}
C_{p}(\Lag_{H}):=\bigoplus_{\mathsf{x}\in\text{crit}_{p}(\Lag_{H})}\Z\langle\mathsf{z}\rangle
\end{equation}
and
\begin{equation}
C_{p}(\Lag_{H};\Z_{2})=C_{p}(\Lag_{H})\otimes\Z_{2},
\end{equation}
where $\text{crit}_{p}(\Lag_{H})=\{\mathsf{x}\in\text{crit}(\Lag_{H}):d\Lag_{H}(\mathsf{x})=0,\,\mu_{H}(\mathsf{x})=p\}$. Thus $\{C_{p}(\Lag_{H};\Z_{2})\}_{p\in\Z}$ is a graded group with grading given by the relative Morse index.

For $\mathsf{x,y}\in\text{crit}(\Lag_{H})$ with $\mu_{H}(\mathsf{x})-\mu_{H}(\mathsf{y})=1$, under the above assumption (i), (ii) and Corollary 3.1, $\hat{\mathcal{M}}(\mathsf{x},\mathsf{y})$ is a finite set. For $p\in\Z$, we define
\begin{equation}
\partial_{p}(\Lag_{H})\langle\mathsf{x}\rangle=\sum_{\mathsf{y}\in\text{crit}_{p-1}(\Lag_{H})}n(x,y)\langle\mathsf{y}\rangle
\end{equation}
for a generator $\mathsf{x}\in\text{crit}(\Lag_{H})$, where $n(\mathsf{x},\mathsf{y})=\#\hat{\mathcal{M}}(\mathsf{x},\mathsf{y})\,(\text{mod 2})$, and extend linearly to define $\partial_{p}(\Lag_{H}):C_{p}(\Lag_{H};\Z_{2})\to C_{p-1}(\Lag_{H};\Z_{2})$.

We note that the sum in (4.3) is finite and (4.3) is well-defined. Indeed, for $\mathsf{y}\in\text{crit}(\Lag_{H})$ with $n(\mathsf{x},\mathsf{y})\ne 0$, there holds $\mathcal{M}(\mathsf{x},\mathsf{y})\ne\emptyset$ and there exists a flow line connecting $\mathsf{x}$ and $\mathsf{y}$. Since $\Lag_{H}$ is non-increasing along the negative gradient flow, we have $\Lag_{H}(\mathsf{y})\le\Lag_{H}(\mathsf{x})$. Then by the condition (i) and the Palais-Smale condition (see Proposition 3.6), the set $\{\mathsf{y}\in\text{crit}(\Lag_{H}):\Lag_{H}(\mathsf{y})\le\Lag_{H}(\mathsf{x})\}$ is a finite set. Thus the set of points $\mathsf{y}\in\text{crit}_{p-1}(\Lag_{H})$ with $n(\mathsf{x},\mathsf{y})\ne 0$ is a finite set for any $\mathsf{x}\in\text{crit}_{p}(\Lag_{H})$ and the sum in (4.3) is finite.
By the standard argument, we have

\begin{proposition}
$\{(C_{p}(\Lag_{H};\Z_{2}),\partial_{p}(\Lag_{H}))\}_{p\in\Z}$ is a chain complex, that is, the following holds
$$\partial_{p-1}(\Lag_{H})\circ\partial_{p}(\Lag_{H})=0$$
for all $p\in\Z$.
\end{proposition}
\textit{Proof.} Let $\mathsf{x}\in\text{crit}_{p}(\Lag_{H})$ be arbitrary. We have
\begin{equation}
\partial_{p-1}\circ\partial_{p}\langle\mathsf{x}\rangle=\sum_{\mathsf{y}\in\text{crit}_{p-2}(\Lag_{H})}\sum_{\mathsf{z}\in\text{crit}_{p-1}(\Lag_{H})}n(\mathsf{x,z})n(\mathsf{z,y})\langle\mathsf{y}\rangle.
\end{equation}
For arbitrary $\mathsf{y}\in\text{crit}_{p-2}(\Lag_{H})$, as we have observed in Corollary 3.1, $\overline{\hat{\mathcal{M}}(\mathsf{x,y})}$ is a compact 1-dimensional manifold with boundary if it is not empty. The components of its boundary consists of broken flow lines with exactly one breaking, i.e., any component is a broken flow line of the form $\ell_{1}(\overline{\R})\cup\ell_{2}(\overline{\R})$, where $\ell_{1}\in\mathcal{M}(\mathsf{x,z})$, $\ell_{2}\in\mathcal{M}(\mathsf{z,y})$ for some $\mathsf{z}\in\text{crit}_{p-1}(\Lag_{H})$. The converse is also true by the gluing result Proposition 3.8. Thus the number $\sum_{\mathsf{z}\in\text{crit}_{p-1}(\Lag_{H})}n(\mathsf{x,z})n(\mathsf{z,y})$ (mod 2) is the number (mod 2) of the connected components of the boundary of $\overline{\hat{\mathcal{M}}(\mathsf{x,y})}$. Since the number of the components of 1-dimensional manifold is even, (4.4) is 0 (mod 2). This completes the proof.
\hfill$\Box$

\medskip

\begin{definition} Given $p\ge 3$. For $H\in\BH^{3}_{p+1}$ satisfying (i) and (ii) above, we denote by $HF_{\ast}(\Lag_{H},\mathcal{F}^{1,1/2}(S^{1},N);\Z_{2})$ the homology of the chain complex $\{(C_{p}(\Lag_{H}),\partial_{p}(\Lag_{H})\}_{p\in\Z}$:
$$HF_{p}(\Lag_{H},\mathcal{F}^{1,1/2}(S^{1},N);\Z_{2}):=\frac{\ker\partial_{p}(\Lag_{H})}{\rm{Im}\,\partial_{p+1}(\Lag_{H})}.$$
\end{definition}

\section{Transversality}
\subsection{The Morse property is generic}

First, we need to construct a separable Banach space of perturbations. So we let $\{p_{i}\}_{i\in \N}$ be a dense family of points in $N$ and similarly, for each $p_{i}$ we consider a dense family of vectors  $\{z_{j}^{i}\}_{j\in \N}\subset T_{p_{i}}N$. For a fixed $r>0$ less than the injectivity radius of $N$, we can identify points in $B_{r}(p_{i})$ with vectors in $T_{p_{i}}N$, that is if $p\in B_{r}(p_{i})$, $v_{p}$ is uniquely determined by $p=exp_{p_{i}}(v_{p})$. So we define now the functions $V_{i,j}:N\to \R$ by
\begin{equation}
V_{i,j}(p)=\left\{\begin{array}{ll}
\rho(|v_{p}|^{2})g_{p_{i}}(v_{p},z_{i}^{j}) \quad\text{ if }\quad |v_{p}|<\frac{r}{2}\\
0 \quad\text{otherwise}
\end{array}
\right.,
\end{equation}
where  $|\cdot|=g(\cdot,\cdot)^{\frac{1}{2}}$ and $\rho$ is a cut-off function supported in $[0,\frac{r^{2}}{2}]$, $\rho=1$ on $[0,\frac{r^{2}}{4}]$. Similarly, we choose a dense family of spinors $\{\varphi_{i}^{j}\}_{j\in \N}\in \mathbb{S}^{1}\otimes T_{p_{i}}N$ and construct the functions $W_{i,j}^{k}$ defined by
\begin{equation}
W_{i,j}^{k}(p,\psi)=\left\{\begin{array}{ll}
\rho(|\psi(p)|^{2})\rho(|v_{p}|^{2})\langle \psi(p),\varphi_{i}^{j}\rangle \quad\text{ if }\quad |v_{p}|<\frac{r}{2} \text{ and } |\psi(p)|<\frac{R_{k}}{2}\\
0 \quad\text{otherwise}
\end{array}
\right.,
\end{equation}
where $\{R_{k}\}_{k\in\N}$ is an increasing dense sequence of positive rational numbers. We define the space of perturbations $\mathbb{H}^{3}_{b}$ by $U\in\mathbb{H}^{3}_{b}$ if and only if there exist two sequences of numbers $\{c_{ij}\}$ and $\{d_{ijk}\}$, a function $\beta \in C^{3}(S^{1},\R)$ so that $$U(t,p,\psi)=\beta(t)\sum_{i,j,k}c_{ij}V_{i,j}(p)+d_{ijk}W_{i,j}^{k}(p,\psi),$$
with $$\|U\|_{3,b}:=\|\beta\|_{C^{3}}+\sum_{i,j,k}|c_{i,j}|\|V_{i,j}\|_{C^{3}}+|d_{ijk}|\|W_{i,j}^{k}\|_{C^{3}}<\infty.$$
One can easily see that $\mathbb{H}^{3}_{b}$ is a separable Banach space with the norm $\|\cdot\|_{3,b}$ since it is isomorphic to $C^{3}(S^{1})\times \ell^{1}$. Now we can state the following:
\begin{proposition}
There exists a residual set $\mathbb{H}_{reg}\subset\mathbb{H}^{3}_{b}$ such that $\Lag_{H+h}$ is a Morse function on $\mathcal{F}^{1,1/2}(S^{1},N)$ for every $h\in \mathbb{H}_{reg}$.
\end{proposition}

Indeed, we consider the map $\mathfrak{G} : \mathcal{F}^{1,1/2}(S^{1},N)\times \mathbb{H}_{b}^{3} \to T\mathcal{F}^{1,1/2}(S^{1},N)$ defined by
$$\mathfrak{G}(\phi,\psi,h)=\nabla_{1,1/2}\mathcal{L}_{H+h}(\phi,\psi).$$
The following lemma holds:
\begin{lemma}
Fix $(\phi_{0},\psi_{0},h)\in \mathfrak{G}^{-1}(0)$, then the map $$d\mathfrak{G}(\phi_{0},\psi_{0},h) : T_{(\phi_{0},\psi_{0})}\mathcal{F}^{1,1/2}(S^{1},N)\times \mathbb{H}_{b}^{3}\to T_{(\phi_{0},\psi_{0})}\mathcal{F}^{1,1/2}(S^{1},N)$$ is surjective.
\end{lemma}
\textit{Proof:}
The proof follows closely the idea in \cite{W}. In order to do this, we will show that $\text{Ran}(d\mathfrak{G}(\phi_{0},\psi_{0},h))$ is dense, which is equivalent to showing that
$$\text{Ran}(d\mathfrak{G}(\phi_{0},\psi_{0},h))^{\perp}=\{0\}$$
So we take $Z=(Z_{1},Z_{2})\in\text{Ran}(d\mathfrak{G}(\phi_{0},\psi_{0},h))^{\perp}$, then
$$\langle Z, d_{(\phi,\psi)}\mathfrak{G}(\phi_{0},\psi_{0},h) \xi \rangle =0$$
for all $\xi \in  T_{(\phi_{0},\psi_{0})}\mathcal{F}^{1,2}(S^{1},N)$ and
$$\langle Z, d_{h}\mathfrak{G}(\phi_{0},\psi_{0},h)V \rangle =0$$
for all $V\in \mathbb{H}_{b}^{3}$. Now since $H\in \mathbb{H}^{3}_{p+1}$ and $h\in \mathbb{H}^{3}_{b}$, it is easy to see from the first equation that $Z$ is a solution to a linear system of differential equations with coefficients in $C^{1}$ and thus $Z$ is $C^{3}$. Now we have
$$d_{h}\mathfrak{G}(\phi_{0},\psi_{0},h)V=\left[\begin{array}{ll}
-(-\Delta+1)^{-1}\nabla_{\phi}V
\\
-(|\D|+1)^{-1}\nabla_{\psi}V
\end{array}
\right].$$

Assume first that there exists $t_{0}\in S^{1}$ such that $Z_{1}(t_{0})\not=0$. We can then consider a small interval $(t_{0}-\delta, t_{0}+\delta)$ for $\delta>0$ and small, and $r>0$ and less than the injectivity radius, so that $\phi_{0}(t)=exp_{\phi_{0}(t_{0})}(\xi(t))$ and $|\xi(t)|<r/2$ for $t\in(t_{0}-\delta,t_{0}+\delta)$. Also, for $p\in B_{r}(\phi_{0}(t_{0}))$, we have the existence of a unique vector $v_{p}\in T_{\phi(t_{0})}N$ so that $p=exp_{\phi_{0}(t_{0})}(v_{p})$. By taking $\delta$ even smaller, we can assume by continuity that
$$g(Z_{1}(t),Z_{1}(t_{0}))>0, \text{ for } t\in (t_{0}-\delta, t_{0}+\delta).$$

Therefore, we can define now the following function 
$$V(t,p)=\left\{ \begin{array}{ll}
\gamma(t)\rho(|v_{p}|^{2}) g(v_{p},Z_{1}(t_{0}))\quad\text{if}\quad|v_{p}|<\frac{r}{2}\\
0\quad\text{otherwise}
\end{array}
\right.,$$
 $\gamma$ a cut-off function supported in $(t_{0}-\frac{\delta}{2}, t_{0}+\frac{\delta}{2})$. It is easy to see that 
\begin{align}
\langle Z,d_{h}\mathfrak{G}(\phi_{0},\psi_{0},h)V\rangle&=-\langle Z_{1},(-\Delta+1)^{-1}\nabla_{\phi}V\rangle_{H^{1}} \notag \\
&=-\langle Z_{1},\nabla_{\phi}V\rangle_{L^{2}}\notag\\
&=-\int_{t_{0}-\frac{\delta}{2}}^{ t_{0}+\frac{\delta}{2}}\gamma(t)g(Z_{1}(t),Z_{1}(t_{0}))dt <0 \notag
\end{align}
This leads to a contradiction and therefore $Z_{1}\equiv 0$. This would lead to the conclusion if $V$ was in $\mathbb{H}_{b}^{3}$, which might not be the case, so by taking $p_{i}$ arbitrarily close to $\phi_{0}(t_{0})$ and $z_{i}^{j}$ arbitrarily close to $Z_{1}(t_{0})$, the same conclusion would hold by replacing $V$ by $\gamma(t)V_{i,j}\in \mathbb{H}_{b}^{3}$. \\
Now we do a similar construction for $Z_{2}$. Indeed, if $Z_{2}(t_{0})\not=0$, one considers the function 
$$W(t,p,\psi)=\left\{ \begin{array}{ll}
\gamma(t)\rho(|\psi(p)|^{2})\rho(|v_{p}|^{2}) \langle \psi (p), Z_{2}(t_{0})\rangle\quad\text{if}\quad|\psi(p)|<\frac{R}{2}\text{ and } |v_{p}|<\frac{r}{2}\\
0\quad\text{otherwise}
\end{array}
\right.,$$
where $\rho$ and $\gamma$ are as above with small $\delta>0$ and $R$ is fixed so that $|\psi_{0}(t)|<\frac{R}{2}$ for $t\in(t_{0}-\delta,t_{0}+\delta)$. Since $Z_{1}\equiv 0$ it is easy to see that
\begin{align*}\langle Z, d_{h}\mathfrak{G}(\phi_{0},\psi_{0},h)W\rangle&=-\langle Z_{2}, (|\D|+1)^{-1}\nabla_{\psi}W\rangle_{H^{1/2}}\\
&=-\langle Z_{2},\nabla_{\psi}W\rangle_{L^{2}}\\
&=-\int_{t_{0}-\frac{\delta}{2}}^{ t_{0}+\frac{\delta}{2}}\gamma(t)\langle Z_{2}(t),Z_{2}(t_{0})\rangle dt <0,
\end{align*}
which is another contradiction. But again, $W$ might not be in $\mathbb{H}_{b}^{3}$, so we make similar choices as in the $Z_{1}$ case to get an arbitrarily close function $\gamma W_{i,j,k}$ and the conclusion still holds. Therefore, $Z=0$ and $\text{Ran}(d\mathfrak{G}(\phi_{0},\psi_{0},h))$ is dense. Since $d\mathfrak{G}$ is a Fredholm operator, it has closed range and hence it is surjective.
\hfill$\Box$

The result of Proposition 5.1 follows then from the classical Sard-Smale's theorem.

\subsection{The Morse-Smale property is generic}
In this part, we will prove that there is a generic set of perturbations of the metric in such a way that the gradient flow with respect to this perturbed metric satisfies the Morse-Smale property. We will closely follow the proof presented in \cite{I2} for the case of the non-linear Dirac equation and we will adapt it to our setting. We also refer the reader to a different global construction of a perturbation set in \cite{AM3}.

Since $\mathcal{F}^{1,\frac{1}{2}}(S^{1},N)$ is a Hilbert manifold, then it is separable, so we consider a dense sequence of points $\{p_{i}\}_{i\in \N}\subset \mathcal{F}^{1,\frac{1}{2}}(S^{1},N)$ and for each $p_{i}$ we consider a dense family of pairs $\{(u^{1}_{ij},u_{ij}^{2})\}_{j\in\N}\subset T_{p_{i}}\mathcal{F}^{1,1}(S^{1},N)$ and a family of linear forms on $T_{p_{i}}\mathcal{F}^{1,\frac{1}{2}}(S^{1},N)$, $\{(A_{ik},B_{ik})\}_{k\in\N}$. We fix, $r$, small enough, so that the parallel translation $P_{p,p_{i}}:T_{p}\mathcal{F}^{1,\frac{1}{2}}\to T_{p_{i}}\mathcal{F}^{1,\frac{1}{2}}$ is well defined for $p\in B_{r}(p_{i})$ and let $\rho:C^{\infty}_{c}([0,+\infty))$ compactly supported in $[0,r]$ and equals $\rho=1$ on  $[0,\frac{r}{2}]$. Then we can define the functions $K_{ijk}$ by

$$K_{ijk}(p)=\left \{\begin{array}{ll}
\rho(\|p_{i}-p\|)P_{p_{i},p}\Big(A_{ik}(P_{p,p_{i}}\cdot)u_{ij}^{1}+B_{ik}(P_{p,p_{i}}\cdot)u_{ij}^{2}\Big) \quad\text{ if }\quad p\in B_{\frac{r}{2}}(p_{i})\\
0 \text{ otherwise }
\end{array}
\right.
$$
where $\|p_{i}-p\|$ means the distance between $p_{i}$ and $p$ in $\mathcal{F}^{1,\frac{1}{2}}$. Without loss of generality, we can assume that the $p_{i}$ and $u_{ij}$ are in $C^{2,\alpha}$. We consider then the space $\mathcal{K}$ defined by $K\in \mathcal{K}$ if and only if $K$ is symmetric and there exist a sequence of numbers $\{c_{ijk}\}\in\R$ such that $K=\sum_{i,j,k} c_{ijk}K_{ijk}$ 
$$\|K\|_{C^{2}}:=\sum_{i,j,k}|c_{ijk}|\|K_{i,j,k}\|_{C^{2}_{b}(\mathcal{F}^{1,2},\text{Sym}(T\mathcal{F}^{1,\frac{1}{2}},T\mathcal{F}^{1,1}))}<\infty,$$
where $T\mathcal{F}^{1,1}(S^{1},N)\to\mathcal{F}^{1,1/2}(S^{1},N)$ is the bundle over $\mathcal{F}^{1,1/2}(S^{1},N)$ whose fiber at $(\phi,\psi)\in\mathcal{F}^{1,1/2}(S^{1},N)$ is defined by
$$T_{(\phi,\psi)}\mathcal{F}^{1,1}(S^{1},N)=H^{1}(S^{1},\phi^{\ast}TN)\times H^{1}(S^{1},\Spin(S^{1})\otimes\phi^{\ast}TN).$$ 
Here, $\text{Sym}(T\mathcal{F}^{1,1/2},\mathcal{F}^{1,1})\to\mathcal{F}^{1,1/2}$ is the bundle of symmetric homomorphisms between $T\mathcal{F}^{1,1/2}$ and $T\mathcal{F}^{1,1}$ and $\|\,\cdot\,\|_{\text{op}(E,F)}$ denotes the operator norm on $L(E,F)$, the set of bounded linear operators between Banach spaces $E,F$.
``Symmetric'' here means that
$$\langle K(p)X,Y\rangle_{T_{p}\mathcal{F}^{1,1/2}}=\langle X,K(p)Y\rangle_{T_{p}\mathcal{F}^{1,1/2}}$$
for any $p\in\mathcal{F}^{1,1/2}$, $X\in T_{p}\mathcal{F}^{1,1/2}$ and $Y\in T_{p}\mathcal{F}^{1,1}$.

Given a continuous map $\theta : \mathcal{F}^{1,\frac{1}{2}}(S^{1},N)\to [0,1]$, the space of perturbation that will be considered is then
$$\mathbb{K}^{2}_{\theta}=\{K\in \mathcal{K}; \text{ there exists } C>0; \|K(p)\|_{op}\leq C\theta(p)\}.$$
If we endow it with the norm 
$$\|K\|_{\mathbb{K}_{\theta}^{2}}=\|K\|_{C^{2}(\mathcal{F}^{1,\frac{1}{2}}}+\sup_{p\in \mathcal{F}^{1,\frac{1}{2}}}\frac{\|K(p)\|_{op}}{\theta(p)}.$$
It can be easily checked, that the space $\mathbb{K}^{k}_{\theta}$ is a Banach space.

  
Now given $\rho>0$ we consider the ball $\mathbb{K}^{k}_{\theta,\rho}=\{K\in\mathbb{K}^{k}_{\theta};\|K\|_{\mathbb{K}^{k}_{\theta}}<\rho\}$. For $K\in\mathbb{K}^{k}_{\theta,\rho}$ we define the metric $G^{K}$ on $\mathcal{F}^{1,1/2}(S^{1},N)$ by
$$G^{K}_{\phi,\psi}(X,Y)=\langle(1+K(\phi,\psi))^{-1}X,Y\rangle_{H^{1}\times H^{1/2}}$$
This metric is well-defined for $\rho<\frac{1}{2S_{0}}$, where $S_{0}$ is the best constant of the embedding $H^{1}\times H^{1}$ in $H^{1}\times H^{1/2}$.
The following is the main result of this section.

\begin{proposition}
Assume that $\Lag_{H}$ is a Morse function on $\mathcal{F}^{1,1/2}(S^{1},N)$. Assume in addition that $\theta$ satisfies the following conditions:\\
1) $\theta(p)=0$ for all $p\in crit(\Lag_{H})$.\\
2) The zero set of $\theta$ is the closure of an open set..\\
3) If $W^{u}_{-\nabla \Lag_{H}}(x)$ and $W^{s}_{-\nabla \Lag_{H}}(y)$ intersect non-transversally at a point $p$ for some $x,y\in crit(\Lag_{H})$, then $\theta$ is not identically zero on the flow line through $p$.

Then there exist $0<\rho_{0}<\frac{1}{2S_{0}}$ and a residual set $\mathbb{K}_{reg}\in \mathbb{K}^{2}_{\theta,\rho_{0}}$ such that for all $K\in \mathbb{K}_{reg}$ the negative gradient flow with respect to the metric $G^{K}$ satisfies the Morse-Smale property up to order $2$.
\end{proposition}

Where $W^{u}_{-\nabla \Lag_{H}}(\mathsf{x})$ and $W^{s}_{-\nabla \Lag_{H}}(\mathsf{y})$ in the above proposition are unstable and stable manifolds of critical points $\mathsf{x,y}\in\text{crit}(\Lag_{H})$ of the vector field $-\nabla_{1,1/2}\Lag_{H}$, respectively:
$$W^{u}_{-\nabla \Lag_{H}}(\mathsf{x}):=\{c\in\mathcal{F}^{1,1/2}:\text{$\ell(t,c)$ exists for $t\le 0$ and $\lim_{t\to-\infty}\ell(t,c)=\mathsf{x}$}\},$$
$$W^{s}_{-\nabla \Lag_{H}}(\mathsf{y}):=\{c\in\mathcal{F}^{1,1/2}:\text{$\ell(t,c)$ exists for $t\ge 0$ and $\lim_{t\to+\infty}\ell(t,c)=\mathsf{y}$}\}$$
and $\ell(t,c)$ denotes the negative gradient flow with $\ell(0,c)=c$. The proof of this Proposition will be made through several steps that we will set as Lemmata.

\begin{lemma}
There exists $0<\delta\le\frac{1}{2S_{0}}$ such that for any $0<\rho_{0}<\delta$ and any $K\in \mathbb{K}^{2}_{\theta,\rho_{0}}$ we have: for $(\phi(t),\psi(t))\in C^{1}(\R,\mathcal{F}^{1,1/2})$ a solution to the $G^{K}$-negative gradient flow equation
\begin{equation}\label{GKF}
\frac{d}{dt}\begin{pmatrix}\phi\\ \psi\end{pmatrix}=-(1+K(\phi,\psi))\nabla_{1,1/2}\Lag_{H}(\phi,\psi)
\end{equation}
which satisfies the condition $\sup_{t\in\R}|\Lag_{H}(\phi(t),\psi(t))|=:C_{0}<+\infty$. Then there exists $C(C_{i},\rho_{0})>0$ such that
$$\sup_{t\in\R}\|\partial_{s}\phi(t)\|_{L^{2}(S^{1})}+\sup_{t\in\R}\|\psi(t)\|_{H^{1/2}(S^{1})}\le C.$$
\end{lemma}
\textit{Proof.} Set $\ell(t)=(\phi(t),\psi(t))$. We write $K(\ell(t))=\begin{pmatrix}K_{11}(t)&K_{12}(t)\\
K_{21}(t)&K_{22}(t)\end{pmatrix}$
according to the canonical decomposition (3.1). 
We first notice that since the metric $G^{K}$ is uniformly equivalent to the standard one, we have a similar result to Lemma 3.1. That is, there exists a positive constant $C=C(C_{0},\rho_{0})$ depending only on $C_{0}$ and $\rho_{0}$ such that the following holds for all $\tau\in\R$:
\begin{equation}\label{estphi}
\int_{\tau}^{\tau+1}\|\phi(t)\|_{H^{1}(S^{1})}^{2}\,dt\le C,
\end{equation}
\begin{equation}
\int_{\tau}^{\tau+1}\|\psi(t)\|_{H^{1/2}(S^{1})}\,dt\le C,
\end{equation}
\begin{equation}\label{estpsi}
\int_{\tau}^{\tau+1}\|\psi(t)\|^{p+1}_{L^{p+1}(S^{1})}\,dt\le C.
\end{equation}
Now following the steps in Proposition 3.2, we have
\begin{align}
\partial_{t}\psi(t)&=-(1+K_{22})(t)(1+|\D|)^{-1}(\D\psi-A(\partial_{s}\phi(t),\partial_{s}\cdot\psi(t))-\nabla_{\psi}H(s,\phi(t),\psi(t)))-\notag \\
&\quad-K_{21}\nabla_{\phi}\Lag_{H}(\phi(t),\psi(t))\notag\\
&=-(1+K_{22})(1+|\D|)^{-1}((\D+\lambda)\psi(t)-(\lambda\psi(t)+A(\partial_{s}\phi(t),\partial_{s}\cdot\psi(t))\notag\\
&\quad+\nabla_{\psi}H(s,\phi(t),\psi(t)))-K_{21}\nabla_{\phi}\Lag_{H}(\phi(t),\psi(t))\notag\\
&=-\mathsf{L}_{\lambda}\psi(t)-K_{22}(t)\mathsf{L}_{\lambda}\psi(t)+(1+K_{22})(1+|\D|)^{-1}(\lambda\psi(t)+A(\partial_{s}\phi(t),\partial_{s}\cdot\psi(t))+\notag \\
&\quad+\nabla_{\psi}H(s,\phi(t),\psi(t)))-K_{21}\nabla_{\phi}\Lag_{H}(\phi(t),\psi(t)).
\end{align}
Thus, if we consider 
\begin{align}
r(t)=&-K_{22}(t)\mathsf{L}_{\lambda}\psi(t)+(1+K_{22}(t))(1+|\D|)^{-1}(\lambda\psi(t)+A(\partial_{s}\phi(t),\partial_{s}\cdot\psi(t))\notag\\
&\quad+\nabla_{\psi}H(s,\phi(t),\psi(t)))\notag-K_{21}(t)\nabla_{\phi}\Lag_{H}(\phi(t),\psi(t)) \notag \\
=& -K_{22}(t)\mathsf{L}_{\lambda}\psi(t)+(1+K_{22}(t))(1+|\D|)^{-1}(\lambda\psi(t)+A(\partial_{s}\phi(t),\partial_{s}\cdot\psi(t))\notag\\
&\quad+\nabla_{\psi}H(s,\phi(t),\psi(t)))\notag+K_{21}(t)(1+K_{11}(t))^{-1}K_{12}(t)\nabla_{\psi}\Lag_{H}(\phi(t),\psi(t))\notag\\
&\quad+K_{21}(t)(1+K_{11}(t))^{-1}\partial_{t}\phi(t).
\end{align}
By using the fundamental solution $G_{\lambda}$ (see (3.33)), the gradient flow equation is written as
\begin{equation}
\psi(t)=\int_{\R}G_{\lambda}(t-\tau)r(\tau)\,d\tau.
\end{equation}
Therefore, we have as in (3.36)
$$\|\psi\|_{H^{1/2}}\leq C \sup_{\eta\in \mathbb{R}}\int_{\eta}^{\eta+1}\|r(\tau)\|_{H^{1/2}}\,d\tau.$$
Hence, it remains to estimate $\|r(s)\|_{H^{1/2}}$. But
\begin{align}
r(t)=&(K_{21}(1+K_{11})^{-1}K_{12}-K_{22})(t) \nabla_{\psi}\mathcal{L}_{H}+(1+|\D|)^{-1}(\lambda\psi(t)+A(\partial_{s}\phi(t),\partial_{s}\cdot\psi(t))+\notag \\
&+\nabla_{\psi}H(s,\phi(t),\psi(t)))+K_{21}(1+K_{11})^{-1}(t)\partial_{t}\phi(t),
\end{align}
and estimate as in (3.40) and (3.41)
\begin{align}
&\|(K_{21}(1+K_{11})^{-1}K_{12}-K_{22})(t) \nabla_{\psi}\mathcal{L}_{H}\|_{H^{1/2}}\notag\\
&\quad\leq C\|(K_{21}(1+K_{11})^{-1}K_{12}-K_{22})(t) \nabla_{\psi}\mathcal{L}_{H}\|_{H^{1}} \notag \\
&\quad\leq C\rho_{0}\theta(\phi(t),\psi(t))\|\nabla_{\psi}\mathcal{L}_{H}\|_{H^{1/2}}\notag \\
&\quad\leq C(1+\|\psi(t)\|_{H^{1/2}(S^{1})}+\|\partial_{s}\phi(t)\|_{L^{2}(S^{1})}\|\psi(t)\|_{H^{1/2}(S^{1})}+\|\psi(t)\|_{L^{p+1}(S^{1})}^{p}).
\end{align}
 Therefore, we obtain, by the H\"older's inequality and the fact that $\int_{-\infty}^{+\infty}\|\partial_{t}\phi(t)\|_{H^{1}}^{2}dt \leq C$,
\begin{align}
\|\psi(t)\|_{H^{1/2}(S^{1})}&\le C\sup_{\eta\in\R}\int_{\eta}^{\eta+1}(1+\|\psi(\tau)\|_{H^{1/2}(S^{1})}+\|\partial_{s}\phi(\tau)\|_{L^{2}(S^{1})}\|\psi(\tau)\|_{H^{1/2}(S^{1})}+\|\psi(\tau)\|_{L^{p+1}(S^{1})}^{p})\,d\tau\notag\\
&\le C\sup_{\eta\in\R}\biggl(1+\biggl(\int_{\tau}^{\tau+1}\|\psi(\tau)\|_{H^{1/2}(S^{1})}^{2}\,d\tau\biggr)^{\frac{1}{2}}+\notag\\
&\qquad+\biggl(\int_{\eta}^{\eta+1}\|\partial_{s}\phi(\tau)\|_{L^{2}(S^{1})}^{2}\,d\tau\biggr)^{\frac{1}{2}}\biggl(\int_{\eta}^{\eta+1}\|\psi(\tau)\|_{H^{1/2}(S^{1})}^{2}\,d\tau\biggr)^{\frac{1}{2}}+\notag \\
&\qquad+\biggl(\int_{\eta}^{\eta+1}\|\psi(\tau)\|_{L^{p+1}(S^{1})}^{p+1}\,d\tau\biggr)^{\frac{p}{p+1}}\biggr).
\end{align}
By combining (\ref{estphi}), (5.3), (\ref{estpsi}) we have,
\begin{equation}
\sup_{t\in\R}\|\psi(t)\|_{H^{1/2}(S^{1})}\le C(C_{0},\rho_{0})
\end{equation}
This gives the estimate for $\psi(t)$.\\

To estimate $\sup_{t\in\R}\|\partial_{s}\phi(t)\|_{L^{2}(S^{1})}$, we first have from (3.21) and (5.11),
\begin{equation}
\|\partial_{s}\phi(t)\|_{L^{2}(S^{1})}^{2}\le C(1+\|\partial_{t}\psi(t)\|_{H^{1/2}(S^{1})}).
\end{equation}
On the other hand, we have
\begin{align}
\partial_{t}\psi(t)&=-K_{21}(t)\nabla_{\phi}\Lag_{H}(\phi(t),\psi(t))-(1+K_{22})(t)\nabla_{\psi}\Lag_{H}(\phi(t),\psi(t)).
\end{align}
Again, since the metric $G^{K}$ is uniformly equivalent to the standard metric, we have
\begin{align}
\|\partial_{t}\psi(t)\|_{H^{1/2}(S^{1})}&\le C\Big (\rho_{0}\|\nabla_{\phi}\Lag_{H}(\phi(t),\psi(t))\|_{H^{1}}+\|\psi(t)\|_{H^{1/2}(S^{1})}\notag\\
&\quad+\|(1+|\D|)^{-1}A(\partial_{s}\phi(t),\partial_{s}\cdot\psi(t))\|_{H^{1/2}(S^{1})}\notag\\
&\quad + \|(1+|\D|)^{-1}\nabla_{\psi}H(s,\phi(t),\psi(t))\|_{H^{1/2}(S^{1})}\Big)\notag\\
&\le C\Big(1+\rho_{0}\|\nabla_{\phi}\Lag_{H}(\phi(t),\psi(t))\|_{H^{1}}+\|\psi(t)\|_{H^{1/2}(S^{1})}\notag\\
&\quad+\|\partial_{s}\phi(t)\|_{L^{2}(S^{1})}\|\psi(t)\|_{H^{1/2}(S^{1})}+\|\psi(t)\|_{L^{p+1}(S^{1})}^{p}\Big).
\end{align}
On the other hand, $\|\nabla_{\phi}\Lag_{H}(\phi(t),\psi(t))\|_{H^{1}}$ is estimated by
\begin{align*}
\|\nabla_{\phi}\Lag_{H}(\phi(t),\psi(t))\|_{H^{1}}&\le\|(-\Delta+1)^{-1}\nabla_{s}\partial_{s}\phi\|_{H^{1}}+\|(-\Delta+1)^{-1}R(\phi)\langle\psi,\partial_{s}\phi\cdot\psi\rangle\|_{H^{1}}\\
&\quad+\|(-\Delta+1)^{-1}\nabla_{\phi}H(s,\phi,\psi)\|_{H^{1}}\\
&\le C\Big(\|\partial_{s}\phi(t)\|_{H^{1}}+\|\partial_{s}\phi(t)\|_{L^{2}}^{2}+\|\psi(t)\|_{H^{1/2}}^{2}\|\partial_{s}\phi(t)\|_{L^{2}}+\|\psi(t)\|_{H^{1/2}}^{q_{1}+1}\Big),
\end{align*}
where we have used (1.5) in the last step.

\noindent
Combining the last inequality with (5.12) and (5.14) and using (5.11), we have as in the proof of Proposition 3.2
$$\|\partial_{s}\phi(t)\|_{H^{1}}^{2}\le C\Big(1+\|\partial_{s}\phi(t)\|_{H^{1}}+\rho_{0}\|\partial_{s}\phi(t)\|_{H^{1}}^{2}\Big).$$
Thus, initially choosing $\rho_{0}$ small such that $C\rho_{0}<\frac{1}{2}$, we have a uniform bound $$\sup_{t\in\R}\|\partial_{s}\phi(t)\|_{L^{2}}\le C(C_{0},\rho_{0}).$$ 
This completes the proof.
\hfill$\Box$

\medskip

In the next proposition, we prove relative compactness of flow lines. The following is an analogue of Proposition 3.5. Note that the proof for Proposition 3.5 does not apply in the present case because the solutions of the perturbed gradient flow equation do not have enough regularity as in the case of unperturbed ones. So we give the alternative argument which is based on an idea of~\cite{AM1}.

\begin{proposition}
Let $0<\rho_{0}<\frac{1}{2S_{0}}$, then for every $K\in \mathbb{K}^{2}_{\theta,\rho_{0}}$, if $x_{+},x_{-}\in crit(\Lag_{H})$, then the moduli space $\mathcal{M}(x_{+},x_{-})\subset W^{1,2}_{\mathsf{{x}_{-},x_{+}}}(\R,\mathcal{F}^{1,1/2}(S^{1},N))$ is relatively compact.
\end{proposition}
The proof consists again of several steps.
\begin{lemma}
There exists $C>0$ such that for every $\ell\in \mathcal{M}(x_{+},x_{-})$, the length of $\ell$ is bounded by $C$. That is
$$L(\ell)=\int_{-\infty}^{+\infty}\|\ell'(t)\| dt\leq C.$$
\end{lemma}
{\it Proof:}
The proof is similar to the one in \cite{AM1}, but for the sake of completeness, we provide it here.
Starting from the fact that $\Lag_{H}$ is Morse and satisfies the (PS) condition, we have the existence of finitely many critical points $x_{1},..., x_{k}$ with energy between $\Lag_{H}(x_{+})$ and $\Lag_{H}(x_{-})$. We let $J=[\Lag_{H}(x_{-}),\Lag_{H}(x_{+})]$. Using the invertibility of $d^{2}\Lag_{H}(x_{i})$, we have the existence of $r_{0}>0$ small enough and $c\geq 1$, such that for $x\in \cup_{i=1}^{k}B_{r_{0}}(x_{i})$,
$$\|d^{2}\Lag_{H}(x)\|\leq c \text{ and } \|d^{2}\Lag_{H}(x)^{-1}\|\leq c.$$
Using a local frame around each $x_{i}$ we can assume that $\Lag_{H}$ is defined on a small ball centered at $x_{i}$ on a Hilbert space. Hence, we can apply the implicit function theorem to have that the map $\nabla_{1,2} \Lag_{H}^{-1}$ is well defined on a neighborhood of $B_{r_{0}}(0)$. Hence again by the implicit function theorem one has that
\begin{equation}\label{nab1}
\|\nabla_{1,1/2}\Lag_{H}(x)\|\geq \frac{1}{c}\|x-x_{i}\| \text{ for } x\in B_{r_{0}/2}(x_{i}).
\end{equation}
Using a Taylor's expansion for $\Lag_{H}$, we get for $x\in B_{r_{0}}(x_{i})$, 
\begin{equation}\label{nab2}
|\Lag_{H}(x)-\Lag_{H}(x_{i})|\leq \|\nabla_{1,1/2}\Lag_{H}(x)\|\|x-x_{i}\|+\frac{c}{2}\|x-x_{i}\|^{2}
\end{equation}
Combining $(\ref{nab1})$ and $(\ref{nab2})$, one has for $c_{1}=c+\frac{c^{3}}{2}$ and $x\in B_{r_{0}/2}(x_{i})$
\begin{equation}\label{nab3}
|\Lag_{H}(x)-\Lag_{H}(x_{i})|\leq c_{1}\|\nabla_{1,1/2}\Lag_{H}(x)\|^{2}.
\end{equation}
The second estimate on the gradient that we will need, comes from the (PS) condition. Indeed, there exists $\delta>0$ such that for $x\in \Lag_{H}^{-1}(J)\setminus \cup_{i=1}^{k}B_{r_{0}/2}(x_{i})$, we have
\begin{equation}\label{nab4}
\|\nabla_{G^{K}}\Lag_{H}(x)\|\geq \delta.
\end{equation}

We define now the $\mathcal{T}:J\to \R$ by
$$\mathcal{T}(w)=\inf\{\|\nabla_{G^{K}}\Lag_{H}(x)\|, \Lag_{H}(x)=w\}.$$
Then using $(\ref{nab3})$ and $(\ref{nab4})$, we have that
$$\mathcal{T}(w)\geq \min\{\frac{1}{\sqrt{c_{1}}}\min_{i=1,\cdots,k}|w-\Lag_{H}(x_{i})|^{\frac{1}{2}},\delta\}$$
Thus,
\begin{equation}\label{nab5}
\frac{1}{\mathcal{T}(w)}\leq \max\{\sqrt{c_{1}}\max_{i=1,\cdots,k}\frac{1}{|w-\Lag_{H}(x_{i})|^{\frac{1}{2}}},\frac{1}{\delta}\}.
\end{equation}
We deduce from this, that $\frac{1}{\mathcal{T}}\in L^{1}(J)$ and

\begin{align}
\int_{J}\frac{1}{\mathcal{T}(w)}dw&\leq \min\{\sqrt{8kc_{1}},|J|^{\frac{1}{2}}+2kc_{1} \frac{\delta}{|J|^{\frac{1}{2}}}\}|J|^{\frac{1}{2}}\notag \\
&\leq C|\Lag_{H}(x_{+})-\Lag_{H}(x_{-})|^{\frac{1}{2}},\notag
\end{align}
where $C=4c\sqrt{kc}+\frac{1}{\delta}|J|^{\frac{1}{2}}$. Now by construction, we have that $$\mathcal{T}(\Lag_{H}(x))\leq \|\nabla_{G^{K}}\Lag_{H}(x)\|,$$
Hence if $\ell \in \mathcal{M}(x_{+},x_{-})$, we have
\begin{align}
\|\ell'(t)\|&=-\frac{\langle \nabla_{G^{K}}\Lag_{H}(\ell(t)),\ell'(t)\rangle}{\|\nabla_{G^{K}}\Lag_{H}(x)\|}\notag \\
&\leq -\frac{\langle \nabla_{G^{K}}\Lag_{H}(\ell(t)),\ell'(t)\rangle}{\mathcal{T}(\ell(t))}.\notag
\end{align}
Using the substitution $w=\Lag_{H}(\ell(t))$, one has
$$L(\ell)=\int_{\R}\|\ell'(t)\|dt \leq \int_{J} \frac{1}{\mathcal{T}(w)}dw \leq C |\Lag_{H}(x_{+})-\Lag_{H}(x_{-})|^{\frac{1}{2}}.$$
\hfill$\Box$

\begin{lemma}
The set $\{\ell(t);\ell\in \mathcal{M}(x_{+},x_{-}), t\in \R\}$ is relatively compact.
\end{lemma}
{\it Proof:} Let $\ell=(\phi,\psi)\in \mathcal{M}(x_{+},x_{-})$. We consider first the $\psi$ part. As in the previous case, if we consider $r(t)$ defined by
\begin{align}
r(t)=&(K_{21}(1+K_{11})^{-1}K_{12}-K_{22})(t) \nabla_{\psi}\mathcal{L}_{H}+(1+|\D|)^{-1}(\lambda\psi(t)+A(\partial_{s}\phi(t),\partial_{s}\cdot\psi(t))+\notag \\
&+\nabla_{\psi}H(s,\phi(t),\psi(t)))-K_{21}(1+K_{11})^{-1}(t)\partial_{t}\phi(t)\notag\\
=&r_{1}(t)+r_{2}(t)
\end{align}
where $r_{1}(t)=(K_{21}(1+K_{11})^{-1}K_{12}-K_{22})(t) \nabla_{\psi}\mathcal{L}_{H}$.
Hence, 
$$\psi(t)=\int_{\R}G_{\lambda}(t-\tau)r(s)\,d\tau =\psi_{1}(t)+\psi_{2}(t).$$
Now the estimate on $\psi_{2}(t)$ is exactly similar to the one proved in Proposition 3.3. For $\psi_{1}$, we have

\begin{align}
\|\psi_{1}(t)\|_{H^{1}}&\leq \int_{\R}\|G_{\lambda}(t-\tau)\|_{H^{1}}\|K_{21}(1+K_{11})^{-1}K_{12}-K_{22}\|_{\text{op}(H^{1/2},H^{1})}\|\nabla_{\psi}\Lag_{H}(\phi,\psi\|_{H^{1/2}}\,d\tau\notag \\
&\leq C\rho_{0} \int_{\R}e^{-\kappa |\tau|}\left(1+\|\psi(t-\tau)\|_{H^{1/2}}+\|\partial_{s}\phi(t-\tau)\|_{L^{2}}\|\psi(t-\tau)\|_{H^{1/2}}+\right.\notag\\
&\left.\qquad +\|\psi(t-\tau)\|_{L^{p+1}}^{p}\right)\,d\tau,
\end{align}
and using Lemma 5.1 we have $\|\psi_{1}(t)\|_{H^{1}}\leq C$.\\
Now we want to show that the $\phi$ component is relatively compact. Recall that
$$\frac{d\ell}{dt}=-(1+K(\ell))\nabla_{1,1/2} \Lag_{H}(\ell)$$
therefore
\begin{align}
\frac{d\phi}{dt}&=-\nabla_{\phi}\Lag_{H}-K_{11}\nabla_{\phi}\Lag_{H}-K_{12}\nabla_{\psi}\Lag_{H}\notag\\
&=-\phi+U(t)\notag
\end{align}
where 
\begin{align}
U(t)&=(-\Delta +1)^{-1}\Big(\phi+\Gamma(\phi(t),\partial_{s}\phi(t),\partial_{s}\phi(t))-B(\phi,\psi) -\notag \\
&\quad-\frac{1}{2}R(\phi)\langle\psi,\partial_{s}\cdot\psi\rangle+\nabla_{\phi}H(s,\phi,\psi)\Big)-K_{11}\nabla_{\phi}\Lag_{H}(\phi,\psi)-K_{12}\nabla_{\psi}\Lag_{H}(\phi,\psi)\notag
\end{align}
First recall that $(\phi(t),\psi(t))\in H^{1}\times H^{\frac{1}{2}}$ are uniformly bounded. Now consider the operator
$$U_{1}(\phi,\psi)=  (-\Delta +1)^{-1}\Big(\phi+\Gamma(\phi, \partial_{s}\phi,\partial_{s}\phi)-B(\phi,\psi)-\frac{1}{2}R(\phi)\langle\psi,\partial_{s}\cdot\psi\rangle+\nabla_{\phi}H(s,\phi,\psi)\Big)$$
Clearly, $U_{1}$ is a compact operator and since, $\sup_{t\in \R}(\|\phi\|_{H^{1}}+\|\psi\|_{H^{\frac{1}{2}}})\leq C<+\infty,$, we have 
$$\{U_{1}(\ell(t)), \ell \in \mathcal{M}(x_{+},x_{-}), t\in \R\}$$
is relatively compact since it is the image of a bounded set via the operator $U_{1}$.\\
Also, we notice first the compactness of the operators $K_{11}$, $K_{12}$ (because they are finite rank locally and in fact they do have higher regularity by the choice of the $p_{i}$ and $u_{ij}$) and we have uniform boundedness of $\ell\in \mathcal{M}(x_{+},x_{-})$. Combining this, with the fact that the $K_{11}(\ell)$ and $K_{12}(\ell)$ are supported in a compact set of the form $[T_\ell,T_{\ell}]$, we conclude the relative compactness of the set
$$\{ K_{11}\nabla_{\phi}\Lag_{H}(\ell(t))-K_{12}\nabla_{\psi}\Lag_{H}(\ell(t)), \ell\in \mathcal{M}(x_{+},x_{-}), t\in [-T_{\ell},T_{\ell}]\}.$$
Now if $G_{1}$ is the fundamental solution of the operator $T$ defined by $T(\phi)=\frac{d\phi}{dt}+\phi$, we have that
$$\phi(t)=G_{1}\ast U(t)=T^{-1}(U(t)).$$
and since $T$ is continuous, we have that the set of such $\phi$ is relatively compact and since the embedding $H^{1}\subset H^{1/2}$ is compact, we conclude that $\{\ell(t),\ell\in \mathcal{M}(x_{+},x_{-}),t\in \R\}$ is relatively compact.

\hfill $\Box$

We consider now the set $\tilde{\mathcal{M}}(x_{+},x_{-})$ consisting of elements $\ell \in \mathcal{M}(x_{+},x_{-})$ that we parametrize by arc length. Notice that we have a map $\Pi:\mathcal{M}(x_{+},x_{-})\to\tilde{\mathcal{M}}(x_{+},x_{-})$.
That is our new set of curves consists of maps $\tilde{\ell}:(0,1)\to \mathcal{F}^{1,1/2}$ such that $\tilde{\ell}$ is the re-parametrization by arc length of an element $\ell\in \mathcal{M}(x_{+},x_{-})$. Observe that we can also include the points $\{0,1\}$ by adding their respective images $x_{+}$ and $x_{-}$. So we can say $\tilde{\mathcal{M}}(x_{+},x_{-})\subset C([0,1],\mathcal{F}^{1,2})$.\\
Notice now that $\tilde{\ell}([0,1])=\ell(\R)\cup\{x_{+},x_{-}\}$, therefore
$$\{\tilde{\ell}(t),\tilde{\ell} \in \tilde{\mathcal{M}}(x_{+},x_{-}),t\in \R\}=\{\ell(t),\ell\in \mathcal{M}(x_{+},x_{-}),t\in \R\}\cup\{x_{+},x_{-}\}. $$
Hence, the set $\{\tilde{\ell}(t),\tilde{\ell} \in \tilde{\mathcal{M}}(x_{+},x_{-}),t\in \R \}$ is relatively compact. Also, since the maps in this set are parametrized by arc length, and the length of curves in $\mathcal{M}(x_{+},x_{-})$ is uniformly bounded, then 
the maps $\tilde{\mathcal{M}}(x_{+},x_{-})\subset C([0,1],\mathcal{F}^{1,2})$ are uniformly Lipschitz. It follows from the Arzela-Ascoli theorem that $\tilde{\mathcal{M}}(x_{+},x_{-})$ is a relatively compact subset of $\subset C([0,1],\mathcal{F}^{1,2})$. In order to conclude, we only need to notice that the map $\Pi:\mathcal{M}(x_{+},x_{-})\to\tilde{\mathcal{M}}(x_{+},x_{-})$, induced a homeomorphism $\tilde{\Pi}:\mathcal{M}(x_{+},x_{-})/\R \to \tilde{\mathcal{M}}(x_{+},x_{-}).$

\hfill$\Box$

In this part we consider two critical points of $\Lag_{H}$, $\mathsf{x}_{-},\mathsf{x}_{+}\in crit(\Lag_{H})$ such that $\mu_{H}(\mathsf{x}_{-})-\mu_{H}(\mathsf{x}_{+}) \leq 2$. For $0<\rho_{0}<\delta$ ($\delta>0$ is as in Lemma 5.2), we define the map 
$$\mathcal{F}_{\mathsf{x_{-},x_{+}}}: W^{1,2}_{\mathsf{x_{-},x_{+}}}(\R,\mathcal{F}^{1,1/2}(S^{1},N))\times \mathbb{K}^{2}_{\theta,\rho_{0}}\to L^{2}(\R,T\mathcal{F}^{1,1/2}(S^{1},N))$$
by
$$\mathcal{F}_{\mathsf{x_{-},x_{+}}}(\ell,K):=\frac{d\ell}{dt}+(1+K(\ell(t)))\nabla_{1,1/2}\Lag_{H}(\ell(t)).$$

We set $\mathcal{M}(\mathsf{{x}_{-},x_{+}})=\mathcal{F}_{\mathsf{x_{-},x_{+}}}^{-1}(0)$.

\begin{lemma}
For $(\ell,K)\in\mathcal{M}(\mathsf{{x}_{-},x_{+}})$ the map $$d_{\ell}\mathcal{F}_{\mathsf{x_{-},x_{+}}}:T_{\ell}W^{1,2}_{\mathsf{{x}_{-},x_{+}}}(\R,\mathcal{F}^{1,1/2}(S^{1},N))\to L^{2}(\R,\ell^{*}T\mathcal{F}^{1,1/2}(S^{1},N))$$
is a Fredholm operator with ${\rm index}(d_{\ell}\mathcal{F}_{\mathsf{x_{-},x_{+}}})=\mu_{H}(\mathsf{x}_{-})-\mu_{H}(\mathsf{x}_{+})$.
\end{lemma}
\textit{ Proof.}
First of all, since $\theta =0$ in the neighborhood of $\mathsf{x}_{+}$ and $\mathsf{x}_{-}$, we can assume as in the proof of Proposition 3.4, that $\ell$ is compactly supported and that $\ell(t)=\mathsf{x}_{-}$ when $t\leq -T$ and $\ell(t)=\mathsf{x}_{+}$ when $t\geq T$. Hence, using the same notations, we just need to prove the index formula for the operator
\begin{equation}
\mathsf{P}^{-1}\circ d_{\ell}\mathcal{F}_{\mathsf{x_{-},x_{+}}}(\ell,K)\circ\mathsf{P}:W^{1,2}(\R,T_{\mathsf{x}_{-}}T\mathcal{F}^{1,1/2}(S^{1},N))\to L^{2}(\R,T_{\mathsf{x}_{-}}T\mathcal{F}^{1,1/2}(S^{1},N)).
\end{equation}
Again we write
$$\mathsf{P}^{-1}\circ d\mathcal{F}_{\mathsf{x_{-},x_{+}}}(\ell)\circ\mathsf{P}=\nabla_{t}+\mathsf{P}^{-1}\circ d\nabla_{1,1/2,K}\Lag_{H}(\ell)\circ\mathsf{P}.$$
We set $A(t):=\mathsf{P}_{t}^{-1}\circ d\nabla_{1,1/2,K}\Lag_{H}(\ell(t))\circ\mathsf{P}_{t}$. 
We observe that
$$A(-\infty)=A(-T)=\mathsf{P}^{-1}_{-T}\circ d\nabla_{1,1/2}\Lag_{H}(\ell(-T))\circ\mathsf{P}_{-T}=d\nabla_{1,1/2}\Lag_{H}(\mathsf{x}_{-})$$
and
$$A(+\infty)=A(T)=\mathsf{P}_{T}^{-1}\circ d\nabla_{1,1/2}\Lag_{H}(\ell(T))\circ\mathsf{P}_{T}=\mathsf{P}^{-1}_{T}\circ d\nabla_{1,1/2}\Lag_{H}(\mathsf{x}_{+})\circ\mathsf{P}_{T}$$
are invertible hyperbolic operators since we have assumed that $\mathsf{x}_{-},\mathsf{x}_{+}\in\text{crit}(\Lag_{H})$ are non-degenerate. In this sense, operator family $\{A(t)\}_{t\in\R}$ is asymptotically hyperbolic. 
In this case, we claim that
$$A(t)=d\nabla_{1,1/2}\Lag_{H}(\mathsf{x}_{-})+k(t),$$
where the $(2,2)$-component $k_{22}(t)$ of 
\begin{align*}
k(t)=\begin{pmatrix}k_{11}(t)&k_{12}(t)\\k_{21}(t)&k_{22}(t)\end{pmatrix}:&T_{\mathsf{x}_{-}}\mathcal{F}^{1,1/2}(S^{1},N)=H^{1}(S^{1},\phi_{-}^{\ast}TN)\times H^{1/2}(S^{1},\Spin(S^{1})\otimes\phi_{-}^{\ast}TN)\\
&\to T_{\mathsf{x}_{-}}\mathcal{F}^{1,1/2}(S^{1},N)=H^{1}(S^{1},\phi_{-}^{\ast}TN)\times H^{1/2}(S^{1},\Spin(S^{1})\otimes\phi_{-}^{\ast}TN)
\end{align*}
is compact for all $t\in\R$.
Indeed, we first write 
$$d_{\ell}\nabla_{G^{K}}\Lag_{H}(\ell)[v]=(1+K(\ell))d_{\ell}\nabla_{1,1/2}\Lag_{H}(\ell)[v]+d_{\ell}K[v]\nabla_{1,1/2}\Lag_{H}(\ell)$$
Thus, as in the proof of Lemma 3.4, we can write
\begin{equation}
d_{l}\nabla_{G^{K}}\Lag_{H}(\ell)=(1+K(\ell))\begin{pmatrix}(-\nabla_{s}^{2}+1)^{-1}(-\nabla_{s}^{2})&O\\
O&(1+|\D_{\phi_{-}}|)^{-1}\D_{\phi_{-}}\end{pmatrix}+\mathcal{K}(t)+d_{\ell}K[\cdot]\nabla_{1,1/2}\Lag_{H}(\ell),
\end{equation}
where $\mathcal{K}(t):T_{\ell(t)}\mathcal{F}^{1,1/2}\to T_{\ell(t)}\mathcal{F}^{1,1/2}$ is compact.

By the construction of the element of $\mathbb{K}_{\theta}^{2}$, $K$, 
$$K(\ell(t))\begin{pmatrix}(-\nabla_{s}^{2}+1)^{-1}(-\nabla_{s}^{2})&O\\
O&(1+|\D_{\phi_{-}}|)^{-1}\D_{\phi_{-}}\end{pmatrix}$$
is compact for all $t$. On the other hand, the compactness of the embedding $H^{1}(S^{1})\subset H^{1/2}(S^{1})$ implies that the $(2,2)$-component of $d_{\ell}K[\cdot]\nabla_{1,1/2}\Lag_{H}(\ell)$ is compact. Thus, we have the assertion as stated above. Under the condition, we claim that the operator (5.22) is Fredholm with index $\mu_{H}(\mathsf{x}_{-})-\mu_{H}(\mathsf{x}_{+})$. Note that the present case is different from the one of Lemma 3.4, since $k(t)$ may not be compact (we only require the compactness of the $(2,2)$-component of $k(t)$). To prove the Fredholm property, we argue as follows. We write $A(t)=\begin{pmatrix}L_{1}&0\\0&L_{2}\end{pmatrix}+\begin{pmatrix}k_{11}(t)&k_{12}(t)\\ k_{21}(t)&k_{22}(t)\end{pmatrix}$, where
$L_{1}=(-\nabla_{s}^{2}+1)^{-1}(-\nabla_{s}^{2})$ and $L_{2}=(1+|\D_{\phi_{-}}|)^{-1}\D_{\phi_{-}}$. We first observe that $\nabla_{t}+A(t)$ has closed range with finite dimensional cokernel. To see this, we have
$$(\nabla_{t}+A(t))\begin{pmatrix}\xi\\0\end{pmatrix}=\begin{pmatrix}\nabla_{t}\xi+L_{1}\xi+k_{11}(t)\xi\\ k_{21}(t)\xi\end{pmatrix}.$$
Note that $L_{1}+k_{11}(\pm\infty)$ are essentially positive hyperbolic operator and in such a case, $\nabla_{t}+L_{1}+k_{11}:W^{1,2}(\R,H^{1}(S^{1},\phi_{-}^{\ast}TN)\to L^{2}(\R,\phi_{-}^{\ast}TN)$ is Fredholm as was proved, for example in~\cite[Theorem 5.1]{AM2}. In particular, the codimension of the range of $\nabla_{t}+L_{1}+k_{11}$ is finite. On the other hand, 
$$(\nabla_{t}+A(t))\begin{pmatrix}0\\X\end{pmatrix}=\begin{pmatrix}k_{12}(t)X\\L_{2}X+k_{22}(t)X\end{pmatrix}.$$
As before, $L_{2}+k_{22}(\pm\infty)$ is hyperbolic. In addition, $k_{22}(t)$ is compact for all $t$. Thus,~\cite[Theorem B]{AM2} implies that $\nabla_{t}+L_{2}+k_{22}:W^{1,2}(\R,H^{1/2}(S^{1},\Spin(S^{1})\otimes\phi_{-}^{\ast}TN))\to L^{2}(\R,H^{1/2}(S^{1},\Spin(S^{1})\otimes\phi_{-}^{\ast}TN))$ is Fredholm. In particular, the codimension of the range of $\nabla_{t}+L_{2}+k_{22}$ is finite. Combining these, the range of $\nabla_{t}+A(t)$ is closed with finite dimensional cokernel. Because the adjoint $-\nabla_{t} +A(t)^{\ast}$ is a similar form, $-\nabla_{t} +A(t)^{\ast}$ has closed range with finite dimensional cokernel. But this means that $\nabla_{t}+A(t)$ has a finite dimensional kernel. This proves that $\nabla_{t}+A(t)$ is Fredholm.\\
\noindent
To prove the index formula, for $0\le\epsilon\le1$ we consider the operator $\nabla_{t}+A_{\epsilon}(t)$, where 
$$A_{\epsilon}(t)=\begin{pmatrix}L_{1}&0\\0&L_{2}\end{pmatrix}+\begin{pmatrix}k_{11}(t)&\epsilon k_{12}(t)\\ \epsilon k_{21}(t)&k_{22}(t)\end{pmatrix}.$$
By the same reasoning as given above, $\nabla_{t}+A_{\epsilon}:W^{1,2}(\R,T_{\mathsf{x}_{-}}T\mathcal{F}^{1,1/2}(S^{1},N))\to L^{2}(\R,T_{\mathsf{x}_{-}}T\mathcal{F}^{1,1/2}(S^{1},N))$ is Fredholm for $0\le\epsilon\le1$. Since the Fredholm index is invariant under continuous deformations, we have
\begin{align*}
\text{index}(\nabla_{t}+A)&=\text{index}(\nabla_{t}+A_{0})\\
&=\text{index}(\nabla_{t}+L_{1}+k_{11})+\text{index}(\nabla_{t}+L_{2}+k_{22}).
\end{align*}
Since $L_{1}+k_{11}(\pm\infty)$ are essentially positive, the first index in the above formula is given by the spectral flow 
$$\text{index}(\nabla_{t}+L_{1}+k_{11})=\mathsf{sf}\{L_{1}+k_{11}(t)\}_{-\infty\le t\le+\infty}.$$
The second index is also given by the spectral flow since $k_{22}(t)$ is compact for all $t$, see~\cite{AM2}:
$$\text{index}(\nabla_{t}+L_{2}+k_{22})=\mathsf{sf}\{L_{2}+k_{22}(t)\}_{-\infty\le t\le+\infty}.$$
Combining these, we have
\begin{align*}\text{index}(\nabla_{t}+A)&=\mathsf{sf}\{L_{1}+k_{11}(t)\}_{-\infty\le t\le+\infty}+\mathsf{sf}\{L_{2}+k_{22}(t)\}_{-\infty\le t\le+\infty}\\
&=\mathsf{sf}\{A_{0}(t)\}_{-\infty\le t\le+\infty}.
\end{align*}
Since the spectral flow is a homotopy invariant, we have 
$$\mathsf{sf}\{A_{0}(t)\}_{-\infty\le t\le+\infty}=\mathsf{sf}\{A(t)\}_{-\infty\le t\le+\infty}.$$
Since $\text{sf}\{A(t)\}_{-\infty\le t\le +\infty}=\mu_{H}(\mathsf{x}_{-})-\mu_{H}(\mathsf{x}_{+})$ as in (3.110), we finally have the index formula as stated in the lemma.
\hfill$\Box$

\begin{lemma}
Let $K\in \mathbb{K}^{2}_{\theta,\rho_{0}}$. The map $d_{\ell}\mathcal{F}_{\mathsf{x_{-},x_{+}}}:T_{\ell}W^{1,2}_{x_{-},x_{+}}(\R,\mathcal{F}^{1,1/2}(S^{1},N))\to L^{2}(\R,\ell^{*}T\mathcal{F}^{1,1/2}(S^{1},N))$ is onto for any $\ell$ such that $(\ell,K)\in \mathcal{M}_{x_{-},x_{+}}$ if and only if $W^{u}_{-\nabla_{G^{K}}\Lag_{H}}(\mathsf{x}_{-})$ and $W^{u}_{-\nabla_{G^{K}}\Lag_{H}}(\mathsf{x}_{+})$ intersect transversally. 
\end{lemma}

\textit{ Proof.}
Using the same notations as in Proposition 3.4, we consider the operator
\begin{equation}
\tilde{\mathcal{F}}=\mathsf{P}^{-1}\circ d_{\ell}\mathcal{F}_{\mathsf{x_{-},x_{+}}}(\ell,K)\circ\mathsf{P}:W^{1,2}(\R,T_{\mathsf{x}_{-}}T\mathcal{F}^{1,1/2}(S^{1},N))\to L^{2}(\R,T_{\mathsf{x}_{-}}T\mathcal{F}^{1,1/2}(S^{1},N)).
\end{equation}
Indeed, we have
$$\tilde{\mathcal{F}}=\mathsf{P}^{-1}\circ d_{\ell}\mathcal{F}_{\mathsf{x_{-},x_{+}}}(\ell)\circ\mathsf{P}=\nabla_{t}+\mathsf{P}^{-1}\circ d\nabla_{G^{K}}\Lag_{H}(\ell)\circ\mathsf{P}.$$
Therefore $\tilde{\mathcal{F}}$ is onto if and only if $d_{\ell}\mathcal{F}_{\mathsf{x_{-},x_{+}}}(\ell)$ is onto. Using the results from \cite{AM3}, we have that $W^{s}(\mathsf{P}^{-1}\circ d\nabla_{G^{K}}\Lag_{H}(\ell)\circ\mathsf{P})$ and $W^{u}(\mathsf{P}^{-1}\circ d\nabla_{G^{K}}\Lag_{H}(\ell)\circ\mathsf{P})$ intersect transversally. But this holds if and only if 
$W^{u}_{-\nabla_{G^{K}}\Lag_{H}}(x_{-})$ and $W^{s}_{-\nabla_{G^{K}}\Lag_{H}}(x_{+})$ intersect transversally.

\hfill$\Box$

\begin{lemma}
Let $(\ell,K)\in \mathcal{M}(\mathsf{x_{-},x_{+}})$ and consider the operator $$d_{\ell}\mathcal{F}_{\mathsf{x_{-},x_{+}}}^{+}(\ell,K):T_{\ell}W^{1,2}_{\mathsf{x_{-},x_{+}}}(\R^{+},\mathcal{F}^{1,1/2}(S^{1},N))\to L^{2}(\R^{+},\ell^{*}T\mathcal{F}^{1,1/2}(S^{1},N))$$ defined by restricting $d_{\ell}\mathcal{F}_{\mathsf{x_{-},x_{+}}}(\ell,K)$ to $T_{\ell}W^{1,2}_{\mathsf{x_{-},x_{+}}}(\R^{+},\mathcal{F}^{1,1/2}(S^{1},N))$. Then $d_{\ell}\mathcal{F}_{\mathsf{x_{-},x_{+}}}^{+}(\ell,K)$ is a left inverse.
\end{lemma}
\textit{ Proof.}
Again by using the same trick of parallel transport, it is enough to show the property for the operator
\begin{equation}
\tilde{\mathcal{F}}^{+}=\mathsf{P}^{-1}\circ d_{\ell}\mathcal{F}_{\mathsf{x_{-},x_{+}}}^{+}(\ell,K)\circ\mathsf{P}:W^{1,2}(\R^{+},T_{\mathsf{x}_{-}}T\mathcal{F}^{1,1/2}(S^{1},N))\to L^{2}(\R^{+},T_{\mathsf{x}_{-}}T\mathcal{F}^{1,1/2}(S^{1},N)).
\end{equation}
Then since it is asymptotically hyperbolic, one have again from \cite{AM2}. That it is a left inverse.

\hfill$\Box$

\begin{lemma}
Let $(\ell,K)\in \mathcal{M}(\mathsf{x_{+},x_{-}})$ and assume that there exists $a$, $b\in\R$ with $a<b$ such that $\theta(\ell(t))\not = 0$ for all $t\in[a,b]$. Then for any $w\in C^{2}(\R,\ell^{*}T\mathcal{F}^{1,1})$ with $supp$ $w\subset [a,b]$ and for $\varepsilon>0$, there exists $k\in \mathbb{K}_{\theta}^{2}$ such that $$\|d_{K}\mathcal{F}_{\mathsf{x_{-},x_{+}}}(\ell,K)[k]-w\|_{\ell^{*}T\mathcal{F}^{1,1}}<\varepsilon.$$
\end{lemma}
\textit{ Proof.}
First we notice that $$d_{K}\mathcal{F}_{\mathsf{x_{-},x_{+}}}(\ell,K)[k] = k(\ell)\nabla_{1,1/2} \Lag_{H}(\ell)$$
and clearly $\ell :\R \to \mathcal{F}^{1,1}$ is an embedding of class $C^{3}$ since it is a solution to the negative gradient flow equation with respect to the gradient $\nabla_{G^{K}} \Lag_{H}$.  We consider the $C^{2}$ curve defined as follow:
$$\tilde{k}(t)=\frac{\langle \cdot, \nabla_{1,1/2}\Lag_{H}(\ell(t))\rangle}{\|\nabla_{1,1/2}\Lag_{H}(\ell(t))\|^{2}}w(t)+\frac{\langle \cdot, w(t)\rangle}{\|\nabla_{1,1/2}\Lag_{H}(\ell(t))\|^{2}}\nabla_{1,1/2}\Lag_{H}(\ell(t))$$
$$-\frac{\langle \nabla_{1,1/2}\Lag_{H}(\ell(t)), w(t)\rangle\langle \cdot ,\nabla_{1,1/2}\Lag_{H}(\ell(t)) \rangle}{\|\nabla_{1,1/2}\Lag_{H}(\ell(t))\|^{4}}\nabla_{1,1/2}\Lag_{H}(\ell(t)).$$
Now notice that since $w(t) \in T_{\ell(t)}\mathcal{F}^{1,1}$ and from Lemma 5.4, $\nabla_{1,1/2}\Lag_{H}(\ell(t))$ is also in $T_{\ell(t)}\mathcal{F}^{1,1}$, we have $\tilde{k}(t)\in Sym(T_{\ell(t)}\mathcal{F}^{1,1/2},T_{\ell(t)}\mathcal{F}^{1,1})$ and if $H$ is regular enough, we have that $\tilde{k}$ is $C^{2}$ in the $t$ variable and $\dot{\tilde{k}}$, $\ddot{\tilde{k}}$ take values in $H^{1}\times H^{1}$. It is also important to notice that $\tilde{k}$ might not have the form of an element in $\mathbb{K}_{\theta}^{2}$ but since the interval $[a,b]$ is compact, we can take a partition $t_{0}=a<t_{1}<\cdots<t_{n}=b$ and replace $w(t)$ by $w(t_{i})$ on the interval $[t_{i},t_{i+1}]$ and using an approximation of identity adapted to the partition, we have indeed an element $\tilde{k}_{\varepsilon}$ having the same form as an operator in $\mathbb{K}_{\theta}^{2}$ such that $\|\tilde{k}(t)-\tilde{k}_{\varepsilon}(t)\|_{op}<\varepsilon$.\\

Now since $\ell:\R\to \mathcal{F}^{1,1/2}(S^{1},N)$ is an embedding, for a given $\delta>0$, there exists a neighborhood $U\in \mathcal{F}^{1,1/2}(S^{1},N)$ such that $\ell((a-\delta,b+\delta))\subset U$ and a submersion $\tau : U\to (a-\delta,b+\delta)$ such that $\tau(\ell(t))=t$ for $t\in (a-\delta,b+\delta)$. But we already assumed that $\theta>0$ on $\ell([a,b])$, thus taking $\delta>0$ and $U$ smaller if necessary, we can assume that $\inf_{U}\theta>0$. Therefore, we can consider a cut-off function $\rho\in C^{2}_{b}(\mathcal{F}^{1,1/2}(S^{1},N),\R)\cap C_{b}(\mathcal{F}^{s,\frac{1}{2}},\R)$ such that $\rho=1$ on $\ell([a,b])$ and $supp$ $\rho\subset U$. We set then
$$k(p)=\rho(p)\tilde{k}(\tau(p)),$$
for $p\in \mathcal{F}^{1,1/2}(S^{1},N)$.

We set $\ell(\tau(p))=\tilde{\tau}(p)$ for $p\in U$ and write $\tilde{\tau}(p)=(\tilde{\phi}_{p},\tilde{\psi}_{p})$. We also write $p=(\phi_{p},\psi_{p})$ for $p\in U$. For $U$ small enough, $\|\phi_{p}-\tilde{\phi}_{p}\|_{H^{1}}$ is small and so is $\|\phi_{p}-\tilde{\phi}_{p}\|_{L^{\infty}}$ by the Sobolev embedding $H^{1}(S^{1})\subset L^{\infty}(S^{1})$. Thus the parallel translation $\mathsf{P}_{\tilde{\phi}_{p},\phi_{p}}(s):T_{\tilde{\phi}_{p}(s)}N\to T_{\phi_{p}(s)}N$ is defined for all $s\in S^{1}$. By~\cite[Lemma 7.4]{I}, it induces a bounded linear operator $\mathsf{P}_{\tilde{\phi}_{p},\phi_{p}}:H^{1}(S^{1},\tilde{\phi}_{p}^{\ast}TN)\to H^{1}(S^{1},\phi_{p}^{\ast}TN)$ which we denote by the same symbol. We finally define
$$\hat{k}(p)=\mathsf{P}_{\tilde{\phi}_{p},\phi_{p}}(k(p)).$$
It is easy to check that $\hat{k}$ belongs to $\mathbb{K}_{\theta}^{2}$ and this finishes the proof.
\hfill$\Box$

\begin{lemma}
For any $(\ell,K)\in \mathcal{M}(\mathsf{x_{-},x_{+}})$, $d\mathcal{F}_{\mathsf{x_{-},x_{+}}}(\ell,K):T_{\ell}W^{1,2}_{\mathsf{x_{-},x_{+}}}(\R,\mathcal{F}^{1,1/2}(S^{1},N))\times \mathbb{K}_{\theta}^{2}\to L^{2}(\R,\ell^{*}T\mathcal{F}^{1,1/2}(S^{1},N))$ is a left inverse.
\end{lemma}
\textit{ Proof.}
We fix $(\ell,K)\in \mathcal{M}(\mathsf{x_{-},x_{+}})$. We already proved in Lemma 5.5 that $d_{\ell}\mathcal{F}_{\mathsf{x_{-},x_{+}}}(\ell,K)$ is Fredholm. So it is enough to show that $d\mathcal{F}_{\mathsf{x_{-},x_{+}}}(\ell,K):T_{\ell}W^{1,2}_{\mathsf{x_{-},x_{+}}}(\R,\mathcal{F}^{1,1/2}(S^{1},N))\times \mathbb{K}_{\theta}^{2}\to L^{2}(\R,\ell^{*}T\mathcal{F}^{1,1/2}(S^{1},N))$ is onto. But since we have $\text{Ran}( d_{\ell}\mathcal{F}_{\mathsf{x_{-},x_{+}}}(\ell,K))\subset \text{Ran}(d\mathcal{F}_{\mathsf{x_{-},x_{+}}}(\ell,K))$, we have that $Codim\text{ } d\mathcal{F}_{\mathsf{x_{-},x_{+}}}(\ell,K)<\infty$ and hence to proof our lemma we just need to show that $\text{Ran} (d\mathcal{F}_{\mathsf{x_{-},x_{+}}}(\ell,K))$ is dense in $L^{2}(\R,\ell^{*}T\mathcal{F}^{1,1/2}(S^{1},N))$. Again we distinguish two cases.\\
\begin{itemize}
\item[i)] \textit{The case $\theta(\ell(t))=0$ for all $t$}.\\
In that case, by construction of the space $\mathbb{K}_{\theta}^{2}$ we have that $K(\ell(t))=0$ for all $t$ and hence $\nabla_{G^{K}}\Lag_{H}=\nabla_{1,1/2}\Lag_{H}$ and thus the stable and unstable manifolds of $-\nabla_{G^{K}}\Lag_{H}$ coincides with the ones of $-\nabla_{1,1/2}\Lag_{H}$ and by assumption they intersect transversally since $\theta=0$ therefore from Lemma 5.6, we have that $d_{\ell}\mathcal{F}_{\mathsf{x_{-},x_{+}}}(\ell,K)$ is onto.

\item[ii)] \textit{The case $\theta(\ell(t))\not=0$ for some $t$}.\\
Let $w\in  L^{2}(\R,\ell^{*}T\mathcal{F}^{1,1/2})$ and $\varepsilon>0$. We consider an interval $[a,b]$ such that $\theta(\ell(t))\not=0$. Now by using Lemma 5.7 we have the existence of $v\in T_{\ell}W^{1,2}_{\mathsf{x_{-},x_{+}}}(\R,\mathcal{F}^{1,1/2}(S^{1},N)) $ such that 
$$d_{\ell}\mathcal{F}_{\mathsf{x_{-},x_{+}}}(\ell,K)[v]=w$$
for $t\in (-\infty,a]\cup [b,+\infty)$. Then the support of $\tilde{w}=w-d_{\ell}\mathcal{F}_{\mathsf{x_{+},x_{-}}}(\ell,K)[v]\in L^{2}(\R,\ell^{*}T\mathcal{F}^{1,1/2})$ is contained in $[a,b]$. Now we need to approximate $\tilde{w}$ by a function in $C^{2}(\R,\ell^{*}T\mathcal{F}^{1,1})$. This can be done first by transporting the space to the point $\ell(a)$. Indeed, using the notations of the proof of Proposition 3.4, we consider the operator $P_{t}$ defining the parallel transport along the path $\phi(\tau,s)$ for $a\leq \tau \leq t$, here $\phi$ is so that $\ell(t)=(\phi(t),\psi(t))$. Then we have 
$$P_{t}(s):T_{\phi(a,s)}N\to T_{\phi(t,s)}N.$$
Thus $P_{t}$ induces a bounded linear operator $$P_{t}:H^{1}(S^{1},\phi(a)^{*}TN)\to H^{1}(S^{1},\phi(t)^{*}TN)$$
Similarly 
$$P_{t}=1\otimes P_{t}:H^{1/2}(S^{1},\mathbb{S}(S^{1})\otimes \phi(a)^{*}TN)\to H^{1/2}(S^{1},\mathbb{S}(S^{1})\otimes \phi(t)^{*}TN)$$
Thus by setting $\underline{w}=P_{t}^{-1}\tilde{w}$, we see that $\underline{w}\in L^{2}([a,b],T_{\ell(a)}\mathcal{F}^{1,1/2})$ and thus since it is valued in a fixed vector space it can be easily approximated by a smooth function $\underline{w}_{\varepsilon}$ so that
$$\|\underline{w}-\underline{w}_{\varepsilon}\|_{L^{2}([a,b],T_{\ell(a)}\mathcal{F}^{1,1/2})}<\varepsilon$$
Thus by taking $w_{\varepsilon}=P_{t}\underline{w}_{\varepsilon}$ we have that $w_{\varepsilon}\in C^{2}(\R,\ell^{*}T\mathcal{F}^{1,1})$ and 
$$\|\tilde{w}-w_{\varepsilon}\|_{ L^{2}(\R,\ell^{*}T\mathcal{F}^{1,1/2})}<C\varepsilon$$
for a fixed constant $C$ depending on the parallel transport norm $P_{t}$. Now using Lemma 5.8, we have the existence of $k$ such that 
$$\|d_{K}\mathcal{F}_{\mathsf{x_{+},x_{-}}}(\ell,K)[k]-w_{\varepsilon}\|_{L^{2}(\R,\ell^{*}T\mathcal{F}^{1,\frac{1}{2}})}<\varepsilon$$
Therefore 
$$\|d\mathcal{F}_{\mathsf{x_{+},x_{-}}}(\ell,K)[v,k]-w\|_{ L^{2}(\R,\ell^{*}T\mathcal{F}^{1,1/2})}\leq C\varepsilon.$$
Since $\varepsilon>0$ is arbitrary, we have that $\text{Ran} (d\mathcal{F}_{\mathsf{x_{+},x_{-}}}(\ell,K))\subset L^{2}(\R,\ell^{*}T\mathcal{F}^{1,1/2})$ is dense and this completes the proof.

\end{itemize}
 \hfill$\Box$

\textit{ Proof of Proposition 5.1}
If $\mathsf{x}_{+}=\mathsf{x}_{-}$ the transversality is satisfied since $\ell(t)=\mathsf{x}_{+}$ is a stationary solution and $d\mathcal{F}_{\mathsf{x_{+},x_{+}}}=\frac{d}{dt} + d\nabla_{1,1/2}\Lag_{H}(\mathsf{x_{+}})$ is invertible. Next, we assume that $\mathsf{x}_{+}\not=\mathsf{x}_{-}$. Using the result of Lemma 5.9, we have that $\mathcal{M}(\mathsf{x_{+},x_{-}})$ is a $C^{2}$-submanifold of $W^{1,2}_{\mathsf{x_{+},x_{-}}}(\R,\mathcal{F}^{1,1/2})\times \mathbb{K}^{2}_{\theta,\rho_{0}}$. If we let $\pi_{\mathbb{K}^{2}}:\mathcal{M}(\mathsf{x_{+},x_{-}})\to \mathbb{K}^{2}_{\theta,\rho_{0}}$ the projection onto the second component. It is a Fredholm map of class $C^{2}$ and its index at $(\ell,K)\in \mathcal{M}(\mathsf{x_{+},x_{-}})$ is 
$$index\text{ }\pi_{\mathbb{K}^{2}}(\ell,K)=index\text{ }d_{\ell}\mathcal{F}_{\mathsf{x_{+},x_{-}}}(\ell,K)=\mu_{H}(\mathsf{x}_{+})-\mu_{H}(\mathsf{x}_{-}).$$
Since $\R$ acts freely on $\mathcal{M}(\mathsf{x_{+},x_{-}})$, the quotient $\hat{\mathcal{M}}(\mathsf{x_{+},x_{-}})=\mathcal{M}(\mathsf{x_{+},x_{-}})/\R$ is a $C^{2}$-manifold of dimension $$dim \tilde{\mathcal{M}}(\mathsf{x_{+},x_{-}})=\mu_{H}(\mathsf{x}_{+})-\mu_{H}(\mathsf{x}_{-})-1.$$
Therefore, for $\mathsf{x}_{+}$ and $\mathsf{x}_{-}\in crit(\Lag_{H})$ with $\mu_{H}(\mathsf{x}_{+})-\mu_{H}(\mathsf{x}_{-})\leq 2$ we have that $dim \hat{\mathcal{M}}(\mathsf{x_{+},x_{-}})\leq 1$.
Since $W^{1,2}(\R,\mathcal{F}^{1,1/2})\times \mathbb{K}^{2}_{\theta,\rho_{0}}$ is second countable, we have that $\hat{\mathcal{M}}(\mathsf{x_{+},x_{-}})$ is also second countable and by restricting the map $\pi_{\mathbb{K}^{2}}$ to $\tilde{\pi}_{\mathbb{K}^{2}}: \hat{\mathcal{M}}(\mathsf{x_{+},x_{-}})\to \mathbb{K}^{2}_{\theta,\rho_{0}}$ which is $C^{2}$, one applies the Sard-Smale theorem to conclude that the regular values of $\tilde{\pi}_{\mathbb{K}^{2}}$ that we denote by $\mathbb{K}_{reg}(\mathsf{x_{+},x_{-}})$ is residual in $\mathbb{K}^{2}_{\theta,\rho_{0}}$. We set then
$$\mathbb{K}_{reg}=\bigcap_{\mathsf{x_{+},x_{-}}\in crit(\Lag_{H}),\mathsf{x}_{+}\not=\mathsf{x}_{-},\mu_{H}(\mathsf{x}_{+})-\mu_{H}(\mathsf{x}_{-})\leq 2}\mathbb{K}_{reg}(\mathsf{x}_{+},\mathsf{x}_{-}).$$
Since the set $\{(\mathsf{x_{+},x_{-}})\in crit(\Lag_{H})\times crit(\Lag_{H});\mu_{H}(\mathsf{x_{+}})-\mu_{H}(\mathsf{x_{-}})\leq 2\}$ is countable, we have that $\mathbb{K}_{reg}\subset \mathbb{K}^{2}_{\theta,\rho_{0}}$ is residual. Hence for any $K\in \mathbb{K}_{reg}$ the operator 
$$d_{\ell}\mathcal{F}_{\mathsf{x_{+},x_{-}}}:T_{\ell}W^{1,2}_{\mathsf{x_{+},x_{-}}}(\R,\mathcal{F}^{1,1/2})\to L^{2}(\R,\ell^{*}T\mathcal{F}^{1,1/2})$$
is onto.
\hfill$\Box$
\section{Invariance of the Morse-Floer homology, I: independence of the metric on $\mathcal{F}^{1,1/2}(S^{1},N)$}

\begin{proposition}
Let $K_{0}$, $K_{1} \in \mathbb{K}_{reg}$ then there exists a chain complex isomorphism
$$\Phi : \{C_{*}(\Lag_{H},K_{0}),\partial_{*}(K_{0})\} \to \{C_{*}(\Lag_{H},K_{1}),\partial_{*}(K_{1})\}.$$
In particular the homology $HF_{*}(\Lag_{H},K_{0})$ is isomorphic to $HF_{*}(\Lag_{H},K_{1})$

\end{proposition}
The proof of this proposition will be made through several steps. First we will define an alternative functional by adding an extra real variable. Consider the function $f:\R\to \R$ defined by $f(s)=2s^{3}-3s^{2}+1$ and we set the space $\tilde{\mathcal{F}}^{1,1/2}=\mathcal{F}^{1,1/2}\times \R,$ and we extend our functional $\Lag_{H}$ by $\tilde{\Lag_{H}}$ as follows :
$$\tilde{\Lag}_{H}(x,s)=\Lag_{H}(x)+f(s).$$
We will assume that $\Lag_{H}$ is Morse, therefore $\tilde{\Lag}_{H}$ is also Morse and
$$crit(\tilde{\Lag}_{H})=\{(x,0),(x,1);x\in crit(\Lag_{H})\}.$$
Recall that the relative Morse index for a critical point $x$ of $\Lag_{H}$ is as follow

$$\mu_{H}(x)=-\mathsf{sf}\{\mathcal{A}_{x_{t},H}\}_{0\le t\le 1}$$
Where here $\mathcal{A}$ is the Hessian of $\Lag_{H}$ and $x_{t}$ is a path from a fixed point $x_{0}\in \mathcal{F}^{1,1/2}$ to the critical point $x$. In the same way we define the relative index for critical points of $\tilde{\Lag}_{H}$ as
$$\tilde{\mu}_{H}(x,i)=-\mathsf{sf}\{\tilde{\mathcal{A}}_{x_{t},s_{t},H}\}_{0\le t\le 1},$$
for $i=0,1$, where here
$$\tilde{\mathcal{A}}_{x,s,H}=\left(\begin{array}{cc}
\mathcal{A}_{x,H}& 0\\
0&2s-1 
\end{array}
\right)
$$
and $s(t)=(1-t)2+ti$. Therefore from the properties of the spectral flow we have that
$$\mathsf{sf}\{\tilde{\mathcal{A}}_{x_{t},s_{t},H}\}_{0\le t\le 1}=\mathsf{sf}\{\mathcal{A}_{x_{t},H}\}_{0\le t\le 1}+\mathsf{sf}\{(2s(t)-1)\}_{0\le t\le 1}.$$
Hence,
$$\tilde{\mu}_{H}(x,0)=\mu_{H}(x)+1, \text{      }\tilde{\mu}_{H}(x,1)=\mu_{H}(x).$$
This gives us a splitting of the set of critical points of $\tilde{\Lag}_{H}$ of the form
$$crit_{k}(\tilde{\Lag}_{H})=\left(crit_{k-1}(\Lag_{H})\times \{0\}\right)\cup \left(crit_{k}(\Lag_{H})\times \{1\}\right).$$
This naturally descends to the chain level giving us
$$C_{k}(\tilde{\Lag}_{H})=C_{k-1}(\Lag_{H})\oplus C_{k}(\Lag_{H}).$$
Now we consider a smooth function $\rho \in C^{\infty}(\R,[0,1])$ such that $\rho(s)=1$ for $s\leq \frac{1}{3}$ an $\rho(s)=0$ for $s\geq \frac{2}{3}$. We define then the metric $G^{\rho,K_{0},K_{1}}$, for $X, Y \in T_{x}\mathcal{F}^{1,1/2}(S^{1},N)$ an $a,b\in \R$, by
$$G^{\rho,K_{0},K_{1}}((X,a),(Y,b))=\rho(s)G^{K_{0}}(X,Y)+(1-\rho(s))G^{K_{1}}(X,Y)+ab.$$

\begin{lemma}
The functional $\tilde{\Lag}_{H}$ satisfies the Palais-Smale condition on $\tilde{\mathcal{F}}^{1,1/2}$ equipped with the metric $G^{\rho,K_{0},K_{1}}$.
\end{lemma}

{\it Proof.}
Assume that $|\tilde{\Lag}_{H}(x_{n},s_{n})|\leq C$ and $\nabla_{G^{\rho,K_{0},K_{1}}} \tilde{\Lag}_{H}(x_{n},s_{n})\to 0$. In particular, we have that $f'(s_{n})\to 0$, therefore $s_{n}\to i$ for $i=0,1$, we will assume for instance that it converges to $1$. Hence, $f(s_{n})$ is also bounded and for $n$ big enough $s_{n}$ will be close to $1$. Therefore we have that  $| \Lag_{H}(x_{n},s_{n})|$ is bounded and $\nabla_{1,1/2,G^{\rho,K_{0},K_{1}}} \Lag_{H}(x_{n})=\nabla_{1,1/2,G^{K_{1}}} \Lag_{H}(x_{n})\to 0$, and $G^{K_{1}}$ is equivalent to the standard metric $G$ that we fixed in the beginning. But $\Lag_{H}$ satisfies the Palais-Smale condition, therefore up to a subsequence, $(x_{n})$ converges in $\mathcal{F}^{1,1/2}(S^{1},N)$. 
\hfill$\Box$

We consider now the negative gradient flow of $\tilde{\Lag}_{H}$ with respect to the metric $G^{\rho,K_{0},K_{1}}$. Then one finds three kinds of flow lines.\\
Type 1 : These flow lines connect critical points of the form $(x_{1},0)$ and $(x_{2},0)$ for $x_{1},x_{2}\in crit(\Lag_{H})$. These flow lines take the form $(\ell(t),0)$ and $x(t)$ satisfies 
$$\frac{d\ell}{dt}=-\nabla_{1,1/2,K_{0}}\Lag_{H}(\ell),\text{ } \ell(+\infty)=x_{2},\text{ } \ell(-\infty)=x_{1}.$$
Type 2 : Similarly to the type 1 flow lines, they connect critical points of the form $(x_{1},1)$ and $(x_{2},1)$ for $x_{1},x_{2}\in crit(\Lag_{H})$, therefore they take the form $(\ell(t),1)$ and $x(t)$ satisfies 
$$\frac{d\ell}{dt}=-\nabla_{1,1/2,K_{1}}\Lag_{H}(\ell),\text{ } \ell(+\infty)=x_{2},\text{ } \ell(-\infty)=x_{1}.$$
Type 3 : These flow lines of the form $(\ell(t),s(t))$ connect critical points of the form $(x_{1},0)$ and $(x_{2},1)$ where $x_{1},x_{2}\in crit(\Lag_{H})$. First we restrict the metric $G^{\rho,K_{0},K_{1}}$ to $\mathcal{F}^{1,1/2}$ that we denote it by $\tilde{G}^{\rho,K_{0},K_{1}}$. We notice now that 
$$\frac{d}{dt}\Lag_{H}(\ell(t))=-\|\nabla_{\tilde{G}^{\rho,K_{0},K_{1}}}\Lag_{H}\|^{2}_{\tilde{G}^{\rho,K_{0},K_{1}}}\leq 0,$$
therefore $\Lag_{H}(x_{2})\leq \Lag_{H}(x_{1})$ and $\Lag_{H}(x_{2})=\Lag_{H}(x_{1})$ if and only if $x_{2}=x_{1}$ and the flow is stationary. Hence, we can split these flow lines to two categories.\\
The first category consists of flow lines of the form $(x,s(t))$ where $s(t)$ solves
$$s'(t)=-f'(s(t)),\text{ } s(-\infty)=0,\text{ } s(+\infty)=1.$$
The second category consists of flow lines of the form $(\ell(t),s(t))$ connecting $(x_{1},0)$ and $(x_{2},1)$, where $x_{1}\not =x_{2}$. 

\begin{lemma}
Along the flow lines of the third type and first category the unstable manifold $W^{u}_{\nabla_{G^{\rho,K_{0},K_{1}}}\tilde{\Lag}_{H}}(x,0)$ and the stable manifold $W^{s}_{\nabla_{G^{\rho,K_{0},K_{1}}}\tilde{\Lag}_{H}}(x,1)$ intersect transversally.
\end{lemma}
{\it Proof.}
We consider a smooth function $s_{0}:\R\to \R$ such that $s_{0}(t)=1$ if $t\geq 1$ and $s_{0}(t)=0$ if $t\leq -1$ and let $W=s_{0}+W^{1,2}(\R)$. We define then the map $\tilde{\mathcal{F}}:W^{1,2}_{x,x}(\R,\mathcal{F}^{1,1/2})\times W\to L^{2}(\R,T\mathcal{F}^{1,1/2})\times L^{2}(\R)$, by
$$\tilde{\mathcal{F}}(x,s)=\left(\frac{dx}{dt}+\nabla_{\tilde{G}^{\rho,K_{0},K_{1}}}\Lag_{H}(x),\dot{s}+f'(s)\right).$$
Now $$d\tilde{\mathcal{F}}(x,s):T_{x}W^{1,2}(\R,\mathcal{F}^{1,2})\times W^{1,2}(\R)\to L^{2}(\R,x^{*}T\mathcal{F}^{1,1/2})\times L^{2}(\R)$$
and one sees that
$$d\tilde{\mathcal{F}}(x,s)=d\mathcal{F}_{x,x}(x)\oplus (\frac{d}{dt}+f''(s)).$$
Clearly the map $\frac{d}{dt}+f''(s):\times W^{1,2}(\R)\to \times L^{2}(\R)$ is clearly Fredholm of index 1. Now it remains to study the first operator $d\mathcal{F}_{x,x}:T_{x}W^{1,2}(\R,\mathcal{F}^{1,2})\times L^{2}(\R,x^{*}T\mathcal{F}^{1,1/2})$. Notice now that $$d\mathcal{F}_{x,x}=\frac{d}{dt}+d\nabla_{\tilde{G}^{\rho,K_{0},K_{1}}(t)}\Lag_{H}(x).$$
But since $K_{0}$ and $K_{1}$ are in $\mathbb{K}_{reg}$, we have that $K_{0}(x)=K_{1}(x)=0$, hence $\tilde{G}^{\rho,K_{0},K_{1}}=G_{0}$ the standard metric in $\mathcal{F}^{1,1/2}$. Therefore, the operator $d\mathcal{F}_{x,x}=\frac{d}{dt}+d\nabla_{1,1/2}\Lag_{H}(x)$ and since $\Lag_{H}$ is assumed to be Morse the operator $d\nabla_{1,1/2}\Lag_{H}(x)$ is hyperbolic and the surjectivity follows as in Proposition 5.2.

\hfill$\Box$

Now transversality needs to be checked only for the type 3 orbits from the second category. In this case, we do again the same construction as in Proposition 5.2 to get the following result 
\begin{lemma}
There exists a metric $\hat{G}$ on $\tilde{\mathcal{F}}^{1,1/2}$ arbitrarily close to $G^{\rho,K_{0},K_{1}}$ such that :
\begin{itemize}
\item $\tilde{\Lag}_{H}$ satisfies the Palais-Smale condition on $(\tilde{\mathcal{F}}^{1,1/2},\hat{G})$
\item for any arbitrarily small neighborhood $U$ of $crit(\Lag_{H})$, $\hat{G}$ coincides with $G^{\rho, K_{0},K_{1}}$ on $(\mathcal{F}^{1,1/2}\times(-\infty,\frac{1}{3}))\cup (\mathcal{F}^{1,1/2}\times(\frac{2}{3},+\infty)) \cup (U\times \R).$
\item $d\Lag_{H}(x)[-\nabla_{\hat{G}}\tilde{\Lag}_{H}(x,s)]<0$ if $x\not \in crit(\Lag_{H})$ and $s\in \R$
\item The negative gradient flow of $\tilde{\Lag}_{H}$ with respect to $\hat{G}$ satisfies the Morse-Smale property up to order 2.
\end{itemize}
\end{lemma}
The idea of the proof is again to consider a function $\tilde{\theta}:\tilde{\mathcal{F}}^{1,1/2}\to [0,1]$ with the following properties :
\begin{itemize}
\item[i)] $\tilde{\theta}=0$ on the set $(\mathcal{F}^{1,1/2}\times(-\infty,\frac{1}{3}))\cup (\mathcal{F}^{1,1/2}\times(\frac{2}{3},+\infty))\cup (U\times \R)$.
\item[ii)] The zero set of $\tilde{\theta}$ is the closure of an open set.
\item[iii)] If $W^{u}_{-\nabla_{G^{\rho,K_{0},K_{1}}}\tilde{\Lag}_{H}}(x,i)$ and $W^{s}_{-\nabla_{G^{\rho,K_{0},K_{1}}}\tilde{\Lag}_{H}}(y,j)$, for $(x,i),(y,j)\in crit(\tilde{\Lag}_{H})$, intersect non transversally at a point $(p,s)$, then $\tilde{\theta}$ is identically zero on the orbit passing from $(p,s)$.
\end{itemize}
The rest of the proof is exactly similar to Lemma 5.3-5.9,  so we will omit it.

{\it Proof of Proposition 11.1}
We consider now the negative gradient flow of $\tilde{\Lag}_{H}$ with respect to the metric $\hat{G}$. Then again we have the same decompositions of the flow lines. That is, flow lines of type 1,2 and 3 and the later ones decompose to first and second category. Moreover, we have that the intersections of stable and unstable manifolds of critical points are transversal.\\

We consider now the boundary operator 
\begin{align*}
\partial_{k}(\tilde{\Lag}_{H},\hat{G}):C_{k}(\tilde{\Lag}_{H},\mathbb{Z}_{2})&=C_{k-1}(\Lag_{H},\mathbb{Z}_{2})\oplus C_{k}(\Lag_{H},\mathbb{Z}_{2})\\
&\to C_{k-1}(\tilde{\Lag}_{H},\mathbb{Z}_{2})=C_{k-2}(\Lag_{H},\mathbb{Z}_{2})\oplus C_{k-1}(\Lag_{H},\mathbb{Z}_{2})
\end{align*}
Then one has that $\partial_{k}(\tilde{\Lag}_{H},\hat{G})$ takes the form 
$$\partial_{k}(\tilde{\Lag}_{H},\hat{G})=\left(\begin{array}{cc}
\partial_{k-1}(\Lag_{H},G^{K_{0}}) & 0\\
\Phi_{k-1} & \partial_{k}(\Lag_{H},G^{K_{1}})
\end{array}
\right),$$
where $\Phi_{k}:C_{k}(\Lag_{H},\mathbb{Z}_{2})\to C_{k}(\Lag_{H},\mathbb{Z}_{2})$ is defined by counting the flow lines of type 3 $(\text{mod 2})$. Notice that since Plais-Smale holds for $\tilde{\Lag}_{H}$, the proof for the compactness of the moduli spaces of flow lines goes exactly as in Subsection 3.2 and 3.4. We already have the transversality. Hence we have that $\partial_{k-1}(\tilde{\Lag}_{H},\hat{G})\circ \partial_{k}(\tilde{\Lag}_{H},\hat{G})=0$. Writing this using matrix notation we get that
$$\Phi_{k-1}\circ \partial_{k}(\Lag_{H},G^{K_{0}}) +\partial_{k}(\Lag_{H},G^{K_{1}})\circ\Phi_{k}=0$$
Therefore, $\Phi_{*} : C_{*}(\Lag_{H},\mathbb{Z}_{2})\to C_{*}(\Lag_{H},\mathbb{Z}_{2})$ is a chain map. Notice now that there is only one flow line connecting $(x,0)$ to $(x,1)$ for all $x\in crit(\Lag_{H})$. Therefore if we order the critical points of $\Lag_{H}$ by increasing critical values of their energy, one sees that the chain map $\Phi$ can be represented by an upper triangular matrix with $1$ in the diagonal, thus it is invertible and it induces a chain isomorphism.
\hfill$\Box$

With this fact, given $H$ satisfying  $(1.1)-(1.5)$ and if we assume that $\Lag_{H}$ is Morse, then the Floer-Morse homology of the pair $(\Lag_{H},\mathcal{F}^{1,1/2}(S^{1},N))$ is defined by $HF_{*}(\Lag_{H},\mathcal{F}^{1,1/2}(S^{1},N))=HF_{*}(\Lag_{H},\mathcal{F}^{1,1/2}(S^{1},N),G_{K})$ for $K\in \mathbb{K}_{reg}$ and its isomorphism class is independent of he choice of $K\in \mathbb{K}_{reg}$.

\section{Invariance of the Morse-Floer homology, II: independence of $H$}
We consider two function $H_{0}$ and $H_{1}$ in $\mathbb{H}^{3}_{p+1}$ such that $\Lag_{H_{0}}$ and $\Lag_{H_{1}}$ are Morse. Let $\eta(t)$ be a smooth function such that $\eta(t)=1$ for $t\geq 1$ and $\eta(t)=0$ for $t\leq 0$. We set $H_{t}=\eta(t)H_{1}+(1-\eta(t))H_{0}$. We will study in this section the non-autonomous gradient flow equation defined by

\begin{equation}\label{equanon}
\frac{d\ell}{dt}=-\nabla_{G^{K}}\Lag_{H_{t}}(\ell(t)).
\end{equation}

\begin{definition}
We define the distance $d_{p+1}$ on $\mathbb{H}^{3}_{p+1}$ as follows: for $H$ and $H'$ in $\mathbb{H}^{3}_{p+1}$,
$$d_{p+1}(H,H')=\sup_{(s,\phi,\psi)\in S^{1}\times \mathbb{S}(S^{1})\otimes TN}\frac{|H(s,\phi,\psi)-H'(s,\phi,\psi)|}{1+|\psi|_{\phi}^{p+1}}.$$

\end{definition}
It is easy to see that $d_{p+1}$ is well defined on $\mathbb{H}^{3}_{p+1}$ and that if $d_{p+1}(H,H')<\varepsilon$ then
$$|H(s,\phi,\psi)-H'(s,\phi,\psi)|<\varepsilon (1+|\psi|^{p+1}).$$
From Subsection 5.2, we can always choose $K\in \mathbb{K}_{\theta,\rho}^{2}$ such that the negative gradient flows of $\Lag_{H_{0}}$ and $\Lag_{H_{1}}$ with respect to $G^{K}$ satisfy the Morse-Smale property.

\begin{proposition}
Let $H_{0}$ and $H_{1}$ in $\mathbb{H}^{3}_{p+1}$, there exists a constant $\varepsilon_{0}$ depending on the constants $C_{i}$, $1\leq i \leq 4$ such that if $d_{p+1}(H_{0},H_{1})<\varepsilon_{0}$ then for any solution of $(\ref{equanon})$ such that $\sup_{t\leq 0} \Lag_{H_{0}}(\ell(t))\leq A$, there exists $C=C(A,C_{1},C_{2},C_{3},C_{4})$ such that $$\sup_{t\geq 1} \Lag_{H_{1}}(\ell(t)) \leq C.$$
\end{proposition}
Again this will be proved through several Lemmata.
\begin{lemma}
Let $\varepsilon>0$ and $H_{1}$, $H_{0}$ in $\mathbb{H}^{3}_{p+1}$ such that $d_{p+1}(H_{0},H_{1})<\varepsilon$, then for every solution $\ell(t)$ of $(\ref{equanon})$ that satisfies $\sup_{t\leq 0}\Lag_{H_{0}}(\ell(t))\leq A$, we have
\begin{equation}\label{est1}
\Lag_{H_{t}}(\ell(t))\leq A+4\pi \varepsilon +2\varepsilon \int_{0}^{1}\|\psi(r)\|^{p+1}_{p+1}dr,
\end{equation}
for all $t\in \R$ and
\begin{equation}\label{est2}
\Lag_{H_{1}}(\ell(1))+\int_{0}^{1}\|\dot{\ell}(r)\|_{G^{K}}^{2}dr \leq A+4\pi \varepsilon +2\varepsilon \int_{0}^{1}\|\psi(r)\|^{p+1}_{p+1}dr
\end{equation}

\end{lemma}

{\it Proof:}
Notice first that for $t\geq 1$, $\Lag_{H_{t}}(\ell(t))\leq \Lag_{H_{1}}(\ell(1))$, therefore it is enough to show the proof for $0\leq t \leq 1$. Now, we have
\begin{align}
\Lag_{H}(\ell(t))&=\Lag_{H_{0}}(\ell(0))+\int_{0}^{t}\frac{d}{dr}\Lag_{H}(\ell(r))dr\notag\\
&=\Lag_{H_{0}}(\ell(0))+\int_{0}^{t}-\|\frac{d\ell}{dt}(r)\|_{G^{K}}^{2} dr +\int_{0}^{t}\eta'(r)\int_{S^{1}}(H_{1}(s,\ell(r))-H_{0}(s,\ell(r)))dsdr\notag\\
&\leq A +2\varepsilon \int_{0}^{1}\int_{S^{1}}(1+|\psi(r,s)|^{p+1})dsdr\notag\\
&\leq A+4\pi \varepsilon +2\varepsilon \int_{0}^{1}\|\psi(r)\|_{p+1}^{p+1}dr,
\end{align}
which yields $(\ref{est1})$. The proof of $(\ref{est2})$ follows by setting $t=1$ and keeping the negative term involving $\dot{\ell}$ in the second inequality.
\hfill$\Box$

\begin{lemma}
Under the assumptions of the previous Lemma, we have
\begin{equation}\label{est3}
\|\psi(t)\|_{p+1}^{p+1} \leq \frac{2A}{C_{2}} +C_{5}+ \frac{4\varepsilon}{C_{2}}\int_{0}^{1}\|\psi(r)\|_{p+1}^{p+1}dr +\frac{C'}{C_{2}}\|\dot{\psi}(t)\|_{H^{1/2}}\|\psi(t)\|_{H^{1/2}}
\end{equation}
Where $C_{5}=\frac{8\pi \varepsilon+C}{C_{2}}$ and $C'$ is the constant showing up in the equivalence of the norm $\|\cdot \|_{H^{1/2}}$ and the one induced by $G^{K}$.

\end{lemma}

{\it Proof:}
As in $(3.16)$ we have that
$$\langle \partial_{t}\psi(t),\psi(t) \rangle_{H^{1/2}} +2\Lag_{H_{t}}(\ell(t)) \geq C_{2} \|\psi(t)\|_{p+1}^{p+1} -C.$$
Hence,
$$\|\psi(t)\|_{p+1}^{p+1}\leq \frac{C}{C_{2}}+\frac{2}{C_{2}}\Lag_{H_{t}}(\ell(t))+\frac{C'}{C_{2}}\|\partial_{t}\psi(t)\|_{H^{1/2}}\|\psi(t)\|_{H^{1/2}}.$$
Using $(\ref{est1})$, we have that 
$$\|\psi(t)\|_{p+1}^{p+1}\leq \frac{C}{C_{2}}+\frac{2A}{C_{2}}+\frac{8\pi \varepsilon}{C_{2}}+\frac{4\varepsilon}{C_{2}}\int_{0}^{1} \|\psi(r)\|_{p+1}^{p+1}dr+\frac{C'}{C_{2}}\|\partial_{t}\psi(t)\|_{H^{1/2}}\|\psi(t)\|_{H^{1/2}},$$
which finishes the proof of the Lemma.

\hfill$\Box$

It follows from the previous lemma that if $\varepsilon<\frac{C_{2}}{8}$ we have that 

\begin{equation}\label{est4}
\int_{0}^{1}\|\psi(t)\|_{p+1}^{p+1}dt\leq \frac{4A}{C_{2}} +2C_{5}+\frac{2C'}{C_{2}}\int_{0}^{1}\|\partial_{t}\psi(t)\|_{H^{1/2}}\|\psi(t)\|_{H^{1/2}}dt,
\end{equation}
and
\begin{equation}\label{est5}
\|\psi(t)\|_{p+1}^{p+1} \leq C(A)+\frac{8\varepsilon C'}{C^{2}_{2}}\int_{0}^{1}\|\partial_{t}\psi(r)\|_{H^{1/2}}\|\psi(r)\|_{H^{1/2}}dr+\frac{C'}{C_{2}}\|\partial_{t}\psi(t)\|_{H^{1/2}}\|\psi(t)\|_{H^{1/2}}.
\end{equation}
In what follows $C_{k}(A)$ will be a constant depending on $C_{i}$, $1\leq i \leq 4$ and $A$ while $C_{k}$ will depend only on $C_{i}$, $1\leq i \leq 4$.

\begin{lemma}
Under the assumptions of Lemma 7.1, we have,
\begin{align}\label{est6}
\int_{0}^{1}\|\psi(t)\|^{2}_{H^{1/2}}dt \leq& C(A)+C_{13}\left(\int_{0}^{1}\|\partial_{t}\psi(t)\|_{H^{1/2}}\|\psi(t)\|_{H^{1/2}}dt+\int_{0}^{1}\|\psi(t)\|_{H^{1/2}}dt\right. \notag\\
 &\left.+\Big(\int_{0}^{1}\|\partial_{t}\psi(t)\|_{H^{1/2}}^{2}dt\Big)^{\frac{1}{2}}+\Big(\int_{0}^{1}\|\partial_{t}\psi(t)\|_{H^{1/2}}^{2}dt\Big)^{\frac{3}{4}}\Big(\int_{0}^{1}\|\psi(t)\|_{H^{1/2}}^{2}dt\Big)^{\frac{1}{4}}\right).
\end{align}
\end{lemma}

{\it Proof:}
We will use the notations of Section 3. From the inequality $(3.24)$ we have that
\begin{align}
\|\psi_{0}^{+}(t)\|^{2}_{H^{1/2}}\leq & C_{9}\Big(\|\partial_{t}\psi(t)\|_{H^{1/2}}\|\psi(t)\|_{H^{1/2}}+\|\psi(t)\|_{H^{1/2}}+\|\psi(t)\|^{p+1}_{p+1}+1\notag\\
&+\Big(1+\|\partial_{t}\psi(t)\|_{H^{1/2}}\|\psi(t)\|_{H^{1/2}}\Big)^{\frac{1}{2}}\|\partial_{t}\psi(t)\|_{H^{1/2}}\Big).
\end{align}
Using $(\ref{est5})$, we have
\begin{align}
\|\psi_{0}^{+}(t)\|^{2}_{H^{1/2}}\leq & C(A)+C_{10}\Big(\|\partial_{t}\psi(t)\|_{H^{1/2}}\|\psi(t)\|_{H^{1/2}}+\|\psi(t)\|_{H^{1/2}}+\int_{0}^{1}\|\partial_{t}\psi(r)\|_{H^{1/2}}\|\psi(r)\|_{H^{1/2}}dr\notag\\ 
&+\Big(1+\|\partial_{t}\psi(t)\|_{H^{1/2}}\|\psi(t)\|_{H^{1/2}}\Big)^{\frac{1}{2}}\|\partial_{t}\psi(t)\|_{H^{1/2}} \Big).
\end{align}
Similarly we have
\begin{align}
\|\psi_{0}^{-}(t)\|^{2}_{H^{1/2}}\leq & C(A)+C_{10}\Big(\|\partial_{t}\psi(t)\|_{H^{1/2}}\|\psi(t)\|_{H^{1/2}}+\|\psi(t)\|_{H^{1/2}}+\int_{0}^{1}\|\partial_{t}\psi(r)\|_{H^{1/2}}\|\psi(r)\|_{H^{1/2}}dr\notag\\ 
&+\Big(1+\|\partial_{t}\psi(t)\|_{H^{1/2}}\|\psi(t)\|_{H^{1/2}}\Big)^{\frac{1}{2}}\|\partial_{t}\psi(t)\|_{H^{1/2}} \Big).
\end{align}
Notice that the term $\|\psi_{0}^{0}\|_{H^{1/2}}^{2}$ can be estimated as follows
\begin{align}
\|\psi_{0}^{0}(t)\|_{H^{1/2}}^{2}&\leq C_{10}\|\psi(t)\|_{p+1}^{2}\notag\\
&\leq C(A)+C_{11}\left(\int_{0}^{1}\|\partial_{t}\psi(r)\|_{H^{1/2}}\|\psi(r)\|_{H^{1/2}}dr+\|\partial_{t}\psi(t)\|_{H^{1/2}}\|\psi(t)\|_{H^{1/2}}\right).
\end{align}
Combining these last three inequalities we get
\begin{align}\label{est8}
\|\psi(t)\|^{2}_{H^{1/2}}\leq & C(A)+C_{12}\Big(\|\partial_{t}\psi(t)\|_{H^{1/2}}\|\psi(t)\|_{H^{1/2}}+\|\psi(t)\|_{H^{1/2}}+\int_{0}^{1}\|\partial_{t}\psi(r)\|_{H^{1/2}}\|\psi(r)\|_{H^{1/2}}dr\notag\\
&+\Big(1+\|\partial_{t}\psi(t)\|_{H^{1/2}}\|\psi(t)\|_{H^{1/2}}\Big)^{\frac{1}{2}}\|\partial_{t}\psi(t)\|_{H^{1/2}}\Big).
\end{align}
After integration and the use of H\"{o}lder's inequality, we have 
\begin{align}
\int_{0}^{1}\|\psi(t)\|^{2}_{H^{1/2}}dt \leq& C(A)+C_{13}\Big(\int_{0}^{1}\|\partial_{t}\psi(t)\|_{H^{1/2}}\|\psi(t)\|_{H^{1/2}}dt+\int_{0}^{1}\|\psi(t)\|_{H^{1/2}}dt \notag\\
 &+\Big(\int_{0}^{1}\|\partial_{t}\psi(t)\|_{H^{1/2}}^{2}dt\Big)^{\frac{1}{2}}+\Big(\int_{0}^{1}\|\partial_{t}\psi(t)\|_{H^{1/2}}^{2}dt\Big)^{\frac{3}{4}}\Big(\int_{0}^{1}\|\psi(t)\|_{H^{1/2}}^{2}dt\Big)^{\frac{1}{4}}\Big).
\end{align}

\hfill$\Box$
\begin{corollary}
Under the assumptions of Lemma 7.1, we have
\begin{equation}\label{est7}
\int_{0}^{1}\|\psi(t)\|^{2}_{H^{1/2}}dt \leq C(A)+C_{14}\int_{0}^{1}\|\partial_{t}\psi(t)\|_{H^{1/2}}^{2}dt.
\end{equation}
\end{corollary}
{\it Proof:}
This corollary follows easily by applying Young's inequality first in the form $ab\leq C_{\epsilon}a^{2}+\epsilon b^{2}$ and a second time in the form $ab\leq C_{\epsilon}a^{\frac{4}{3}}+\epsilon b^{4}$. Then grouping the terms yields the inequality.
\hfill$\Box$

{\it Proof of Proposition:}
First, from $(\ref{est2})$ we have
$$
\int_{0}^{1}\|\partial_{t}\psi(t)\|_{H^{1/2}}^{2}dt\leq C(A)-2\Lag_{H_{1}}(\ell(1))+4\varepsilon \int_{0}^{1}\|\psi(t)\|^{p+1}_{p+1}dt.$$
Then from $(\ref{est4})$, it follows
$$\int_{0}^{1}\|\partial_{t}\psi(t)\|_{H^{1/2}}^{2}dt\leq C(A)-2\Lag_{H_{1}}(\ell(1))+\varepsilon C_{15} \int_{0}^{1}\|\partial_{t}\psi(t)\|_{H^{1/2}}\|\psi(t)\|_{H^{1/2}}dt.$$
Again by using Young's inequality, and taking $\varepsilon<\frac{1}{C_{15}}$ we have 
\begin{equation}\label{est8}
\int_{0}^{1}\|\partial_{t}\psi(t)\|_{H^{1/2}}^{2}dt\leq C(A)-4\Lag_{H_{1}}(\ell(1))+\varepsilon C_{16} \int_{0}^{1}\|\psi(t)\|_{H^{1/2}}^{2}dt.
\end{equation}
Combining this inequality with $(\ref{est7})$ one has
$$\int_{0}^{1}\|\psi(t)\|^{2}_{H^{1/2}}dt \leq C(A)-C_{17}\Lag_{H_{1}}(\ell(1))+\varepsilon C_{18} \int_{0}^{1}\|\psi(t)\|_{H^{1/2}}^{2}dt.$$
Hence, for $\varepsilon<\frac{1}{2C_{18}}$ we get
\begin{equation}
\int_{0}^{1}\|\psi(t)\|^{2}_{H^{1/2}}dt \leq C(A)-C_{19}\Lag_{H_{1}}(\ell(1)).
\end{equation}
Plugging it back in $(\ref{est8})$ we have
\begin{equation}
\int_{0}^{1}\|\partial_{t}\psi(t)\|_{H^{1/2}}^{2}dt\leq C(A)-C_{20} \Lag_{H_{1}}(\ell(1)).
\end{equation}
Using H\"{o}lder's inequality in $(\ref{est4})$ we have
\begin{equation}
\int_{0}^{1}\|\psi(t)\|_{p+1}^{p+1}dt\leq C(A)-C_{21}\Lag_{H_{1}}(\ell(1)).
\end{equation}
Inserting it back in $(\ref{est1})$, we get 
$$\Lag_{H_{1}}(\ell(1))\leq C(A)-2\varepsilon C_{21}  \Lag_{H_{1}}(\ell(1)).$$
Hence if we take $\varepsilon_{0}\leq \min\{\frac{C_{2}}{8},\frac{1}{C_{15}},\frac{1}{2C_{18}}\}$ we have that
$$\Lag_{H_{1}}(\ell(1))\leq C(A,C_{1},C_{2},C_{3},C_{4}).$$
\hfill$\Box$

\begin{lemma}
Let $H_{0}$ and $H_{1}$  in $\mathbb{H}_{p+1}^{3}$ and consider $\varepsilon_{0}>0$ as in Proposition 7.1. Assume that $\sup_{t\leq 0}\Lag_{H_{0}}(\ell(t))\leq A$ and $\inf_{t\geq 1}\Lag_{H_{1}}(\ell(t))\geq B$, then there exists a constant $C(A,B,C_{1},C_{2},C_{3},C_{4})$ such that the length of $\ell$, $L(\ell)$ satisfies 
$$L(\ell)\leq C(A,B).$$
\end{lemma}

{\it Proof:}
Similarly to $(\ref{est2})$ we have that 
\begin{equation}\label{esth1}
\int_{\R}\|\frac{d\ell}{dt}(t)\|_{G^{K}}^{2}dt \leq C(A,B) +2\varepsilon \int_{0}^{1}\|\psi(t)\|^{p+1}_{p+1}dt.
\end{equation}
Using the same estimates as in $(7.17-7.19)$
We get that 
$$\int_{\R}\|\frac{d\ell}{dt}(t)\|_{G^{K}}^{2}dt\leq C(A,B).$$
Now it is enough to write that
$$L(\ell)=L(\ell,(-\infty,0])+L(\ell,[1,+\infty))+L(\ell,[0,1]).$$
The first two terms are uniformly bounded as in Lemma 5.3. For the last term we have
\begin{align}
\int_{0}^{1}\|\frac{d\ell}{dt}(t)\|dt&\leq \left(\int_{0}^{1}\|\frac{d\ell}{dt}(t)\|^{2}dt\right)^{\frac{1}{2}}\notag \\
&\leq \left(\int_{\R}\|\frac{d\ell}{dt}(t)\|^{2}dt\right)^{\frac{1}{2}}\leq C(A,B).\notag
\end{align}
\hfill $\Box$

\begin{proposition}
Under the assumptions of Lemma 7.4, we have that $\mathcal{M}(x_{+},x_{-},H_{t})$ is relatively compact.
\end{proposition}

{\it Proof: }
Now as in $(3.17)$ we have
$$\|\partial_{t}\phi(t)\|^{2}_{L^{2}}\leq C(A,B)(1+\|\psi(t)\|_{H^{1/2}}\|\partial_{t}\psi(t)\|_{H^{1/2}}).$$
Integrating in $[\tau,\tau+1]$ and using H\"{o}lder's inequality, we have
$$\int_{\tau}^{\tau+1}\|\partial_{t}\phi(t)\|^{2}_{L^{2}} dt \leq C(A,B)\Bigg(1+\Big(\int_{\tau}^{\tau+1}\|\psi(t)\|^{2}_{H^{1/2}}dt\Big)^{\frac{1}{2}}\Bigg).$$
On the other hand, if we integrate $(7.13)$ we have
$$\int_{\tau}^{\tau+1}\|\psi(t)\|_{H^{1/2}}^{2}dt \leq C(A,B)\Bigg(1+\Big(\int_{\tau}^{\tau+1}\|\psi(t)\|^{2}_{H^{1/2}}dt\Big)^{\frac{1}{2}}\Bigg).$$
Therefore,
$$\sup_{\tau \in \R} \int_{\tau}^{\tau+1}\|\psi(t)\|_{H^{1/2}}^{2}dt \leq C(A,B),$$
and 
$$\sup_{\tau \in \R} \int_{\tau}^{\tau+1}\|\partial_{t}\phi(t)\|^{2}_{L^{2}} dt \leq C(A,B).$$
Now going back to $(7.7)$ and integrating, we have again
$$\sup_{\tau \in \R} \int_{\tau}^{\tau+1}\|\psi(t)\|_{p+1}^{p+1}dt \leq C(A,B).$$
Hence we have an equivalent result to the one in Lemma 3.1 which has the main estimates to prove Proposition 3.2 and Proposition 3.3. The same procedure can be repeated here using these estimates to get the desired result of the proposition.
\hfill$\Box$

We consider now two functions $H_{0}$ and $H_{1}$ in $\mathbb{H}_{p+1}^{3}$ such that $d_{p+1}(H_{0},H_{1})<\varepsilon_{0}$ and $\Lag_{H_{i}}$ is a Morse function for $i=0,1$. We want to construct a chain homomorphism between $\{C_{*}(\Lag_{H_{0}};\Z_{2}),\partial_{*}(\Lag_{H_{0}},G^{K})\}$ and $\{C_{*}(\Lag_{H_{1}};\Z_{2}),\partial_{*}(\Lag_{H_{1}},G^{K})\}$. For this purpose, we will study the moduli spaces of flow lines of the non-autonomous gradient flow. That is, for $x\in crit_{k}(\Lag_{H_{0}})$ and $y\in crit_{k}(\Lag_{H_{1}})$, we consider the solutions $\ell(t)$ of 
\begin{equation}\label{equaf}
\frac{d\ell}{dt}(t)=-\nabla_{G^{K}}\Lag_{H_{t}}(\ell(t)), \text{ and } \ell(-\infty)=x,\text{ } \ell(+\infty)=y.
\end{equation}
Again we want to perturb the metric as in the previous section in order to achieve transversality. The main trajectories here causing problems are coming from stationary paths. That is, $\ell(t)=x_{0}$ where $x_{0}\in crit_{k}(\Lag_{H_{0}})\cap crit_{k}(\Lag_{H_{1}})$. Therefore we will assume for now that 

\begin{equation}\tag{A}
 crit_{k}(\Lag_{H_{0}})\cap crit_{k}(\Lag_{H_{1}})=\emptyset.
\end{equation}
This last assumption is not as restrictive as it seems, in fact it is generic as we show in the following
\begin{lemma}
Let $H_{0}$ and $H_{1}$ in $\mathbb{H}^{3}_{p+1}$ such that $\Lag_{H_{0}}$ and $\Lag_{H_{1}}$ are Morse then there exists a residual set $\mathbb{H}_{b}(H_{1},H_{0})\subset \mathbb{H}^{3}_{b}$ such that $crit_{k}(\Lag_{H_{0}})\cap crit_{k}(\Lag_{H_{1}+h})=\emptyset$, for all $h\in \mathbb{H}_{b}(H_{1},H_{0})$.
\end{lemma}
{\it Proof:}
The proof of this result is straight forward. In fact if $x\in crit_{k}(\Lag_{H_{0}})$, the set of perturbations $h\in \mathbb{H}_{b}^{3}$ such that $x\not \in crit(\Lag_{H_{1}+h})$ is an open dense set. Indeed, this is equivalent to the equation $$\nabla_{G^{K}}\Lag_{H_{1}}(x)+\nabla_{G^{K}}\tilde{h}(x)\not =0,$$
where  $\tilde{h}(x)=\int_{S^{1}}h(s,x(s))ds$. Then since $\Lag_{H_{0}}$ is Morse, $crit(\Lag_{H_{0}})$ is countable and hence the set $\mathbb{H}_{b}(H_{1},H_{0})$ is $G_{\delta}$-dense.
\hfill$\Box$

Regarding the transversality issue, it is similar to the previous case, that is up to a perturbation of the metric $G^{K}$, the stable and unstable manifolds intersect transversally. Indeed, one starts first by taking a function $\theta :\mathcal{F}^{1,1/2}\to [0,1]$ such that $\theta=0$ on $crit(\Lag_{H_{0}})\cup crit(\Lag_{H_{1}})$ and $\theta >0$ elsewhere. Since we do not have stationary orbits, $\dot{\ell}\not =0$, hence there exists an interval $[a,b]$ such that if $\ell(t)$ is a flow line, $\ell:[a,b]\to \mathcal{F}^{1,1/2}$ is an embedding. The rest of the argument can be carried out exactly like in Section 5.2.

Now we denote by $\mathcal{M}(x,y,H_{t})\in W^{1,2}_{x,y}(\R,\mathcal{F}^{1,1/2})$ the set of solutions of $(\ref{equaf})$. From the transversality and in a similar way as in Section 3, we have that $\mathcal{M}(x,y,H_{t})$ is a manifold of dimension $\mu_{H_{0}}(x)-\mu_{H_{1}}(y)=0$. In fact, using Propositions 7.1 and 7.2, we show that $\mathcal{M}(x,y,H_{t})$ is relatively compact and it can be compactified by adding $\partial \overline{\mathcal{M}}(x,y,H_{t})$ formed by broken trajectories. These trajectories can be described as follow: there exist $x=x_{0},...,x_{p}\in crit(\Lag_{H_{0}})$, $\ell_{i} \in \mathcal{M}(x_{i},x_{i+1},H_{0})$, with $0\leq i \leq p-1$ and $y_{1},..,y_{q}=y \in  crit(\Lag_{H_{1}})$, $\overline{\ell}_{i}\in \mathcal{M}(y_{i},y_{i+1},H_{1})$ for $1\leq i \leq q-1$ such that 
$$\mu_{H_{0}}(x_{0})>\mu_{H_{0}}(x_{1})>...>\mu_{H_{0}}(x_{p})\geq \mu_{H_{1}}(y_{1})>\mu_{H_{1}}(y_{2})>...>\mu_{H_{1}}(y_{q}).$$
Therefore in our case, since $\mu_{H_{0}}(x)=\mu_{H_{1}}(y)$, there is no broken trajectories and the manifold is just a finite number of points. Hence we denote $$n_{0,1}(x,y)=\sharp \mathcal{M}(x,y,H_{t})(\text{mod 2}).$$
We define therefore the homomorphism $\Phi_{01}:C_{k}(\Lag_{H_{0}},\Z_{2})\to C_{k}(\Lag_{H_{1}},\Z_{2})$ on the generators by
\begin{equation}
\Phi_{01}(x)=\sum_{y\in crit_{k}(\Lag_{H_{1}})}n_{01}(x,y)y.
\end{equation}
\begin{lemma}
 Let $H_{1}$ and  $H_{0}$ in $\mathbb{H}^{3}_{p+1}$ such that $\Lag_{H_{1}}$ and $\Lag_{H_{0}}$ are Morse and $d_{p+1}(H_{0},H_{1})<\varepsilon_{0}$, then under the assumption $(A)$, the homomorphism $\Phi_{01}:\{C_{*}(\Lag_{H_{0}},\Z_{2}),\partial_{*}(\Lag_{H_{0}}G^{K})\}\to \{C_{*}(\Lag_{H_{1}},\Z_{2}),\partial_{*}(\Lag_{H_{1}}G^{K})\}$ defines a chain homomorphism.

\end{lemma}
{\it Proof:}
We consider two arbitrary critical points $x_{0,k}\in crit_{k}(\Lag_{H_{0}})$ and $x_{1,k-1}\in crit_{k-1}(\Lag_{H_{1}})$. By transversality (up to perturbation of the metric $G^{K}$ ) we have that $\mathcal{M}(x_{0,k},x_{1,k-1},H_{t})$ is a one dimensional manifold that is precompact and its compactification $\overline{\mathcal{M}}(x_{0,k},x_{1,k-1},H_{t})$ consists of the following two types of orbits :
\begin{itemize}
\item[i)] Broken trajectories connecting $x_{0,k}$ to $x_{0,k-1}\in crit_{k-1}(\Lag_{H_{0}})$, that is flow lines $\ell_{0}\in \mathcal{M}(x_{0,k},x_{0,k-1},H_{0})$ then $x_{0,k-1}$ to $x_{1,k-1}$, that is a flow line $\ell_{01}\in \mathcal{M}(x_{0,k-1},x_{1,k-1},H_{t})$.
\item[ii)]Broken trajectories connecting $x_{0,k}$ to $x_{1,k}\in crit_{k}(\Lag_{H_{1}})$, that is flow lines $\ell_{01}\in \mathcal{M}(x_{0,k},x_{1,k},H_{t})$ then $x_{1,k}$ to $x_{1,k-1}$, that is a flow line $\ell_{1}\in \mathcal{M}(x_{1,k},x_{1,k-1},H_{1})$.
\end{itemize}
This of course require a gluing statement with the compactness result ( Proposition 7.1 and 7.2). But this relies mainly on the Fredholm property of the linearized flow equation, which hold by transversality and the argument of Subsection 3.5 applies to prove gluing. This yields in particular that since this boundary is the one of a one dimensional manifold, its cardinal is zero mod 2. So in terms of counting we have
\begin{align*}
&\sum_{x_{0,k-1}\in crit_{k-1}(\Lag_{H_{0}})}n_{0}(x_{0,k},x_{0,k-1})n_{01}(x_{0,k-1},x_{1,k-1})\\
&+ \sum_{x_{1,k}\in crit_{k}(\Lag_{H_{1}})}n_{01}(x_{0,k},x_{1,k})n_{1}(x_{1,k},x_{1,k-1})=0\,(\text{mod 2})
\end{align*}
and this is equivalent to
$$\Phi_{01}\circ \partial_{k}(\Lag_{H_{0}},G^{K})=\partial_{k}(\Lag_{H_{1}},G^{K})\circ \Phi_{01}.$$
\hfill$\Box$

Now we will prove a naturality result for the homomorphism $\Phi_{01}$ at the homology level.
\begin{proposition}
There exists $0<\varepsilon<\varepsilon_{0}$ such that if we consider three function $H_{0}$, $H_{1}$ and $H_{2}$ in $\mathbb{H}^{3}_{p+1}$ with $\Lag_{H_{i}}$ is Morse for $0\leq i \leq 2$ and $d_{p+1}(H_{i},H_{j})<\varepsilon$ for $0\leq i,j\leq 2$. We assume furthermore that $(A)$ holds for the pair $H_{i}$, $H_{j}$ with $i\not =j$. Then for a generic $K\in \mathbb{K}^{2}_{\rho,\theta}$, with $\theta:\mathcal{F}^{1,1/2} \to [0,1]$ is adequately chosen and adapted to the three functions $H_{i}$, the chain homomorphism $\Phi_{ij}:\{C_{*}(\Lag_{H_{i}},\Z_{2}),\partial_{*}(\Lag_{H_{i}}G^{K})\}\to \{C_{*}(\Lag_{H_{j}},\Z_{2}),\partial_{*}(\Lag_{H_{j}}G^{K})\}$ are well defined (Lemma 7.4) and satisfies the following :
\begin{itemize}
\item[1)] $\Phi_{00}=id$.
\item[2)] $\Phi_{02}=\Phi_{12}\circ\Phi_{01}$ in the homology level.
\end{itemize}

\end{proposition}

{\it Proof of statement 1)}

First , by suitably choosing $\theta$, we can assume that for a generic $K\in \mathbb{K}^{2}_{\rho,\theta}$, the Morse-Smale, property up to order 2 holds for the negative gradient flows of $\Lag_{H_{i}}$, $i=0,1,2$. To prove $1)$ in the Proposition, we construct $\Phi_{00}$ by taking the constant homotopy $H_{t}(t,s,x)=H_{0}(s,x)$. In this case the flow $(\ref{equanon})$ becomes the negative gradient flow of $\Lag_{H_{0}}$ which satisfies the Morse-Smale property up to order 2. Hence, for $x,y\in crit_{k}(\Lag_{H_{0}})$, $\mathcal{M}(x,y,H_{t})=\emptyset$ unless $x=y$. In that case it contains the constant Flow line $\ell(t)=x$, therefore $n_{00}(x,x)=1$ and $n_{00}(x,y)=0$ for $x\not =y$. Hence, $\Phi_{00}=id$.
\hfill$\Box$

Now we move to the proof of the second point, which is more involved and technical. We consider the function $H_{ij}$ defined by
$$H_{ij}(t,s,x)=\eta(t)H_{j}(s,x)+(1-\eta(t))H_{i}(s,x),$$
for $0\leq i,j \leq 2$ and also for $R>1$,
$$H_{02,R}(t,s,x)=\eta(t)H_{12}(t-R,s,x)+(1-\eta(t))H_{01}(t+R,s,x).$$
We consider now the flow
 \begin{equation}\label{equar}
\frac{d\ell}{dt}=-\nabla_{G^{K}}\Lag_{H_{02,R}} (\ell(t)).
\end{equation}
Since the flow is autonomous except in the intervals $[-R,-R+1]$ and $[R,R+1]$, all the compactness estimates hold as in Propositions 7.1 and 7.2. Indeed it is easy to notice that 
$$\Big|\frac{\partial H_{02,R}}{\partial t}\Big|\leq 4\varepsilon (1+|\psi|^{p+1}).$$ 
Therefore if $\varepsilon<\frac{\varepsilon_{0}}{4}$ we have the results of the previously cited propositions by integrating between $[R, R+1]$ instead of $[0,1]$ in the proofs. Moreover, the estimates are independent of $R$ since the epsilon that we picked does not depend in $R$. Therefore one then can define the chain homomorphism $\Phi_{02,R}:\{C_{*}(\Lag_{H_{0}},\Z_{2}),\partial_{*}(\Lag_{H_{0}}G^{K})\}\to \{C_{*}(\Lag_{H_{2}},\Z_{2}),\partial_{*}(\Lag_{H_{2}}G^{K})\}$ in the same way as in Lemma 7.5. Notice that we did not exclude the case $H_{i}=H_{j}$ for some $i$ and $j$, in that case we are back at the situation where we have only two functions $H_{0}$ and $H_{1}$ such that $H_{0}\not =H_{1}$, in view of the assumption (A), we have that for any flow line, solution to $(\ref{equar})$ linking two critical points of the same index, $\frac{d\ell}{dt}\not=0$. Hence nothing changes in the construction of the morphism $\Phi_{02,R}$
\hfill$\Box$

\begin{lemma}
The chain homomorphism $\Phi_{02,R}$ is chain homotopy equivalent to $\Phi_{02}$.
\end{lemma}

{\it Proof:}
We consider again the following homotopy defined for $0\leq r \leq 1$ by 
$$H_{r}(t,s,x)=rH_{02,R}(t,s,x)+(1-r)H_{02}(t,s,x).$$
For $x_{0}\in crit(\Lag_{H_{0}})$ and $x_{2}\in crit(\Lag_{H_{2}})$, we define the operator $\mathcal{F}_{x_{0},x_{2}}:[0,1]\times W^{1,2}_{x_{0},x_{2}}(\R,\mathcal{F}^{1,1/2}) \to L^{2}(\R,T\mathcal{F}^{1,1/2})$ by
$$\mathcal{F}_{x_{0},x_{2}}(r,\ell)=\frac{d\ell}{dt}+\nabla_{G^{K}}\Lag_{H_{r}} (\ell(t))$$
We may assume, by a perturbation again of $G^{K}$, that $0$ is a regular value of $\mathcal{F}_{x_{0},x_{2}}$. Also, notice that 
$$\Big|\frac{\partial H_{r}}{\partial t}\Big|\leq 6\varepsilon (1+|\psi|^{p+1}).$$ 
Therefore, for $\varepsilon <\frac{\varepsilon_{0}}{6}$, the compactness estimates of Propositions 7.1 and 7.2 hold, since the flow is autonomous except on the intervals of length 1, $[-R,-R+1]$, $[0,1]$ and $[R, R+1]$. It follows then that the set $\mathcal{M}(x_{0},x_{2}, H_{r})=\mathcal{F}_{x_{0},x_{2}}^{-1}(0) \subset [0,1]\times W^{1,2}_{x_{0},x_{2}}(\R,\mathcal{F}^{1,1/2})$, is  relatively compact. Moreover, by transversality, we have that 
$$dim \mathcal{M}(x_{0},x_{2}, H_{r})=\mu_{H_{0}}(x_{0})-\mu_{H_{2}}(x_{2})+1.$$
We focus now on the zero dimensional case. So we take $x_{0,k}\in crit_{k}(\Lag_{H_{0}})$ and $x_{2,k+1} \in crit_{k+1}(\Lag_{H_{2}})$. By compactness, we can define the number $\tilde{n}(x_{0,k},x_{2,k+1})$ as the cardinal of $\mathcal{M}(x_{0,k},x_{2,k+1}, H_{r})$ mod $2$. That is,
$$\tilde{n}(x_{0,k},x_{2,k+1})=\sharp \mathcal{M}(x_{0,k},x_{2,k+1}, H_{r}) (\text{mod 2}).$$
Then we define the morphism $\tilde{\Phi}:C_{k}(\Lag_{H_{0}},\Z_{2})\to C_{k+1}(\Lag_{H_{0}},\Z_{2})$ on the generators  by 
\begin{equation}
\tilde{\Phi}(x_{0,k})=\sum_{x_{2,k+1}\in crit_{k+1}(\Lag_{H_{2}})}\tilde{n}(x_{0,k},x_{2,k+1})x_{2,k+1}.
\end{equation}
Notice that by the usual compactness and gluing argument, we have that $\mathcal{M}(x_{0,k},x_{2,k}, H_{r})$ is a one dimensional manifold and has a natural compactification to $\overline{\mathcal{M}}(x_{0,k},x_{2,k}, H_{r})$ and the boundary $\partial \mathcal{M}(x_{0,k},x_{2,k+1}, H_{r})$ consists of broken trajectories as follows 
\begin{itemize}
\item Flow lines connecting $x_{0,k}$ to $x_{2,k+1}$ via the flow $(\ref{equar})$ then $x_{2,k+1}$ to $x_{2,k}$ via the negative gradient flow of $\Lag_{H_{2}}$.
\item Flow lines connecting $x_{0,k}$ to $x_{0,k-1}$ via the negative gradient flow of $\Lag_{H_{0}}$ then $x_{0,k-1}$ to $x_{2,k}$ via the the flow $(\ref{equar})$.
\item Flow lines connecting $x_{0,k}$ to $x_{2,k}$ via the flow of the vector field $-\nabla_{G^{K}}\Lag_{H_{02}}$.

\item Flow lines connecting $x_{0,k}$ to $x_{2,k}$ via the flow of the vector field $-\nabla_{G^{K}}\Lag_{H_{02,R}}$.

\end{itemize}
Therefore the summation 
\begin{align}
&\sum_{x_{2,k+1}\in crit_{k+1}(\Lag_{H_{2}})}\tilde{n}(x_{0,k},x_{2,k+1})n_{2}(x_{2,k+1},x_{2,k})\notag\\
&+\sum_{x_{0,k-1}\in crit_{k-1}(\Lag_{H_{0}})}n_{0}(x_{0,k},x_{0,k-1})\tilde{n}(x_{0,k-1},x_{2,k}) \notag \\
&+n_{02}(x_{0,k},x_{2,k}) + n_{02,R}(x_{0,k},x_{2,k}) = 0(\text{mod 2})
\end{align}

But notice that 
\begin{align}
&\left(\partial_{k+1}(\Lag_{H_{2}},G^{K})\circ \tilde{\Phi} +\tilde{\Phi}\circ \partial_{k}(\Lag_{H_{0}},G^{K})+\Phi_{02}+\Phi_{02,R}\right)(x_{0,k})=\notag \\
&\sum_{x_{2,k}\in crit_{k}(\Lag_{H_{2}})}\sum_{x_{2,k+1}\in crit_{k+1}(\Lag_{H_{2}})}\tilde{n}(x_{0,k},x_{2,k+1})n_{2}(x_{2,k+1},x_{2,k})x_{2,k}\notag\\
&+\sum_{x_{2,k}\in crit_{k}(\Lag_{H_{2}})}\sum_{x_{0,k-1}\in crit_{k-1}(\Lag_{H_{0}})}n_{0}(x_{0,k},x_{0,k-1})\tilde{n}(x_{0,k-1},x_{2,k})x_{2,k} \notag \\
&+\sum_{x_{2,k}\in crit_{k}(\Lag_{H_{2}})} n_{02}(x_{0,k},x_{2,k})x_{2,k} + \sum_{x_{2,k}\in crit_{k}(\Lag_{H_{2}})}n_{02,R}(x_{0,k},x_{2,k})x_{2,k}
\end{align}
Therefore, from (7.25), we have
$$\partial_{k+1}(\Lag_{H_{2}},G^{K})\circ \tilde{\Phi} +\tilde{\Phi}\circ \partial_{k}(\Lag_{H_{0}},G^{K})+\Phi_{02}+\Phi_{02,R}=0,$$
Hence $\tilde{\Phi}$ is a chain homotopy between $\Phi_{0,2}$ and $\Phi_{02,R}$.
\hfill$\Box$

\begin{lemma}
Consider an arbitrary $x_{0,k}\in crit(\Lag_{H_{0}})$, then for $R>1$ large enough, $$\Phi_{02,R}(x_{0,k})=\Phi_{12}\circ \Phi_{01}(x_{0,k}).$$
\end{lemma}
{\it Proof:}
The idea of the proof lies mainly in studying the moduli space $\mathcal{M}(x_{0,k},x_{2,k},H_{02,R})$ for $x_{0,k}\in crit_{k}(\Lag_{H_{0}})$ and $x_{2,k}\in crit_{k}(\Lag_{H_{2}})$, when $R\to \infty$. It is similar to the case of Floer homology in \cite{AD} and \cite{S}. We will give a sketch of it in our case. First, using the same compactness estimates as in Propositions 7.1 and 7.2, we can show that the space $\mathcal{M}(x_{0,k},x_{2,k},H_{02,R_{n}})$ is relatively compact in the Hausdorff topology when $R_{n}\to \infty$. Moreover its limit consists of the following trajectories:
\begin{itemize}
\item[i)] The gradient flow lines of $\Lag_{H_{0}}$
\item[ii)] The gradient flow lines of $\Lag_{H_{01}}$
\item[iii)] The gradient flow lines of $\Lag_{H_{1}}$
\item[iv)] The gradient flow lines of $\Lag_{H_{12}}$
\item[v)] The gradient flow lines of $\Lag_{H_{2}}$
\end{itemize}
Notice that in the cases $i)$, $iii)$ and $v)$ the flow lines are constant because of the index restriction. So the important flow lines are in the cases $ii)$ and $iv)$. And it can be shown that $\ell_{n}(\cdot-R_{n})\to \ell_{01}$ and $\ell_{n}(\cdot+R_{n})\to \ell_{12}$ as $R_{n}\to \infty$, where $\ell_{01}\in \mathcal{M}(x_{0,k},x_{1,k},H_{01})$ and $\ell_{12}\in \mathcal{M}(x_{1,k},x_{2,k},H_{12})$ for some $x_{1,k}\in crit_{k}(\Lag_{H_{1}})$. Conversely, given broken trajectory from $x_{0,k}$ to $x_{1,k}$ then from $x_{1,k}$ to $x_{2,k}$, by the compactness argument and transversality, we can do the same gluing construction to associate for $R>1$ big enough an element in $\mathcal{M}(x_{0,k},x_{2,k},H_{02,R})$. Therefore we have the following identification for $R$ large enough :
\begin{equation}
\bigcup_{x_{1,k}\in crit_{k}(\Lag_{H_{1}})}\mathcal{M}(x_{0,k},x_{1,k},H_{01})\times \mathcal{M}(x_{1,k},x_{2,k},H_{12}) \to \mathcal{M}(x_{0,k},x_{2,k},H_{02,R}).
\end{equation}
It follows that 
\begin{equation}
n_{02,R}(x_{0,k},x_{2,k})=\sum_{x_{1,k}\in crit_{k}(\Lag_{H_{1}})}n_{01}(x_{0,k},x_{1,k})n_{12}(x_{1,k},x_{2,k})
\end{equation}
and 
\begin{align}
\Phi_{02,R}(x_{0,k})&=\sum_{x_{2,k}\in crit_{k}(\Lag_{H_{2}})}\sum_{x_{1,k}\in crit_{k}(\Lag_{H_{1}})}n_{01}(x_{0,k},x_{1,k})n_{12}(x_{1,k},x_{2,k})x_{2,k}\notag\\
&=\Phi_{12}\circ \Phi_{01}(x_{0,k}).
\end{align}
Which yields the proof of the lemma.
\hfill$\Box$

Now we go back to the proof of Proposition 7.3. Consider a cycle $c_{0,k}=\sum_{i}a_{i}x_{0,k}^{i} \in C_{k}(\Lag_{H_{0}},\Z_{2})$, from Lemma 7.6, we have that 
$$\Phi_{02}([c_{0,k}])=\Phi_{02,R}([c_{0,k}]).$$
Now for $R>1$ big enough, possibly depending on $c_{0,k}$, we have from Lemma 7.5
$$\Phi_{02,R}([c_{0,k}])=\Phi_{12}\circ\Phi_{01}([c_{0,k}])$$
and this finishes the proof of the second point of Proposition 7.3.\\

As a corollary from the previous construction, we can drop the assumption $(A)$ to define a homomorphism $\Phi_{01}:HF_{*}(\Lag_{H_{0}},\mathcal{F}^{1,1/2},\Z_{2})\to HF_{*}(\Lag_{H_{1}},\mathcal{F}^{1,1/2},\Z_{2})$. Indeed, we consider $H_{0}$ and $H_{1}$ such that $d_{p+1}(H_{0},H_{1})<\varepsilon$, then we can perturb $H_{1}$ by $h\in \mathbb{H}_{b}^{3}$ so that the pairs $(H_{0},H)$ and $(H,H_{1})$ satisfy $(A)$ for $H=H_{1}+h$ and such that $d_{p+1}(H_{0},H)<\varepsilon$, $d_{p+1}(H,H_{1})<\varepsilon$. Then by Proposition 7.3, we have that $\Phi_{H_{0},H}:HF_{*}(\Lag_{H_{0}},\mathcal{F}^{1,1/2},\Z_{2})\to HF_{*}(\Lag_{H},\mathcal{F}^{1,1/2},\Z_{2})$ and $\Phi_{H,H_{1}}:HF_{*}(\Lag_{H},\mathcal{F}^{1,1/2},\Z_{2})\to HF_{*}(\Lag_{H_{1}},\mathcal{F}^{1,1/2},\Z_{2})$ are well defined and hence we can define $\Phi_{01}=\Phi_{H,H_{1}}\circ \Phi_{H_{0},H}$. 
We claim that $\Phi_{01}$ defined this way is independent of the perturbation $H$. Indeed, consider another function $H'=H_{1}+h'$ such that $d_{p+1}(H_{0},H')<\varepsilon$ and $d_{p+1}(H',H_{1})<\varepsilon$, then by Lemma 7.4, there exists two generic sets $\mathbb{H}(H)$ and $\mathbb{H}(H')$ in $\mathbb{H}_{b}^{3}$ such that for $k\in\mathbb{H}(H)$ and $k'\in\mathbb{H}(H')$ the condition $(A)$ is satisfied for the pairs $(H_{0},H_{1}+k)$, $(H_{1},H_{1}+k)$, $(H,H_{1}+k)$, $(H_{0},H_{1}+k')$, $(H_{1},H_{1}+k')$ and $(H',H_{1}+k')$. Then we have in the homology level
\begin{align}
\Phi_{H,H_{1}}\circ \Phi_{H_{0},H}&= \left(\Phi_{H_{1}+k,H_{1}}\circ \Phi_{H,H_{1}+k}\right)\circ \left(\Phi_{H_{1}+k,H}\circ \Phi_{H_{0},H_{1}+k}\right)\notag\\
&=\Phi_{H_{1}+k,H_{1}}\circ \left(\Phi_{H,H_{1}+k}\circ \Phi_{H_{1}+k,H}\right)\circ \Phi_{H_{0},H_{1}+k}\notag\\
&=\Phi_{H_{1}+k,H_{1}}\circ\Phi_{H_{0},H_{1}+k}.
\end{align}
Similarly we have
\begin{equation}
\Phi_{H',H_{1}}\circ \Phi_{H_{0},H'}=\Phi_{H_{1}+k',H_{1}}\circ\Phi_{H_{0},H_{1}+k'}.
\end{equation}
Now notice that since the sets $\mathbb{H}(H)$ and $\mathbb{H}(H')$ are generic, then $\mathbb{H}(H)\cap \mathbb{H}(H')$ is also generic, hence we can pick $k=k'\in  \mathbb{H}(H)\cap \mathbb{H}(H')$ and this leads to 
$$\Phi_{H,H_{1}}\circ \Phi_{H_{0},H}=\Phi_{H',H_{1}}\circ \Phi_{H_{0},H'}.$$
In a similar way we can show that the naturality also holds, that is $\Phi_{00}=id$ and $\Phi_{02}=\Phi_{12}\circ \Phi_{01}$. Therefore we can write the following
\begin{corollary}
Let $H_{0},H_{1}$ and $H_{2}$ in $\mathbb{H}^{3}_{p+1}$ such that $\Lag_{H_{i}}$ is Morse for $i=0,1,2$ and $d_{p+1}(H_{i},H_{j})<\varepsilon$ for $0\leq i,j \leq 2$. Then we have a natural isomorphism 
$$\Phi_{ij}:HF_{*}(\Lag_{H_{i}},\mathcal{F}^{1,1/2}(S^{1},N),\Z_{2})\to HF_{*}(\Lag_{H_{j}},\mathcal{F}^{1,1/2}(S^{1},N),\Z_{2}).$$
Moreover,
$$\Phi_{00}=id \text{ and } \Phi_{02}= \Phi_{12}\circ \Phi_{01}.$$
\end{corollary}

This allows us to see that the homology $HF_{*}(\Lag_{H},\mathcal{F}^{1,1/2}(S^{1},N),\Z_{2})$ can be defined even when $\Lag_{H}$ is not a Morse function. Indeed, if $H\in \mathbb{H}^{3}_{p+1}$, we know from SubSection 5.1 that there exists a generic $h\in \mathbb{H}_{b}^{3}$ such that $d_{p+1}(H,H+h)<\varepsilon$ and that $\Lag_{H+h}$ is Morse. Hence we define $$HF_{*}(\Lag_{H},\mathcal{F}^{1,1/2}(S^{1},N),\Z_{2}):=HF_{*}(\Lag_{H+h},\mathcal{F}^{1,1/2}(S^{1},N),\Z_{2}).$$
From Corollary 7.2, this construction is independent from the generic perturbation $h$, therefore the identification is well defined for any function $H$ in $\mathbb{H}^{3}_{p+1}$ and this proves the first point of Theorem 1.1.

Now to prove the second point, we need the following Lemma.
\begin{lemma}
The metric space $(\mathbb{H}^{3}_{p+1},d_{p+1})$ is contractible (in fact, it is convex).
\end{lemma}
{\it Proof:}
We fix $H_{0}\in \mathbb{H}_{p+1}^{3}$ and we consider the map $\mathcal{H}:[0,1]\times \mathbb{H}_{p+1}^{3} \to \mathbb{H}_{p+1}^{3}$ defined by 
$$\mathcal{H}(\lambda,H)=\lambda H_{0}+(1-\lambda)H.$$
Now this map is well defined since the assumptions $(1.1)-(1.3)$ are stable under convex combination. In remains to show that $\mathcal{H}$ is continuous. Indeed,
\begin{align}
&|\mathcal{H}(\lambda,H)-\mathcal{H}(\lambda',H')|\notag \\
&\leq |\mathcal{H}(\lambda,H)-\mathcal{H}(\lambda,H')|+|\mathcal{H}(\lambda,H')-\mathcal{H}(\lambda',H')|\notag \\
&\leq |\lambda-\lambda'||H_{0}-H|+|1-\lambda'||H-H'|\notag \\
&\leq |\lambda-\lambda'|(1+|\psi|^{p+1})d_{p+1}(H_{0},H)+(1+|\psi|^{p+1})d_{p+1}(H,H') 
\end{align}
This yields to 
$$d_{p+1}(\mathcal{H}(\lambda,H),\mathcal{H}(\lambda',H'))\leq |\lambda-\lambda'|d_{p+1}(H_{0},H)+d_{p+1}(H,H'),$$
which gives us the continuity and clearly $H$ is a contraction from $\mathbb{H}_{p+1}^{3}$ to $\{H_{0}\}$.
\hfill$\Box$

{\it Proof of Theorem 1.1 part 2)}

We consider two functions $H$ and $H'$ in $\mathbb{H}_{p+1}^{3}$. Since this last space is contractible, it is connected in particular, so we take a continuous path $(H_{\lambda})_{0\leq \lambda \leq 1}$ in $\mathbb{H}_{p+1}^{3}$ connection $H$ to $H'$. We fix $\varepsilon>0$ as in Proposition 7.3. Then there exists a subdivision $0=\lambda_{0}<\lambda_{1}<\cdots<\lambda_{n}=1$ such that $d_{p+1}(H_{\lambda_{i-1}},H_{\lambda_{i}})<\varepsilon$ for $1\leq i \leq n$. From Corollary 7.2 we know that there exist natural isomorphisms $\Phi_{i-1,i}:HF_{*}(\Lag_{H_{i-1}},\mathcal{F}^{1,1/2},\Z_{2})\to HF_{*}(\Lag_{H_{i}},\mathcal{F}^{1,1/2},\Z_{2})$. We define then the isomorphism $\Phi_{H,H'}=\Phi_{n-1,n}\circ \cdots \circ \Phi_{1,2} \circ \Phi_{0,1}$. The isomorphism $\Phi_{H,H'}$ only depends on the homotopy class of $H_{\lambda}$, but since the space $\mathbb{H}_{p+1}^{3}$ is contractible, we have that this isomorphism is natural and this finishes the proof of the second point of Theorem 1.1.

\hfill$\Box$

\section{Computation of the homology}
We consider the following functionals $E_{1} : H^{1}(S^{1},N)\to \R$ and $E_{2}:\mathcal{F}^{1,1/2}(S^{1},N)\to \R$ defined by
$$E_{1}(\phi)=\frac{1}{2}\int_{S^{1}}|\dot{\phi}(s)|^{2}ds-v(\phi),$$
and
$$E_{2}(\phi,\psi)=\frac{1}{2}\int_{S^{1}}\langle \D_{\phi} \psi (s), \psi(s) \rangle ds - \int_{S^{1}}H(s,\phi(s),\psi(s))ds,$$
where $v$ is a bounded smooth function and $H\in \mathbb{H}^{3}_{p+1}$. Notice that, if $\pi : \mathcal{F}^{1,1/2}(S^{1},N)\to H^{1}(S^{1},N)$ is the canonical projection, then 
$$E_{1}\circ \pi + E_{2} =\Lag_{H+v}.$$
Without loss of generality, we will write $E_{1}$ instead of $E_{1}\circ \pi$. We consider then the vector field $V$ on $T\mathcal{F}^{1,1/2}$ defined by 
$$V(x)=-\Big(\nabla^{H}E_{1}(x)+\nabla^{N}E_{2}(x)\Big),$$
where $\nabla^{H}$ stands for the horizontal gradient, that is $\nabla_{\phi}$ and $\nabla^{N}$ is the vertical gradient, that is the gradient on the fibers, i.e. $\nabla_{\psi}$. So $V$ can be written as

\begin{equation}
V(x)=-\left(\begin{array}{cc}
(-\Delta+1)^{-1}\Big(-\nabla_{s}\partial_{s}\phi-\nabla_{\phi}v(s,\phi)\Big)\\
(1+|\D|)^{-1}\Big(\D_{\phi}-\nabla_{\psi}H(s,\phi,\psi)\Big)
\end{array}
\right).
\end{equation}
From this decomposition, we see that we have the following

\begin{lemma}
The rest points of the vector field $V$ are of the form $(\phi,\psi)$ where $\phi$ is a perturbed geodesic satisfying
$$-\nabla_{s}\partial_{s}\phi=\nabla_{\phi}v(s,\phi),$$
and $\psi$ is a critical point of the functional $E_{2}$ restricted to the fiber satisfying
$$\D_{\phi}\psi = \nabla_{\psi} H(s,\phi,\psi).$$
\end{lemma}

 We will denote by $rest(V)$ the rest points of the vector field $V$. Now notice that if $(\phi,\psi)\in rest(V)$ then $\phi$ has a natural index as a critical point of the Morse function $E_{1}$. Similarly, for $\psi$, we can define a relative index for it. Indeed, by considering the differential $dV$ we see that we can write

\begin{equation}
dV(\phi,\psi)=\begin{pmatrix}\mathcal{A}_{11}&0\\
\mathcal{A}_{21}&\mathcal{A}_{22}\end{pmatrix}
\end{equation}
so that for all $u=(X,\xi)\in T_{\phi,\psi}\mathcal{F}^{1,1/2}(S^{1},N)$,
\begin{align} 
\mathcal{A}_{11}[X]&=\nabla_{\phi}^{2}E_{1}[X]\notag \\
&=(-\nabla_{s}\nabla_{s}+1)^{-1}\Big(-\nabla_{s}\nabla_{s}X-R(X,\partial_{s}\phi)\partial_{s}\phi-v_{\phi\phi}(\,\cdot\,,\phi)[X]\Big),
\end{align}

\begin{align}
&\mathcal{A}_{21}[X]=(1+|\D|)^{-1}\Big(e_{1}\cdot\partial_{s}\psi^{i}X^{m}\gamma_{mj}^{k}\otimes\frac{\partial}{\partial y^{k}}(\phi)+\nabla_{X}\Gamma_{jk}^{i}(\phi)\partial_{s}\phi^{j}e_{1}\cdot\psi^{k}\otimes\frac{\partial}{\partial y^{i}}(\phi )\notag\\
&+\Gamma_{jk}^{i}(\phi)\partial_{s}X^{j}e_{1}\cdot\psi^{k}\otimes\frac{\partial}{\partial y^{i}}(\phi)+\Gamma_{jk}^{i}(\phi)\partial_{s}\phi^{j}e_{1}\cdot\psi^{k}\Gamma_{li}^{m}X^{l}\otimes\frac{\partial}{\partial y^{m}}(\phi)-H_{\psi\phi}(\,\cdot\,,\phi,\psi)[X]\Big),
\end{align}
\begin{align}
\mathcal{A}_{22}[\xi]&=\nabla^{2}_{\psi,\psi}E_{2}\notag \\
&=(1+|\D|)^{-1}(\D_{\phi}\xi-H_{\psi\psi}(\,\cdot\,\phi,\psi)[\xi]).
\end{align}
Again we use the same indexing process as in Section 2, that is, we fix $x_{0}=(\phi_{0},\psi_{0})\in \mathcal{F}^{1,1/2}$ and a path $x_{t}=(\phi_{t},\psi_{t})_{0\leq t \leq 1}$ connecting it to $x=(\phi,\psi)\in rest(V)$. The relative Morse index of $x$ is then
\begin{equation}
\mu_{0}(x)=-\mathsf{sf}\{dV(x_{t})\}_{0\le t\le 1}.
\end{equation}
Notice then by the above formula given for $dV$, that
$$
\mathsf{sf}\{dV(x_{t})\}_{0\le t\le 1}=\mathsf{sf}\{\mathcal{A}_{11}(x_{t})\}_{0\le t\le 1}+\mathsf{sf}\{\mathcal{A}_{22}(x_{t})\}_{0\le t\le 1}$$

Therefore, we have that 
\begin{equation}
\mu_{0}(x)=\mu_{1}(\phi)+\mu_{2}(\phi,\psi)
\end{equation}
where
$$\mu_{1}(\phi)=-\mathsf{sf}\{\mathcal{A}_{11}(x_{t})\}_{0\le t\le 1} \text{ and } \mu_{2}=-\mathsf{sf}\{\mathcal{A}_{22}(x_{t})\}_{0\le t\le 1}.$$
Notice that if we choose $\phi_{0}$ to be a minimizer of $E_{1}$, then $\mu_{1}$ coincides with the Morse index of the $\phi$ as a critical point of $E_{1}$. Based on this splitting of the index, we can define then
\begin{definition}
For $(r,q)\in \N \times \Z$ we define the set $rest_{r,q}(V)\subset rest(V)$ by
$$rest_{r,q}(V)=\{(\phi,\psi)\in rest(V);\mu_{1}(\phi)=r, \mu_{2}(\phi,\psi)=q\},$$
and $$rest_{k}(V)=\{(\phi,\psi)\in rest(V);\mu_{0}(\phi,\psi)=k\}=\bigcup_{r+q=k}rest_{r,q}(V).$$

\end{definition}

\begin{lemma}
The vector field $V$ satisfies the generalized Palais-Smale condition. That is if $V(x_{k})\to 0$ and $|E_{1}(\phi_{k})|+|E_{2}(x_{k})|\leq C$ then $x_{k}$ have a convergent subsequence.
\end{lemma}
{\it Proof:}
This proof is similar to the classical (PS) condition for the functional $\Lag_{H}$. In fact it is even easier since we have the compactness and convergence for the $\phi_{k}$ part of $x_{k}$ so the part that needs to be checked is for the $\psi$ part. But this follows from the local trivialization around $\phi$ and it is similar in nature to proof of the Palais-Smale in \cite{I}.

\hfill$\Box$

It is a classical result that for a generic set of perturbations $v$, $E_{1}$ is Morse, and from \cite{I2}, we have that for every critical point $\phi$ of $E_{1}$, the functional $E_{2}(\phi,\cdot)$ is Morse for a generic set of perturbations $H\in \mathbb{H}^{3}_{p+1}$. Since the critical points of $E_{1}$ are countable, then we have the existence of a generic subset of $\mathbb{H}^{3}_{p+1}$ such that $E_{2}(\phi,\dot)$ is Morse for every $\phi\in crit(E_{1})$. Notice that we can also perturb the vector field $V$ in a slightly different way from the perturbation of the gradient flow of $\Lag_{H}$ to achieve transversality, as follow
\begin{equation}
V_{K}(x)=-(1+K)\left(\begin{array}{cc}
(-\Delta+1)^{-1}\Big(-\nabla_{s}\partial_{s}\phi-\nabla_{\phi}v(s,\phi)\Big)\\
(1+|\D|)^{-1}(\D_{\phi}-\nabla_{\psi}H(s,\phi,\psi))
\end{array}
\right),
\end{equation}
where $K \in \tilde{\mathbb{K}}_{\theta,\rho_{0}}^{2}\subset \mathbb{K}_{\theta,\rho_{0}}^{2}$  for $\theta$ chosen in a suitable way and we add the fact that $K$ is lower triangular. That is $K\in \tilde{\mathbb{K}}_{\theta,\rho_{0}}^{2}$ if and only if $K\in \mathbb{K}_{\theta,\rho_{0}}^{2}$ and 
$$K=\left(\begin{array}{cccc}
K_{11}&0\\
K_{21}& K_{22} 
\end{array}
\right).$$
This choice is made in such a way to preserve the decoupling of the first equation of the flow. We claim that with this choice, transversality can be achieved. Notice that the main issue one needs to worry about is the result of Lemma 5.8. First we start by two rest points $x_{+}=(\phi_{+},\psi_{+})$ and $x_{-}=(\phi_{-},\psi_{-})$. Notice that we have two kinds of flows for $V$.
The first kind is when $\phi_{+}=\phi_{-}$, hence $\ell(t)=(\phi_{+},\psi(t))$. In this case the stable and unstable manifold depend mainly on the $\psi$ component and transversality can be achieved via a perturbation $K$ with $K_{11}=K_{21}=0$ as in \cite{I2}.

The second kind of flows is when $\phi_{+}\not= \phi_{-}$ in this case we have always that $\nabla_{\phi} E_{1}(\ell(t))\not =0$. So we consider $w\in C^{2}(\R, \ell^{*}T\mathcal{F}^{1,1})$ with $ supp(w)\subset [a,b]$. We want to show the existence of $k\in \tilde{\mathbb{K}}_{\theta}$ such that $k(\ell)V(\ell)=w.$
First, if we set $w=\left[ \begin{array}{cc} 
w_{1}\\
w_{2}
\end{array}
\right]$.
We want to find a $k$ such that $k\left[ \begin{array}{cc} 
0\\
Y
\end{array}
\right]=\left[ \begin{array}{cc} 
0\\
\tilde{Y}
\end{array}
\right]$.
Notice that we can always, generate $w_{1}$ by taking 
\begin{align}
\tilde{k}_{11}(t)=&\frac{\langle \cdot, \nabla_{\phi}E_{1}(\ell(t))\rangle}{\|\nabla_{\phi}E_{1}(\ell(t))\|^{2}}w_{1}(t)+\frac{\langle \cdot, w_{1}(t)\rangle}{\|\nabla_{\phi}E_{1}(\ell(t))\|^{2}}\nabla_{\phi}E_{1}(\ell(t))\notag\\
&-\frac{\langle \nabla_{\phi}E_{1}(\ell(t)), w_{1}(t)\rangle\langle \cdot ,\nabla_{\phi}E_{1}(\ell(t)) \rangle}{\|\nabla_{\phi}E_{1}(\ell(t))\|^{4}}\nabla_{\phi}E_{1}(\ell(t)).
\end{align}

Now, if we take the $\tilde{w}=\left[ \begin{array}{cc} 
0\\
w_{2}
\end{array}
\right]$, then by taking $\tilde{k}^{1}$ as in Lemma 5.8, that is 
\begin{align}
\tilde{k}^{1}(t)=&\frac{\langle \cdot, V(\ell(t))\rangle}{\|\nabla_{1,1/2}V(\ell(t))\|^{2}}\tilde{w}(t)+\frac{\langle \cdot, \tilde{w}(t)\rangle}{\|V(\ell(t))\|^{2}}V(\ell(t))\notag\\
&-\frac{\langle V(\ell(t)), \tilde{w}(t)\rangle\langle \cdot ,V(\ell(t)) \rangle}{\|V(\ell(t))\|^{4}}V(\ell(t)),
\end{align}
we have that $\tilde{k}^{1}V=\tilde{w}$ but $$\tilde{k}^{1}\left[ \begin{array}{cc} 
0\\
Y
\end{array}
\right]=\left[ \begin{array}{cc} 
\tilde{X}\\
\tilde{Y}
\end{array}
\right]$$
But we can always compensate the component on $\tilde{X}$ using the first process and adding a term $\tilde{k}_{11}$. Now the rest of the proof stays unchanged in order to achieve transversality. From now on we will write just $V$ instead of $V_{K}$ since these properties are independent of the perturbation $K$. We consider then, the flow lines of $V$. That is solutions of 
\begin{equation}\label{equav}
\left\{\begin{array}{ll}
\frac{d\ell}{dt}=V(\ell(t))\\
\ell(0)=\ell_{0}\in \mathcal{F}^{1,1/2}
\end{array}
\right.
\end{equation}

\begin{lemma}
The solutions of the flow of $V$, $(\ref{equav})$, exist for all time.
\end{lemma}
{\it Proof:}
First, since $V$ is a vector field tangent to $\mathcal{F}^{1,1/2}(S^{1},N)$, the local existence is guarantied via the classical Cauchy-Lipschitz for ODEs. We will assume that $\ell$, the solution to $(\ref{equav})$, exists in $[0,T)$ where $T>0$ is the maximal time of existence and assume for the sake of contradiction that $T<+\infty$.  It is important to notice that $\phi(t)$, on the other hand, exists for all time and bounded in $C^{2,\alpha}$. Moreover, it converges to a critical point of the functional $E_{1}$ and this makes many estimates involving $\phi$ easier. Therefore the main difficulty resides in the behavior of the component $\psi$. Notice that if $\|\psi\|_{L^{p+1}}$ is bounded then the same estimates in the proof of Propositions 7.1 and 7.2 hold and hence we have a uniform bound on $\|\psi\|_{C^{1,\alpha}}$ and the flow can be extended beyond $T$.\\
 We claim that $\|\psi(t)\|_{p+1}$ is bounded on $[0,T)$. Indeed, the main ingredient here is the inequality $a^{2}\leq C_{\epsilon}+\varepsilon a^{p+1}$ for all $a>0$. Therefore
\begin{align}
\frac{d}{dt}E_{2}(\ell(t))&=-\|\partial_{t}\psi(t)\|_{H^{1/2}}^{2}+\langle \partial_{t}\phi(t), \frac{1}{2}R(\phi)\langle \psi, \partial_{s}\phi \cdot \psi \rangle \rangle + \langle \partial_{t}\phi(t),\nabla_{\phi}H(\psi,\psi)\rangle \notag \\
&\leq -\|\partial_{t}\psi(t)\|_{H^{1/2}}^{2} +\varepsilon C_{5}\|\psi(t)\|^{p+1}_{L^{p+1}}+C_{2},
\end{align}
where $C_{5}>0$ is a constant depending on $\sup_{t\in [0,T)}\|\phi(t)\|_{C^{2,\alpha}}$.We are back to a similar situation of the proof of Propositions 7.1 and 7.2. In fact the estimates are even easier because of the boundedness of the $\phi$ term. Thus $\|\psi(t)\|_{p+1}^{p+1}$ is bounded and $T=+\infty$.
\hfill$\Box$

\begin{proposition}
Let $\ell$ be a flow line of $V$ such that $\ell(-\infty)=x_{-}$ and $\ell(+\infty)=x_{+}$, where $x_{+},x_{-}\in rest(V)$, then then exists a constant $C$ depending on $\Lag_{H+v}(x_{+})$, $\Lag_{H+v}(x_{-})$ and the usual constants $C_{i}$, $1\leq i \leq 4$, such that 
\begin{equation}
\sup_{t\in \R} |\Lag_{H+v}(\ell(t))|\leq C.
\end{equation}
\end{proposition}
{\it Proof:}
Again, here we only need to investigate the boundedness of $E_{2}$ along the flow. This is very similar to  the invariance of the homology of the functional $E_{2}$ under homotopy that was investigated in \cite{I2}. But the main difference here is that the perturbation in \cite{I2} is compactly supported. So if for instance $\phi(t)$ was constant outside an interval $[-T,T]$ then the proof works exactly the same. But here we need to take care of the part outside this compact set.
\begin{lemma}
There exists $C=C(x_{+},x_{-})$ and $T=T(x_{+},x_{-})>0$ such that
$$E_{2}(\ell(-t))\leq C(x_{+},x_{-})$$ and $$-C(x_{+},x_{-})\leq E_{2}(\ell(t)),$$
for all $t\geq T$.
\end{lemma}
{\it Proof:}
We recall that
\begin{align}
\frac{d}{dt}E_{2}(\ell(t))&=-\|\partial_{t}\psi(t)\|_{H^{1/2}}^{2}+\langle \partial_{t}\phi(t), \frac{1}{2}R(\phi)\langle \psi, \partial_{s}\phi \cdot \psi \rangle \rangle + \langle \partial_{t}\phi(t),\nabla_{\phi}H(\psi,\psi)\rangle \notag \\
&\leq -\|\partial_{t}\psi(t)\|_{H^{1/2}}^{2} + C_{5}\|\partial_{t}\phi(t)\|_{H^{1}}\Big(\|\psi(t)\|^{p+1}_{L^{p+1}}+C_{2}\Big).\notag
\end{align}
It is a well known result now, since $E_{1}$ is Morse, we have the existence of $\alpha>0$ such that
$$\|\partial_{t}\phi(t)\|_{H^{1}}\leq Ce^{-\alpha|t|}.$$
We take $T>>1$ and we investigate the energy in the interval $(-\infty,-T]$. We have then

$$\frac{d}{dt}E_{2}(\ell(t))\leq -\|\partial_{t}\psi(t)\|_{H^{1/2}}^{2} + Ce^{-\alpha|t|}\Big(\|\psi(t)\|^{p+1}_{L^{p+1}}+1\Big).$$
Hence,
\begin{equation}\label{e2}
E_{2}(\ell(t))\leq E_{2}(x_{+})-\int_{-\infty}^{t}\|\partial_{t}\psi(s)\|_{H^{1/2}}^{2}ds+C\int_{-\infty}^{t}e^{-\alpha|s|}(\|\psi(s)\|^{p+1}_{L^{p+1}}+1)ds.
\end{equation}
But
$$\langle \partial_{t}\psi, \psi \rangle +2E_{2}(\ell(t))\geq C_{2}\|\psi(t)\|^{p+1}_{p+1}-C_{2}.$$
Hence, if we set $A=E_{2}(x_{+})$, we have
$$\|\psi(t)\|^{p+1}_{p+1}\leq C+2A+C\int_{-\infty}^{-T}e^{-\alpha|s|}\|\psi(s)\|^{p+1}_{p+1} ds +C\|\partial_{t}\psi\|_{H^{1/2}}\|\psi\|_{H^{1/2}}.$$
Therefore, if we multiply by $e^{-\alpha |t|}$ and integrate, we have that for $T$ big enough,
$$\int_{-\infty}^{-T}e^{-\alpha|t|}\|\psi(t)\|^{p+1}_{p+1}dt \leq Ce^{-\alpha T}(1+A)+\int_{-\infty}^{-T}e^{-\alpha |t| }\|\partial_{t}\psi\|_{H^{1/2}}\|\psi\|_{H^{1/2}}dt$$
and
$$
\|\psi(t)\|^{p+1}_{p+1}\leq C(A)+C\int_{-\infty}^{-T}e^{-\alpha |t| }\|\partial_{t}\psi\|_{H^{1/2}}\|\psi\|_{H^{1/2}}dt +C\|\partial_{t}\psi\|_{H^{1/2}}\|\psi\|_{H^{1/2}}.$$
Next, we estimate $\|\psi\|_{H^{1/2}}^{2}$. Similarly to (7.13-7.14), we have that
$$\|\psi\|^{2}_{H^{1/2}}\leq C\|\psi\|_{H^{1/2}}+C(A)+C\int_{-\infty}^{-T}e^{-\alpha |t|}\|\partial_{t}\psi\|_{H^{1/2}}\|\psi\|_{H^{1/2}}dt+C\|\partial_{t}\psi\|_{H^{1/2}}\|\psi\|_{H^{1/2}}.$$
Therefore, if one multiplies by the exponential term and integrates, he finds
\begin{align}
\int_{-\infty}^{-T}e^{-\alpha|t|}\|\psi\|^{2}_{H^{1/2}}dt &\leq C(A)+C\int_{-\infty}^{-T}e^{-\alpha |t| }\|\psi\|_{H^{1/2}}dt +Ce^{-\alpha T}\int_{-\infty}^{-T}e^{-\alpha |t| }\|\partial_{t}\psi\|_{H^{1/2}}\|\psi\|_{H^{1/2}}dt\notag\\
&\quad+C\int_{-\infty}^{-T}e^{-\alpha |t| }\|\partial_{t}\psi\|_{H^{1/2}}\|\psi\|_{H^{1/2}}dt\notag \\
&\leq C(A)+Ce^{-\alpha T}\left(\int_{-\infty}^{-T}e^{-\alpha|t|}\|\psi\|_{H^{1/2}}^{2}dt\right)^{\frac{1}{2}}\notag \\
&\quad+ C\left(\int_{-\infty}^{-T}e^{-\alpha |t| }\|\partial_{t}\psi\|_{H^{1/2}}^{2}dt\right)^{\frac{1}{2}}\left(\int_{-\infty}^{-T}e^{-\alpha|t|}\|\psi\|_{H^{1/2}}^{2}dt\right)^{\frac{1}{2}}.\notag
\end{align}
Again, from $(\ref{e2})$, we have that
\begin{align}
\int_{-\infty}^{-T}\|\partial_{t}\psi(t)\|_{H^{1/2}}^{2}dt&\leq C(A)-E_{2}(\ell(-T))+C\int_{-\infty}^{-T}e^{-\alpha|s|}\|\psi(s)\|^{p+1}_{p+1}ds\notag \\
&\leq C(A)-E_{2}(\ell(-T))+C\int_{-\infty}^{-T}e^{-\alpha |t| }\|\partial_{t}\psi\|_{H^{1/2}}\|\psi\|_{H^{1/2}}dt\notag \\
&\leq C(A)-E_{2}(\ell(-T))+C\Big(\int_{-\infty}^{-T}e^{-\frac{\alpha}{2} |t| }\|\partial_{t}\psi\|_{H^{1/2}}^{2}dt+\int_{-\infty}^{-T}e^{-\frac{3\alpha}{2} |t| }\|\psi\|_{H^{1/2}}^{2}dt\Big).\notag
\end{align}
Taking $T$ even bigger, we have
$$\int_{-\infty}^{-T}\|\partial_{t}\psi(t)\|_{H^{1/2}}^{2}dt\leq C(A)-4E_{2}(\ell(-T))+C\int_{-\infty}^{-T}e^{-\frac{3\alpha}{2} |t| }\|\psi\|^{2}dt.$$
Using this estimate, we have
\begin{align}
\int_{-\infty}^{-T}e^{-\alpha|t|}\|\psi\|_{H^{1/2}}^{2}dt &\leq C(A)+Ce^{-\alpha T}\left(\int_{-\infty}^{-T}e^{-\alpha|t|}\|\psi\|_{H^{1/2}}^{2}dt\right)^{\frac{1}{2}}+ C\left(\int_{-\infty}^{-T}e^{-\alpha |t| }\|\partial_{t}\psi\|_{H^{1/2}}^{2}dt\right.\notag\\
&\quad\left.+\frac{1}{4}\int_{-\infty}^{-T}e^{-\alpha|t|}\|\psi\|_{H^{1/2}}^{2}dt\right).\notag
\end{align}
Thus,
\begin{align}
\int_{-\infty}^{-T}e^{-\alpha|t|}\|\psi\|_{H^{1/2}}^{2}dt &\leq C(A)+Ce^{-\alpha T}\left(\int_{-\infty}^{-T}e^{-\alpha|t|}\|\psi\|_{H^{1/2}}^{2}dt\right)^{\frac{1}{2}}+ C\int_{-\infty}^{-T}e^{-\alpha |t| }\|\partial_{t}\psi\|_{H^{1/2}}^{2}dt\notag \\
&\leq C(A)+Ce^{-\alpha T}\left(\int_{-\infty}^{-T}e^{-\alpha|t|}\|\psi\|_{H^{1/2}}^{2}dt\right)^{\frac{1}{2}}+Ce^{-\alpha T}\biggl(C(A)-4E_{2}(\ell(-T))\notag\\
&\quad +C\int_{-\infty}^{-T}e^{-\frac{3\alpha}{2} |t| }\|\psi\|^{2}_{H^{1/2}}dt\biggr).\notag
\end{align}
Therefore,
$$\int_{-\infty}^{-T}e^{-\alpha|t|}\|\psi\|^{2}_{H^{1/2}}dt\leq C(A)-CE_{2}(\ell(-T))$$
The rest of the proof now follows exactly Proposition 7.1 to yield $$E_{2}(\ell(-T))\leq C(E_{2}(\ell(-\infty))).$$
A similar inequality holds when we start the integrations in $[T,+\infty)$.
\hfill$\Box$

For the bound on $[-T,T]$ we follow exactly the procedure done in Proposition 7.1, since
$$\frac{d}{dt}E_{2}(\ell(t))\leq-\|\partial_{t}\psi(t)\|_{H^{1/2}}^{2}+ C_{5}\varepsilon\|\psi(t)\|^{p+1}_{L^{p+1}}+C$$
\hfill$\Box$

In order to achieve compactness, as in the previous case, we need to show that the length of the orbits is uniformly bounded.
\begin{lemma}
Let $\ell$ be a flow line of $V$ such that $\ell(-\infty)=x_{-}$ and $\ell(+\infty)=x_{+}$, where $x_{+},x_{-}\in rest(V)$, then then exists a constant $C$ depending on $\Lag_{H+v}(x_{+})$, $\Lag_{H+v}(x_{-})$ such that the length of $\ell$, $L(\ell)$ satisfies
$$L(\ell)\leq C.$$
\end{lemma}
{\it Proof:}
We first bound the ends of the flow at infinity. We will use the fact that $V$ is asymptotically hyperbolic to get exponential convergence of the flow near the critical points and hence integrability.
We consider the function defined for $t>T>0$ by $$f(t)=\frac{1}{2}\|\ell'(t)\|^{2}.$$
If we compute $f''(t),$ we find
$$f''(t)=\|\ell''(t)\|^{2}+\langle \ell'(t),dV(\ell(t))\ell''(t)\rangle +\langle \ell'(t),d^{2}V(\ell(t))[\ell'(t),\ell'(t)]\rangle.$$
But in a local frame around $\ell(+\infty)=x_{-}$, we can write that $dV(\ell(t))=dV(x_{-})+g(t)$  where $g(t)\to 0$ as $t\to \infty$. Hence, using the symmetry of $dV(\ell(t)),$ we have that
$$\langle \ell'(t),dV(\ell(t))\ell''(t)\rangle = \|dV(x_{-})\ell'(t)\|^{2}+\epsilon(t)\|\ell'(t)\|^{2}.$$
Also, one has $$\langle \ell'(t),d^{2}V(\ell(t))[\ell'(t),\ell'(t)]\rangle\geq -c\|\ell'(t)\|^{3}.$$
Hence, using the invertibility of $dV(x_{-})$, one has
$$f''(t)\geq \|\ell'(t)\|^{2}(c_{1}-\epsilon(t)-c\|\ell'(t)\|).$$
Hence by taking the neighborhood even smaller we have that
$$f''(t)\geq \alpha f(t).$$
Now based on the proof of the previous Proposition, we have furthermore that $f\in L^{1}([T,+\infty))$ and it is a classical argument to see that 
$$\|\ell'(t)\|\leq Ce^{-\delta |t|},$$
for $\delta$ depending only on the end points. Notice here that in the interval $[-T,T]$ we used the boundedness of $\int_{-T}^{T} \|\ell'(t)\|^{2}dt$ as shown in the previous Proposition. This shows now that the length $L$ is uniformly bounded.
\hfill$\Box$

In particular, this previous Lemma  and Proposition  tell us that the moduli space $\mathcal{M}(x_{-},x_{+},V)$ is a relatively compact manifold of dimension $dim \mathcal{M}(x_{-},x_{+},V)=\mu_{0}(x_{-})-\mu_{0}(x_{+})$. Therefore, by considering critical points with difference of index $1$, we can perturb the metric in a generic way to have transversality up to order 2 and use the gluing argument to construct the homology $HF_{*}(V,\mathcal{F}^{1,1/2},\Z_{2})$.

\begin{proposition}
The homology $HF_{*}(V,\mathcal{F}^{1,1/2},\Z_{2})$ is equal to the homology $DG^{p+1}HF_{*}(\mathcal{F}^{1,1/2},\Z_{2})$.
\end{proposition}
{\it Proof:}
The idea here, is similar to the part related to the invariance of the homology with respect to the perturbation of $H\in \mathbb{H}_{p+1}^{3}$. Indeed, we construct a homomorphism $\Phi_{V,\Lag}:HF_{*}(V,\mathcal{F}^{1,1/2},\Z_{2})\to HF_{*}(\Lag_{H+v},\mathcal{F}^{1,1/2},\Z_{2})$ by taking a cut-off function $\eta:\R\to [0,1]$ such that $\eta(t)=0$ for $t\leq 0$ and $\eta(t)=1$ for $t\geq 1$ and consider the vector field $$V_{t}=\eta(t)(-\nabla_{G^{K}}\Lag_{H+v})+(1-\eta(t))V.$$
Then we consider the flow lines of $V_{t}$, that are solutions to 
$$\frac{d\ell}{dt}=V_{t}(\ell(t)), \text{ } \ell(-\infty)=x_{-} \text{ and }\ell(+\infty)=x_{+}.$$
The compactification of the moduli space $\mathcal{M}(x_{-},x_{+},V_{t})$ works exactly the same as in Section 7 and by counting the flow lines for moduli spaces of dimension zero, we get the desired result.
\hfill$\Box$

\subsection{Leray-Serre type Floer spectral sequence}

We will construct now a spectral sequence converging to the homology generated by $V$ and this will lead us to the computation of our homology. 
We first start by defining the chain complex generated by the rest points of $V$, that is $(C_{*}(V,\Z_{2}),\partial_{*}^{V})$. Then one can define a filtration $(F_{r}(C_{*}(V,\Z_{2}))_{r\geq 0}$ as follows :
$$F_{r}C_{q}(V,\Z_{2})=span_{\Z_{2}}\{x=(\phi,\psi)\in rest_{q}(V);\mu_{1}(\phi)\leq r\}.$$
Clearly 
\begin{itemize}
\item[i)] $F_{r}C_{*}(V,\Z_{2})=\{0\}$ for $r<0$,
\item[ii)] $F_{r}C_{*}(V,\Z_{2})\subset F_{r+1}C_{*}(V,\Z_{2})$,
\item[iii)] $\bigcup_{r\geq0} F_{r}C_{*}(V,\Z_{2}) = C_{*}(V,\Z_{2})$.
\end{itemize}
 Moreover, since the flow of $V$ projects to a negative gradient flow of $E_{1}$, then the boundary operator $\partial^{V}$ decreases the index $\mu_{1}$. Therefore $\partial^{V}(F_{r}C_{*}(V,\Z_{2}))\subset F_{r}C_{*}(V,\Z_{2})$. So the filtration we have is in fact a chain filtration. Thus, we consider its associate spectral sequence.
It follows then from $i)-iii)$ and \cite[Chap 9, Theorem 2]{SP}, that there exists a spectral sequence, with first page 
\begin{equation}\label{p1}
E^{1}_{r,q}=H_{r+q}(F_{r}C(V,\Z_{2})/F_{r-1}C(V,\Z_{2})).
\end{equation}
With boundary operator $\partial_{1}$ corresponding to the triple $(F_{r}C(V,Z_{2}),F_{r-1}C(V,Z_{2}),F_{r-2}C(V,Z_{2}))$. That is the boundary operator coming from the short exact sequence
$$0\to F_{r-1}C(V,Z_{2}) \to F_{r}C(V,Z_{2}) \to F_{r}C(V,Z_{2})/F_{r-1}C(V,Z_{2})\to 0.$$

\begin{lemma}
$$E^{1}_{r,q}=\bigoplus_{\phi \in crit_{r}(E_{1})} D^{p+1}HF_{q}(H^{1/2}(S^{1}),\Z_{2}).$$
\end{lemma}
{\it Proof:}
We first notice that  $$F_{r}C(V,\Z_{2})/F_{r-1}C(V,\Z_{2})=span_{\Z_{2}}\{(\phi,\psi)\in rest(V);  \phi \in crit_{r}(E_{1}); \psi \in crit(E_{2}(\phi, \cdot))\}$$
with the boundary operator $$\partial : F_{r}C_{r+q}(V,\Z_{2})/F_{r-1}C_{r+q}(V,\Z_{2})\to F_{r}C_{r+q-1}(V,\Z_{2})/F_{r-1}C_{r+q-1}.$$
 Now consider two rest points $x_{-}=(\phi_{-},\psi_{-}),x_{+}=(\phi_{+},\psi_{+})\in \{(\phi,\psi)\in rest(V);  \phi \in crit_{r}(E_{1}); \psi \in crit(E_{2}(\phi, \cdot))\}$, and let $\ell=(\phi(t),\psi(t)) \in \mathcal{M}(x_{-},x_{+},V)$, $\ell$ projects down to a flow line of the negative gradient flow of $E_{1}$. Since $\mu_{1}(\phi_{-})=\mu_{1}(\phi_{+})=r$, we have by transversality, that the path $\phi(t)$ is constant and that $\phi(t)=\phi_{-}=\phi_{+}$. This implies that $\psi_{-}$ and $\psi_{+}$ are both critical points of $E_{2}$ restricted to the same fiber, that is $\psi_{-}, \psi_{+}\in crit(E_{2}(\phi,\cdot))$ and the component $\psi(t)$ of $\ell$ is then a flow line of the negative gradient vector of $E_{2}(\phi,\cdot)$. It follows then that 
$$\mathcal{M}(x_{-},x_{+},V)\cong \mathcal{M}(\psi_{-},\psi_{+},E_{2}(\phi,\cdot)).$$
So for every fixed critical point $\phi\in crit_{r}(E_{1})$, we obtain the homology $HF_{*}(E_{2}(\phi,\cdot), H^{1/2}(S^{1}),\Z_{2})$. This homology is independent of the fiber, indeed, if $\phi_{1}$ and $\phi_{2}$ are critical points of $E_{1}$ of index $r$, then $$HF_{*}(E_{2}(\phi_{1},\cdot), H^{1/2}(S^{1}),\Z_{2})=HF_{*}(E_{2}(\phi_{2},\cdot), H^{1/2}(S^{1}),\Z_{2})=D^{p+1}HF_{*}(H^{1/2}(S^{1}),\Z_{2}).$$
This follows from the result in $\cite{I2}$ of the invariance of the homology with respect to functions $H\in \mathbb{H}^{3}_{p+1}$. This can be seen by taking a geodesic path $\gamma(s,t)_{0\leq t \leq 1} $ connecting $\phi_{1}(s)$ to $\phi_{2}(s)$. Then using the parallel transport we pull back the problem to one single fiber above $\phi_{1}$ and thus we will have an equivalent problem with a homotopy between two perturbations in $\mathbb{H}^{3}_{p+1}$, hence they have the same homology as shown in $\cite{I2}$  and this finishes the proof.
\hfill$\Box$

With all these ingredient now we can state the third claim in Theorem 1.

\begin{corollary}
The homology $DG^{p+1}H_{*}(\mathcal{F}^{1,1/2}(S^{1},N),\Z_{2})$ vanishes.
\end{corollary}

{\it Proof:}
So far we constructed a first quadrant spectral sequence, with first page $E^{1}_{r,q}$ defined in $(\ref{p1})$. Again, from \cite{SP}, we have that this spectral sequence converges to $GH_{*}(C(V,\Z_{2}))$. But from $\cite{I2}$, $D^{p+1}HF_{*}(H^{1/2}(S^{1}),\Z_{2})=0$. Therefore, by Lemma 8.4, the spectral sequence collapses at the first page and $$GH_{*}(C(V,\Z_{2}))=0.$$
Therefore one has $$HF_{*}(V,\mathcal{F}^{1,1/2}(S^{1},N),\Z_{2})=0.$$
Then by Proposition 8.2,
$$DG^{p+1}HF_{*}(\mathcal{F}^{1,1/2}(S^{1},N),\Z_{2})=HF_{*}(V,\mathcal{F}^{1,1/2}(S^{1},N),\Z_{2})=0.$$

\hfill$\Box$

\section{Applications to the existence of superquadratic Dirac-geodesics}
In this section we assume that $\nabla_{\psi}H(\phi,0)=\nabla_{\psi,\psi}H(\phi,0)=\nabla_{\phi,\psi}H(\phi,0)=0$. We will suppose furthermore that $H$ is even with respect to $\psi$. We will assume for the sake of contradiction that the only critical points of $\Lag_{H}$ are the trivial one, i.e. $(\phi,0)$ with $\phi$ a perturbed geodesic. Since $\Lag_{H}$ is Morse, the critical points are isolated. We fix a critical point $(\phi_{0},0)$ with $\phi_{0}$ a critical point of $E_{1}$ of index zero and we will compute the relative index $\mu_{H}$ starting from it.

For every critical point $(\phi_{1},0)$ of $\Lag_{H}$ we will associate two indices, that is 
$$\mu_{H}(\phi_{1},0)=\mu_{E_{1}}(\phi_{1})+\mu_{1}(\phi_{1}),$$
where $\mu_{1}(\phi)=-\mathsf{sf}\{\D_{\phi_{t}}\}_{0\leq t \leq 1}$. We will write that $(\phi_{1},0)$ is of index $(\mu_{E_{1}},\mu_{1})$.
We assume now that $0\not=[\phi_{0}]\in H_{0}(\Lambda N)$, where $\Lambda N$ is the loop space of $N$. Now since $DG^{p+1}H(\mathcal{F}^{1,1/2}(S^{1},N),\Z_{2})=0$, either $(\phi_{0},0)$ is a boundary in $C_{*}(\Lag_{H},Z_{2})$ or it maps to a boundary. We consider the first case, that is there exists $\phi_{1},\phi_{2},\cdots, \phi_{\ell}\in crit(E_{1})$ such that $\partial\sum_{i=1}^{\ell} (\phi_{i},0)=(\phi_{0},0)$. We can do the same reasoning for the second case, therefore we will omit it.  So $\mu_{H}(\phi_{i},0)=1$ for $1\leq i \leq \ell$ and we will set $\mu_{E_{1}}(\phi_{1})=m_{i}+1$ and $\mu_{1}(\phi_{i})=-m_{i}$. We distinguish three cases:

\textbf{Case 1:} If $m_{i}>0$, then we have already that $dim \hat{\mathcal{M}}(\phi_{i},\phi_{0},E_{1})=m_{i}\geq 1$. But it is easy to see that a gradient flow line of $E_{1}$ corresponds to a flow line of $\Lag_{H}$ in the neighborhood of the critical points (since the perturbation of the metric occurs outside neighborhoods of critical points), via the map $\phi\mapsto (\phi,0)$. Therefore $\mathcal{M}((\phi_{i},0),(\phi_{0},0),H)$ has to be empty or else the dimensions will not match since in that case one would have $dim\hat{\mathcal{M}}((\phi_{i},0),(\phi_{0},0),H)=0<m_{i}$. Therefore, $(\phi_{0},0)$ cannot be canceled by a critical point with negative $\mu_{1}$-index.

\textbf{Case 2:} If $m_{i}=0$ for all $1\leq i \leq \ell$, then the following holds

\begin{proposition}
Let $(\phi_{+},0)$ and $(\phi_{-},0)$ be two critical points of $\Lag_{H}$, with $\mu_{1}(\phi_{+})=\mu_{1}(\phi_{-})$ then
$$\mathcal{M}((\phi_{+},0),(\phi_{-},0),H)\cong \mathcal{M}(\phi_{+},\phi_{-},E_{1})$$
\end{proposition}
\textit{Proof:}
Assume that there is a flow line of the form $x(t)=(\phi(t),\psi(t))\in\mathcal{M}((\phi_{+},0),(\phi_{-},0),H)$ with $\psi(t)\ne 0$. Under our assumption, we note that $\nabla_{\psi}\Lag_{H}(\psi(t))=-\dot{\psi}(t)\not\equiv 0$. Thus, as in the proof of Proposition 5.2 (in particular Lemma 5.8), further perturbing the metric with perturbation of the form $K=\begin{pmatrix}0&0\\0&K_{22}\end{pmatrix}$ if necessary, we may assume that the $\psi$-part of the flow $\dot{\psi}=-\nabla_{\psi}\Lag_{H}$ is regular along $x(t)$, i.e., its linearization with respect to $\psi$-variable is onto. Under this regularity condition, by linearizing the negative gradient flow equations at both ends $(\phi_{\pm},0)$, we see that along the path $x(t)$, there exists a path of eigenvalues $\{\lambda(t)\}$ of $d^{2}_{\psi}\Lag_{H}(x(t))$ such that $\lambda(-\infty)<0$ and $\lambda(+\infty)>0$ and there is no path of eigenvalues $\{\mu(t)\}$ of $d^{2}_{\psi}\Lag_{H}(x(t))$ with $\mu(-\infty)>0$ and $\mu(+\infty)<0$. Thus the spectral flow $\mathsf{sf}\{d^{2}_{\psi}\Lag_{H}(x(t))\}_{-\infty\le t\le+\infty}$ is positive under the assumption. But this contradicts the equality of the $\mu_{1}$-index.
Thus, there does not exist flow lines of the form $(\phi(t),\psi(t))\in\mathcal{M}((\phi_{+},0),(\phi_{-},0),H)$ with $\psi(t)\ne 0$. This proves that the natural injection $\mathcal{M}(\phi_{+},\phi_{-},E_{1})\ni\phi(t)\mapsto(\phi(t),0)\in\mathcal{M}((\phi_{+},0),(\phi_{-},0),H)$ is also a surjection. This completes the proof.
\hfill$\Box$

By applying this Proposition to each critical points $(\phi_{i},0)$ and $(\phi_{0},0)$ we see that since $(\phi_{0},0)$ is a boundary of $\sum_{i=1}^{\ell}(\phi_{i},0)$, then $\phi_{0}$ is a boundary of $\sum_{i=1}^{\ell}\phi_{i}$ in $C_{*}(E_{1},\Z_{2})$ therefore, $[\phi_{0}]=0$ which is a contradiction.

\textbf{Case 3:} If $m_{i}=-1$ for $1\leq k \leq \ell$ and $m_{i}=0$ for $k<1 \leq \ell$.

\begin{theorem}

If $H$ satisfies  $\nabla_{\psi}H(\phi,0)=\nabla_{\psi,\psi}H(\phi,0)=\nabla_{\phi,\psi}H(\phi,0)=0$, then either we have at least three (perturbed) geodesics, two minimals and one geodesic of index $1$ or we have a non-trivial Dirac Geodesic.
\end{theorem} 
{\it Proof:}
We focus on case 3, and assume for the sake of generality that $\ell=1$ and $m_{1}=-1$. That means that $\phi_{0}$ and $\phi_{1}$ are two geodesics of index zero, in the same connected component of $\Lambda N$. Thus, by connectedness, we need to have another (perturbed) geodesic of index $1$. Now, if $\ell>1$ then clearly we have at least three (perturbed) geodesics as described in the theorem.
\hfill$\Box$

If we assume furthermore that $H$ is even in the $\psi$ component, then we have 
\begin{corollary}
Assume that $H\in \mathbb{H}^{3}_{p+1}$ is even in $\psi$ and satisfies  $\nabla_{\psi}H(\phi,0)=\nabla_{\psi,\psi}H(\phi,0)=\nabla_{\phi,\psi}H(\phi,0)=0$, then we have the existence of at least one non-trivial perturbed Dirac Geodesic.
Moreover, if there exists a sequence of numbers $a_{k}\to \infty$ such that $b_{a_{k}}(\Lambda N,\Z_{2})\not=0$ then we have infinitely many solutions.
\end{corollary}
{\it Proof:}
We see from case 3, that if $m_{i}=-1$ then generically, there is no $E_{1}$-flow lines between $\phi_{i}$ and $\phi_{0}$ since $\mu_{E_{1}}(\phi_{1})=\mu_{E_{1}}(\phi_{0})=0$. Thus the only flow lines we have in $\mathcal{M}((\phi_{i},0),(\phi_{0},0),H)$ are of the form $(\phi(t),\psi(t))$ where $\psi(t)\not=0$. But near the critical points, since the perturbation of the metric vanishes, we have that if $(\phi(t),\psi(t))$ is a flow line then $(\phi(t),-\psi(t))$ is also a flow line. Thus, the flow lines come in pair and $\sharp \hat{\mathcal{M}}((\phi_{i},0),(\phi_{0},0),H)=0(\text{ mod 2})$ leading to another contradiction.
Notice now that this process involves (perturbed) geodesics of difference of index at most one. Hence, we can repeat the same reasoning for another generator $[\phi_{0,a_{k}}]\in H_{a_{k}}(\Lambda N)$ and since $a_{k}\to \infty$ we have indeed the existence of infinitely many non-trivial Dirac-geodesics. 
\hfill$\Box$

In the previous results the word perturbed can be removed to deal with actual geodesics in case $g$, the metric on $N$ is bumpy, which is a generic condition on the metric.

\medskip

Finally, we consider a spacial case of flat tori $N=\T$. As before, we argue by contradiction and assume that there are no non-trivial Dirac-geodesics. In this case, $\D_{\phi}=\D$ and $\mu_{1}(\phi)=-\mathsf{sf}\{\D_{\phi_{t}}\}_{0\le t\le1}=0$. Then by Proposition 9.1, the canonical embedding $\mathcal{M}(\phi_{+},\phi_{-},E_{1})\subset\mathcal{M}((\phi_{+},0),(\phi_{-},0),H)$ is an isomorphism:
$$\mathcal{M}(\phi_{+},\phi_{-},E_{1})\cong \mathcal{M}((\phi_{+},0),(\phi_{-},0),H).$$
Thus, the homology $DG^{p+1}H_{\ast}(\mathcal{F}^{1,1/2}(S^{1},N);\Z_{2})$ is isomorphic to the Morse to the Morse homology $MH_{\ast}(\Lambda N,E_{1};\Z_{2})$ of the perturbed geodesic functional. The latter is isomorphic to the singular homology of the loop space $H_{\ast}(\Lambda N;\Z_{2})$. Thus we arrive at a contradiction $0=DG^{p+1}H_{\ast}(\mathcal{F}^{1,1/2}(S^{1},N);\Z_{2})\cong MH_{\ast}(\Lambda N,E_{1};\Z_{2})\cong H_{\ast}(\Lambda N;\Z_{2})\ne 0$. Thus, we have

\begin{theorem}
Let $N=\T$ be a flat tori. Under the assumption $\nabla_{\psi}H(\phi,0)=\nabla_{\psi,\psi}H(\phi,0)=\nabla_{\phi,\psi}H(\phi,0)=0$, there exists a non-trivial perturbed Dirac geodesics in each connected component.
\end{theorem}

{\bf Concluding Remarks:}\\

We point out that in Corollary 9.1, the condition on the Betti numbers is satisfied in many situations. For instance, in \cite{Serre}, Serre proved that if $H_{*}(N,\Z_{2})\not = 0$ for any nonzero degree, then $H_{i}(\Lambda M, \Z_{2})$ is non-trivial for infinitely many $i$. Then extending the work of  D. Sullivan \cite{Sul} and M. Vigu\'{e}-Poirrier, D. Sullivan \cite{Sul2}, J. McCleary \cite{Mc}, proved that if $N$ is simply connected and $H^{*}(N,\Z_{2})$ as an algebra, is generated by at least two generators, then $b_{k}(\Lambda N,\Z_{2})$ grow unbounded. So in these cases, we have that the result of Corollary 9.1 holds.

Now if we look back at Theorem 9.1, we can see that we are not using Theorem 1.1 in full strength since we are just using one non-trivial generator $[\phi_{0}]$ in the homology of the loop space. But one can investigate further the different generators. The non-existence of non-trivial Dirac geodesics would impose strong restriction on the loop space. Another angle that was not used here was the $S^{1}$ action. One could investigate the index growth after iteration  and other existence results can be extracted from Theorem 1.1.

\section{Appendix}

In this section, we collect some technical results which were used in the proofs of previous sections.

\begin{lemma}[{\cite[Lemma 6.1]{BM}}]
Let $1<p<\infty$, $0<s<\infty$, $1<r<\infty$, $0<\theta<1$, $1<t<\infty$ be such that $\frac{1}{r}+\frac{\theta}{t}=\frac{1}{p}$. For $f\in W^{s,t}(S^{1})\cap L^{\infty}(S^{1})$, $g\in W^{\theta s,p}(S^{1})\cap L^{r}(S^{1})$, we have $fg\in W^{\theta s,p}(S^{1})$ and
\begin{equation}
\|fg\|_{W^{\theta s,p}(S^{1})}\le C(\|f\|_{L^{\infty}(S^{1})}\|g\|_{W^{\theta s,p}(S^{1})}+\|g\|_{L^{r}(S^{1})}\|f\|_{W^{s,t}(S^{1})}^{\theta}\|f\|_{L^{\infty}(S^{1})}^{1-\theta}).
\end{equation}
\end{lemma}
For the proof, see~\cite{BM}. Note that the conclusion of the lemma holds without the dimension restriction. However, we state it only for the $1$-dimensional case since we only use it for this special case.

In particular, for $p=2$, $s=1$, $r=4$, $\theta=\frac{1}{2}$, $t=2$, Lemma 10.1 implies the following: For $f\in W^{1,2}(S^{1})$ and $g\in H^{1/2}(S^{1})$, we have $fg\in H^{1/2}(S^{1})$ and
\begin{equation}
\|fg\|_{H^{1/2}(S^{1})}\le C(\|f\|_{L^{\infty}(S^{1})}\|g\|_{H^{1/2}(S^{1})}+\|g\|_{L^{4}(S^{1})}\|f\|_{H^{1}(S^{1})}^{\frac{1}{2}}\|f\|_{L^{\infty}(S^{1})}^{\frac{1}{2}}.
\end{equation}
Note that, in $1$-dimension, we have $H^{1}(S^{1})\subset L^{\infty}(S^{1})$ and $H^{1/2}(S^{1})\subset L^{4}(S^{1})$ by the Sobolev embedding theorem.

We also have the following corollary:

\begin{corollary}
For any $\frac{1}{2}<s<1$, we have the continuous multiplication:
\begin{equation}
H^{-1/2}(S^{1})\times H^{1/2}(S^{1})\ni(g,h)\mapsto gh\in H^{-s}(S^{1}).
\end{equation}
\end{corollary}
{\it Proof.} For $f\in H^{-1/2}(S^{1})$, $g\in H^{1/2}(S^{1})$ and $h\in C^{\infty}(S^{1})$, the distribution $fg$ is defined as
$$\langle fg,h\rangle=\langle f,gh\rangle.$$
Thus, to prove the assertion, we need to prove $gh\in H^{1/2} (S^{1})$ for $g\in H^{1/2}(S^{1})$ and $h\in H^{s}(S^{1})$. This follows form Lemma 10.1 by taking $p=2$, $s=s$, $r=\frac{4s}{2s-1}$, $\theta=\frac{1}{2s}$ and $t=2$.
Or, we can directly prove this as follows: We note that the continuity of
\begin{equation}
L^{2}(S^{1})\times H^{s}(S^{1})\ni(f,g)\mapsto fg\in L^{2}(S^{1})
\end{equation}
and
\begin{equation}
H^{s}(S^{1})\times H^{s}(S^{1})\ni(f,g)\mapsto fg\in H^{s}(S^{1}),
\end{equation}
where the continuity of (10.5) follows from $H^{s}(S^{1})\subset L^{\infty}(S^{1})$ for $s>\frac{1}{2}$ and $H^{s}(S^{1})$ is an algebra. Thus, for any fixed $h\in H^{s}(S^{1})$, the following maps are continuous:
\begin{equation}
L^{2}(S^{1})\ni g\mapsto gh\in L^{2}(S^{1}),
\end{equation}
\begin{equation}
H^{s}(S^{1})\ni g\mapsto gh\in H^{s}(S^{1}).
\end{equation}
Since $H^{1/2}(S^{1})$ is obtained as an interpolation space between $L^{2}(S^{1})$ and $H^{s}(S^{1})$, the continuity of the multiplication $H^{1/2}(S^{1})\times H^{s}(S^{1})\ni(f,g)\mapsto fg\in H^{1/2}(S^{1})$ follows by interpolating (10.6) and (10.7).
\hfill$\Box$

\begin{lemma} Let $\phi\in H^{1}(S^{1},N)$. Then $P_{\phi}$ defined by $P_{\phi}(\psi)(s)=({\bm 1}\otimes P_{\phi(s)})\psi(s)$ for $\psi\in H^{1/2}(S^{1},\Spin(S^{1})\otimes\underline{\R}^{k})$ and $s\in S^{1}$ defines a map $P_{\phi}:H^{1/2}(S^{1},\Spin(S^{1})\otimes\underline{\R}^{k})\to H^{1/2}(S^{1},\Spin(S^{1})\otimes\phi^{\ast}TN)$. Moreover, we have the following estimate:
\begin{equation}
\|P_{\phi}\psi\|_{H^{1/2}(S^{1})}\le C(\|\psi\|_{H^{1/2}(S^{1})}+\|\partial_{s}\phi\|_{L^{2}(S^{1})}^{1/2}\|\psi\|_{L^{4}(S^{1})}).
\end{equation}
\end{lemma}
{\it Proof.} Let $\{U_{\alpha}\}_{\alpha=1}^{k}$ be a covering of $N$ consisting of local coordinates such that $TN|_{U_{\alpha}}\cong U_{\alpha}\times\R^{n}$. For each $1\le\alpha\le k$, let $\{e_{\alpha,j}\}_{j=1}^{n}$ be an orthonormal frame fields of $TN|_{U_{\alpha}}$. For $s\in S^{1}$ such that $\phi(s)\in N$, $P_{\phi}$ is given as
\begin{equation}
P_{\phi}(\psi)(s)=\sum_{j=1}^{n}(\psi(s),e_{\alpha,j}(\phi(s)))_{T_{\phi(s)}N}\otimes e_{\alpha,j}(\phi(s))
\end{equation}
and the expression is independent of the choices of $\alpha$ and the frame $\{e_{\alpha,j}\}_{j=1}^{n}$.

Let $\{\rho_{\alpha}\}_{\alpha=1}^{k}$ be a partition of unity subordinate to the covering $\{U_{\alpha}\}_{\alpha=1}^{k}$. By (10.9), we have
\begin{equation}
P_{\phi}(\psi)(s)=\sum_{\alpha=1}^{k}\sum_{j=1}^{n}\rho_{\alpha}(\phi(s))(\psi(s),e_{\alpha,j}(\phi(s)))_{T_{\phi(s)}N}\otimes e_{\alpha,j}(\phi(s))
\end{equation}
for any $s\in S^{1}$.

Note that (10.10) is written as the sum of the form $m(\phi)\psi$, where $m(\,\cdot\,)$ is a smooth function on $N$. Thus, to prove (10.8), it suffices to prove the inequality for the multiplication operator of the form $H^{1/2}(S^{1})\ni\psi\mapsto m(\phi)\psi$, where $m$ is a smooth function on $N$. This follows form 10.2 as follows: Note that $m(\phi)\in H^{1}(S^{1})$ and $\|m(\phi)\|_{H^{1}(S^{1})}\le\|m\|_{C^{1}}(1+\|\partial_{s}\phi\|_{L^{2}(S^{1})})$. Thus by (10.2), we have $m(\phi)\psi\in H^{1/2}(S^{1})$ and
\begin{align*}
\|m(\phi)\psi\|_{H^{1/2}(S^{1})}&\le C(\|m(\phi)\|_{L^{\infty}(S^{1})}\|\psi\|_{H^{1/2}(S^{1})}+\|\psi\|_{L^{4}(S^{1})}\|m(\phi)\|_{H^{1}(S^{1})}^{\frac{1}{2}}\|m(\phi)\|_{L^{\infty}(S^{1})}^{\frac{1}{2}})\\
&\le C(\|\psi\|_{H^{1/2}(S^{1})}+\|\partial_{s}\phi\|_{L^{2}(S^{1})}^{\frac{1}{2}}\|\psi\|_{L^{4}(S^{1})}).
\end{align*}
This completes the proof.
\hfill$\Box$

\end{document}